\documentclass[11pt,reqno]{amsart}

\usepackage{amsfonts}
\usepackage{eurosym}
\usepackage{amssymb}
\usepackage{amsthm}
\usepackage{amsmath}
\usepackage{bm}
\usepackage{cite}
\usepackage{mathrsfs}
\usepackage{xcolor}
\usepackage[OT1]{fontenc}
\usepackage[left=2cm, right=2cm, top=3cm]{geometry}
\usepackage{hyperref}
\hypersetup{colorlinks=true, linkcolor=blue, citecolor=red}

\usepackage{times}
\usepackage[shortlabels]{enumitem}
\usepackage{enumitem}

\numberwithin{equation}{section}
\usepackage[square,numbers,sectionbib]{natbib}
\usepackage{varioref}

\usepackage{listings}
\usepackage{color}
\definecolor{dkgreen}{rgb}{0,0.6,0}
\definecolor{gray}{rgb}{0.5,0.5,0.5}
\definecolor{mauve}{rgb}{0.58,0,0.82}

\allowdisplaybreaks
\newtheorem{theorem}{Theorem}[section]

\newtheorem{definition}[theorem]{Definition}
\newtheorem{lemma}[theorem]{Lemma}

\theoremstyle{remark}
\newtheorem{remark}[theorem]{Remark}

\def\d{{\rm d}}

\def\uu{\mathbf{v}}
\def\vv{\mathbf{v}}
\def\ww{\mathbf{w}}

\def\ddt{\frac{\d}{\d t}}

\def\pphi{\boldsymbol{\varphi}}

\def\Psip{\Psi_{,\pphi}}
\def\PsipA{\Psi^1_{,\pphi}}

\newcommand{\numberset}{\mathbb}
\newcommand{\N}{\numberset{N}}
\newcommand{\R}{\numberset{R}}
\def\AA{\mathbb{A}}
\everymath{\displaystyle}
\def \no#1#2#3 {{\bf #1} (#3), #2.}
\def \eds#1#2#3 {#1, #2, #3.}
\def\RevA#1{{\color{black}#1}}
\def\RevB#1{{\color{black}#1}}
\author[H.~Abels]{Helmut Abels}\author[H.~Garcke]{Harald Garcke}
\address[H.~Abels, H.~Garcke]{Fakult\"{a}t f\"{u}r Mathematik\\
Universit\"{a}t Regensburg\\
93040 Regensburg, Germany}
\email[H.~Abels]{helmut.abels@ur.de}
\email[H.~Garcke]{harald.garcke@ur.de}
\author[A.~Poiatti]{Andrea Poiatti}
\address[A.~Poiatti]{Dipartimento di Matematica\\
Politecnico di Milano\\
Milano 20133, Italy}
\email[A.~Poiatti]{andrea.poiatti@polimi.it}

\begin{document}

\title[Incompressible Multi-Phase Flows]{Mathematical Analysis of a Diffuse Interface Model for Multi-Phase Flows of Incompressible Viscous Fluids with Different Densities}
\date{\today }
\keywords{Multi-phase flow, multi-component Cahn-Hilliard equation, weak solutions, strong solutions, convergence to equilibrium}

\subjclass[2010]{35B40 · 35B65. 35Q30 · 35Q35 · 76D03 · 76D05 ·
76D45}

%\date{\today}
\maketitle
\begin{abstract}
   We analyze a diffuse interface model for multi-phase flows of $N$ incompressible, viscous Newtonian fluids with different densities. In the case of a bounded and sufficiently smooth domain existence of weak solutions in two and three space dimensions and a singular free energy density is shown. Moreover, in two space dimensions global \RevA{existence} for sufficiently regular initial data is proven. In three space dimension, \RevA{existence of strong solutions} locally in time is shown as well as regularization for large times in the absence of exterior forces. 
   % \RevA{Local well-posedness of strong solutions is also shown under the assumption of strictly separated initial data.}
   \RevA{Moreover, in both two and three dimensions,} convergence to stationary solutions as time tends to infinity is proved.
\end{abstract}
\section{Introduction}

In the present contribution we analyze a diffuse interface model for the multi-phase flow of $N$ incompressible, viscous Newtonian fluids with different densities ($N\geq 2$). In such a diffuse interface model a partial mixing of the macroscopically immiscible fluids is taken into account. This has the advantage that e.g.\ topological singularities of the (diffuse) interface can be described easily. More precisely, we consider the the following inhomogeneous, multi-component Navier-Stokes/Cahn-Hilliard system, which was derived by Dong~\cite{DongJCP2014} with minor modifications:
\begin{alignat}{2}\label{eq:NSCH1}
\rho\partial_t\mathbf{v}+((\rho\mathbf{v}+\mathbf{J}_\rho)\cdot\nabla)\mathbf{v}-\operatorname{div}(2\nu(\pphi)D\mathbf{v})+\nabla p&=-\operatorname{div}(\nabla\pphi^T\Gamma\nabla\pphi), &&\quad\text{ in }\Omega\times(0,T),\\ \label{eq:NSCH2}
 \operatorname{div}\mathbf{v}&=0,&&\quad\text{ in }\Omega\times(0,T),\\ \label{eq:NSCH3}
 \partial_t\pphi+(\nabla\pphi)\mathbf{v}&=\text{div}(\mathbf{M}(\pphi)\nabla\ww)
 ,&&\quad\text{ in }\Omega\times(0,T),\\ \label{eq:NSCH4}
 \ww&= -\boldsymbol\Gamma\Delta\pphi+\mathbf{P}(\Psi_{,\pphi}^1(\pphi)-\mathbf{A}\pphi),&&\quad\text{ in }\Omega\times(0,T),\\
 \mathbf{J}_\rho&=-\nabla\ww^T{\mathbf{M}}(\pphi)\widetilde{\boldsymbol\rho},&&\quad\text{ in }\Omega\times(0,T),\label{eq:NSCH4bis}
 \end{alignat}
 together with boundary and initial conditions
 \begin{alignat}{2}
     \label{eq:NSCH5}
 \mathbf{v}&=\mathbf{0}, (\boldsymbol\Gamma\nabla\pphi) \mathbf{n}=\mathbf{0},&&\quad\text{ on }\partial\Omega\times(0,T),\\ \label{eq:NSCH6}
 \partial_\mathbf{n}\ww&=(\nabla\ww)\mathbf{n}=\mathbf{0},&&\quad\text{ on   }\partial\Omega\times(0,T),\\\label{eq:NSCH7} 
 (\mathbf{v},\pphi)|_{t=0}&=(\mathbf{v}_0,\pphi_0),&&\quad\text{ in }\Omega.
 \end{alignat} 
Here $\Omega\subseteq \R^n$, $n=2,3$, is a bounded domain with $C^4$-boundary, $0<T<\infty$, $\mathbf{v}\colon \overline{\Omega}\times [0,T)\to \R^n$ is the mean velocity of the fluid mixture, $\pphi\colon \overline{\Omega}\times [0,T)\to \R^N$ is the vector of volume fractions of the $N$ fluids/components, i.e., $\varphi_j(x,t)$ is the volume fraction of fluid $j$ at a point $x$ in space and time $t$, $\mathbf{w}\colon \overline{\Omega}\times [0,T)\to \R^N$ is the vector of chemical potentials for the $N$ fluids and $\rho=\sum_{j=1}^N\tilde{\rho}_j \varphi_j$ is density of the mixture, where $\tilde{\rho}_j>0$ is the (constant) specific density of fluid $j$. We also have set $\widetilde{\boldsymbol{\rho}}:=(\widetilde{\rho}_1,\ldots,\widetilde{\rho}_N)^T$. Note that $\sum_{j=1}^N \varphi_j(x,t)=1$ for all $x\in\overline{\Omega}$, $t\in [0,T)$. Moreover, $\nu(\pphi)$ describes the viscosity of the fluid mixture with volume fraction $\pphi$, which is assumed to be bounded and strictly positive, \RevA{$D\vv=\tfrac{\nabla\vv+(\nabla\vv)^T}{2}$ represents the symmetric strain rate tensor}, $\mathbf{M}(\pphi)\in \R^{N\times N}$ is a suitable mobility matrix (depending on $\pphi$), $\boldsymbol\Gamma\in\R^{N\times N}$ is a positive definite matrix and $\RevA{\mathbf A}\in \R^{N\times N}$. Finally, $\Psi^1$ is a homogeneous free energy density, which will be specified later and $\RevA{\mathbf P}\colon \R^N\to \R^N$ is the projection on the orthogonal complement of $(1,\ldots,1)^T$.

In the literature mostly the case of two fluids, i.e., $N=2$, is considered. In this case the model \eqref{eq:NSCH1}-\eqref{eq:NSCH7} reduces to the model derived by the first two authors and Gr\"un~\cite{AGG} and, if additionally $\tilde{\rho}_1=\tilde{\rho}_2$ (matched densities), to the well-known "Model H", cf.~Hohenberg and Halperin~\cite{HH77} and Gurtin et al.~\cite{GPV96}. Since $\varphi_1(x,t)+\varphi_2(x,t)=1$ in this case, it is sufficient to consider $\varphi(x,t)=\varphi_2(x,t)-\varphi_1(x,t)$ and the multi-component convective Cahn-Hilliard system \eqref{eq:NSCH3}-\eqref{eq:NSCH4} reduces to a scalar convective Cahn-Hilliard equation. We refer to Abels, Garcke and Giorgini~\cite{AGGio} for recent results on well-posedness and long-time behavior in case of two components and further references. An alternative class of diffuse interface models is given by (conserved) Navier-Stokes/Allen-Cahn type systems. For these models we refer to Giorgini, Grasselli, and Wu~\cite{GGW} for recents results and further references.  Moreover, there is also a variant of \eqref{eq:NSCH1}-\eqref{eq:NSCH4}, where the (convective) Cahn-Hilliard part is replaced by its ``non-local variant'', cf.\ Gal et al.~\cite{GGGP} for the details, recent results and further references. 

%%% Multi-component NSCH -- numerics, few analtic results
So far diffuse interface models for multi-phase flows with more than two components/fluids have been studied mainly numerically and from the modeling point of view. A first diffuse interface model for the case of three fluids with same densities was presented by Kim, Kang and Lowengrub in \cite{KimKangLowengrub}, where also a conservative multigrid method was discussed. A derivation of a thermodynamiscally consistent diffuse interface model for a multi-phase flows with different densities was given by Kim and Lowengrub~\cite{KimLowengrub}. We note that the later model is different from \eqref{eq:NSCH1}-\eqref{eq:NSCH4} since it is based on a barycentric/mass-averaged mean velocity, which does not yield ``$\operatorname{div} \mathbf{v}=0$'', while \eqref{eq:NSCH1}-\eqref{eq:NSCH4} is based on a volume averaged mean velocity. The model in \cite{KimLowengrub} is the multi-component counterpart of the model derived by Lowengrub and Truskinovski~\cite{LowengrubQuasiIncompressible} in the case of two fluids. We refer to ten Eikelder et al.~\cite{UnifiedNSCH} for a unified approach and comparision of different kinds of diffuse interface models for two-phase flows with different densities. The present multi-component Navier-Stokes/Cahn-Hilliard system is a minor modification of the model derived by Dong in \cite{DongJCP2014}. It can be considered as a multi-component extension of the model in~\cite{AGG}.   
 We refer to Yang and Dong~\cite{YangDongJCP2018} as well as Dong~\cite{DongJCP2018} for numerical algorithms for this model and further references. In the case of three components a non-smooth variant of the model was presented in Banas and N\"urnberg~\cite{Banasnurnberg}. They presented a fully discrete finite element approximation and proved convergence to a weak solution of the continuous model in the case of matched densites, proving in particular existence of weak solutions.  An extension of the present model \eqref{eq:NSCH1}-\eqref{eq:NSCH4} to include surfactants was derived by Dunbar, Lam and Stinner~\cite{DunbarLamStinner}, where also its relation to sharp interface models is discussed. 
\RevB{ In fact, a sharp interface limit of the system \eqref{eq:NSCH1}-\eqref{eq:NSCH4bis}, in the case where the mobilities in the pure phases vanish to leading order in the interfacial thickness, is derived as a specific case using formally matched asymptotic expansions in \cite{DunbarLamStinner}.
 The sharp interface limit is stated in Section 2.7 of \cite{DunbarLamStinner} and by neglecting the influence of surfactants, one obtains the sharp interface limit
 of \eqref{eq:NSCH1}-\eqref{eq:NSCH4bis}. As a new feature in the multi-component case, one gets an angle condition at triple junctions.}

To the best of our knowledge there is no systematic study of the Navier-Stokes/Cahn-Hilliard system \eqref{eq:NSCH1}-\eqref{eq:NSCH4bis} (or similar models) in the case of more than two components -- not even in the case of matched densities, i.e., $\tilde{\rho}_1= \ldots = \tilde{\rho}_N$. It is the aim of this contribution to fill this gap and extend the results on existence of weak solutions by the first two authors and Depner~\cite{ADG} as well as the regularity, well-posedness result recently obtained in \cite{AGGio} in the case of two-component to the present multi-component model. Since $\pphi$ and $\ww$ are vector-valued (instead of scalar) functions, a careful adaption of the arguments in the scalar case and the coupling between the components is needed. Moreover, much fewer results are already available in the multi-component case.  Recent analytic results for the multi-component Cahn-Hilliard equation and further references can be found in Gal et al.~\cite{GGPS}.

The structure of the manuscript is as follows: In Section~\ref{sec:Framework} we summarize some basic notation, assumptions and results. Afterwards several auxiliary results and basic defintions are presented in Section~\ref{sec:Preliminaries}. The main results on existence of weak and strong solutions as well as regularity and asymptotic behavior of weak solutions are stated in Section~\ref{sec:main}. Section~\ref{timediscr} is devoted to solvability of a semi-implicit time discretization of \eqref{eq:NSCH1}-\eqref{eq:NSCH4} with an additional regularization term in \eqref{eq:NSCH1}, which is a discretization of $-\alpha\Delta\partial_t\vv$, $\alpha\geq 0$. This regularization term is needed to prove the existence of strong solutions. In Section~\ref{sec:Weak} the main result on existence of weak solutions in two and three space dimensions is proved by passing to the limit in the time discretization \RevA{without} regularization term, i.e., $\alpha=0$. In order to prove existence of strong solutions, globally in time in two space dimensions and locally in time in three space dimensions, we first choose $\alpha>0$ and pass to the limit in the time discretized system (with additional regularization term) to obtain existence of (weak) solutions of \eqref{eq:NSCH1}-\eqref{eq:NSCH4} with an additional regularization term $-\alpha\Delta\partial_t\vv$ in \eqref{eq:NSCH1}. This regularization term is needed to justify higher-order estimates (globally in time in two space dimensions and locally in time in three space dimensions) to finally obtain a strong solution in the limit as $\alpha\to 0$. We note that we were not able to get the corresponding higher order-estimates for the time discrete system. An essential point for the regularity study of weak solutions is a novel regularity result for the multi-component, convective Cahn-Hilliard equation, which is proved in Section~\ref{sec:CH}. Finally, in Section~\ref{sec:long} it is shown that any weak solution becomes regular for large times with concentrations strictly positive and converges to a stationary solution of \eqref{eq:NSCH1}-\eqref{eq:NSCH4}, which consists of $\vv = \mathbf{0}$ and $\pphi$ a stationary solution of the multi-component Cahn-Hilliard equation.

% \emph{To do:
% \begin{itemize}
% \item NSAC: Qing Xia, Junxiang Yang, Yibao Li~\cite{XYL} (conservative Navier-Stokes/Allen-Cahn system, I do not see the relevance)
% \end{itemize}
% }

 \section{The mathematical framework}\label{sec:Framework}
 
 \label{setting} Let $\Omega\subset \R^n$, $n=2,3$, be a bounded domain with boundary of class $C^4$. The Sobolev spaces are denoted as usual by $W^{k,p}(\Omega )$%
 , where $k\in \mathbb{N}$ and $1\leq p\leq \infty $, with norm $\Vert \cdot
 \Vert _{W^{k,p}(\Omega )}$. The Hilbert space $W^{k,2}(\Omega )$ is denoted
 by $H^{k}(\Omega )$ with norm $\Vert \cdot \Vert _{H^{k}(\Omega )}$.
 Moreover, given a (real) vector space $X$, we denote by $\mathbf{X}$ the generic space of vectors or matrices, with each component belonging to $X$. In this case $\vert \mathbf{v} \vert$ is the Euclidean norm
 of $\mathbf{v}\in  \mathbf{X}$, i.e., $\vert\mathbf{v}\vert^2=\sum_{j}\|v_j\|_X^2$.
 We then denote by $(\cdot,\cdot )_{\Omega}$ the inner product in $\mathbf{L}^{2}(\Omega )$ and by $\Vert \cdot \Vert $
 the induced norm. We indicate by $(\cdot ,\cdot )_{H}$ and $\Vert \cdot
 \Vert _{H}$ the canonical inner product and its induced norm in the (real) Hilbert
 space $H$, respectively. We also define the spatial average of an integrable function
 $f:\Omega \to\mathbb{R}$ as
 \begin{equation*}
 \overline{f}:=|\Omega |_n^{-1}\int_{\Omega }f(x)dx,
 \end{equation*}%
 where $|\Omega |_n$ stands for the $n$-dimensional Lebesgue measure of
 $\Omega $. Further, we introduce the
 affine hyperplane
 \begin{equation}
 \Sigma :=\left\{ \mathbf{c}^{\prime }\in \mathbb{R}^{N}:%
 \sum_{i=1}^{N}c_{i}^{\prime }=1\right\} ,  \label{sigma}
 \end{equation}%
 and since only the nonnegative values for the $\varphi_i$ are
 physically relevant, we also define the Gibbs simplex
 \begin{equation}
 \mathbf{G}:=\left\{ \mathbf{c}^{\prime }\in \mathbb{R}^{N}:%
 \sum_{i=1}^{N}c_{i}^{\prime }=1,\quad c_{i}^{\prime }\geq 0,\quad i=1,\ldots
 ,N\right\} ,  \label{Gibbs}
 \end{equation}%
 and the tangent space to the affine hyperplane $\Sigma $
 \begin{equation}
 T\Sigma :=\left\{ \mathbf{d}^{\prime }\in \mathbb{R}^{N}:%
 \sum_{i=1}^{N}d_{i}^{\prime }=0\right\} .  \label{sigma2}
 \end{equation}%
 We introduce the following useful notation:
 \begin{align*}
 &\mathbf{H}_{0}:=\left\{\mathbf{f}\in \mathbf{L}^{2}(\Omega ):\ \int_{\Omega
 }\mathbf{f}\ dx=0\text{ and } \mathbf{f}(x)\in
 T\Sigma \;\text{ for a.a. }x\in \Omega \right\},
\\&\widetilde{\mathbf{H}}_{0}:=\left\{\mathbf{f}\in \mathbf{L}^{2}(\Omega ):\mathbf{f}(x)\in T\Sigma \,\text{ for a.a. }x\in
 \Omega \right\},
\\&\mathbf{V}_{0}:=\left\{\mathbf{f}\in \mathbf{H}^{1}(\Omega ):\int_{\Omega }%
 \mathbf{f}\ dx=0\text{ and }\mathbf{f}(x)\in T\Sigma\, \text{ for a.a. }x\in
 \Omega \right\},
\\&\widetilde{\mathbf{V}}_{0}:=\left\{\mathbf{f}\in \mathbf{H}^{1}(\Omega ):\ \mathbf{f}%
 (x)\in T\Sigma \, \text{ for a.a. }x\in \Omega \right\}.
 \end{align*}
 
 Notice that the spaces above are still Hilbert spaces with the same
 inner products defined in $\mathbf{L}^2(\Omega)$ and in $\mathbf{H}^1(\Omega)$, respectively. Furthermore we have the Hilbert triplets $\mathbf{V}_0\hookrightarrow\hookrightarrow \mathbf{H}_0\hookrightarrow \mathbf{V}_0'$ and $\widetilde{\mathbf{V}}_0\hookrightarrow\hookrightarrow \widetilde{\mathbf{H}}_0\hookrightarrow \widetilde{\mathbf{V}}_0'$. We then introduce the Euclidean projection $\mathbf{P}$ of $\mathbb{R}^{N}$
 onto $T\Sigma $, by setting, for $l=1,\ldots ,N$,
 \begin{equation*}
 \left( \mathbf{P}\mathbf{v}\right) _{l}=\left( \mathbf{v}-\left( \frac{1}{N}%
 \sum_{i=1}^{N}v_{i}\right) \mathbf{e}\right) _{l}=\frac{1}{N}%
 \sum_{m=1}^{N}(v_{l}-v_{m}),\quad \forall \mathbf{v}\in \R^N,
 \end{equation*}%
 where $\mathbf{e}=(1,\ldots,1)^T$.
 Notice that the projector $\mathbf{P}$ is also an orthogonal $\mathbf{L}%
 ^{2}(\Omega )$-projector on $\widetilde{\mathbf{H}}_{0}$, being symmetric and idempotent.
 
 For the velocity field we set
 \begin{equation*}
 \mathbf{H}_\sigma=\overline{\{\mathbf{u}\in\ \mathbf{C}^\infty_0(\Omega):\ \operatorname{div}\ \mathbf{u}=0\}}^{\mathbf{L}^2(\Omega)}\ \ \ \ \ 	\mathbf{V}_\sigma=\overline{\{\mathbf{u}\in\ \mathbf{C}^\infty_0(\Omega):\ \operatorname{div}\ \mathbf{u}=0\}}^{\mathbf{H}^1(\Omega)}\quad
 \mathbf{W}_\sigma={\mathbf{H}}^2(\Omega)\cap\mathbf{V}_\sigma.
 \end{equation*}
 We set $\mathbb{P}$ as the Leray projector from $\mathbf{L}^2(\Omega)$ to $\mathbf{H}_\sigma$ and we also introduce $\mathbb{A}=-\mathbb{P}\Delta: \mathfrak{D}(\mathbb{A})=\mathbf{W}_\sigma\subset \mathbf{H}_\sigma\to \mathbf{H}_\sigma$ as the standard Stokes operator. Notice that $\mathfrak{D}(\mathbb{A})=\mathbf{W}_\sigma$.
 Now we consider the following problem: given $\mathbf f\in \mathbf{H}_\sigma$, find $\uu\in \mathbf{W}_\sigma$ such that 
 $$
 -\Delta\uu+\nabla p= \mathbf f.
 $$
 By standard theory there exists a unique $\uu\in \mathbf{W}_\sigma$ and $p\in H^1(\Omega)$ with $\overline{p}=0$, such that 
 \begin{align}
 \Vert \uu\Vert_{\mathbf{W}_\sigma}+\Vert p\Vert_{H^1(\Omega)}\leq C\Vert \mathbf f\Vert.
 \label{reg}
 \end{align}
 Moreover, thanks to \cite[Lemma 3.1, Remark 3.2]{GGGP} we also have 
 \begin{align}
 \Vert p\Vert_{L^s(\Omega)}\leq C\Vert \nabla \AA^{-1}\mathbf f\Vert^\frac 1 2 \Vert \mathbf f\Vert^\frac 1 2, 
 \label{pr}
 \end{align}
 where $s=4$ if $n=2$ and $s=3$ if $n=3$.
 
 Given a measurable set $\mathcal{M} \subset \R^n$, we let $L^q(\mathcal{M}),\ 1 \leq q \leq \infty$ denote the usual
 Lebesgue-space and $\Vert\cdot\Vert_q$ its norm. Let $X$ be a Banach space. We also denote by $L^q (\mathcal{M}; X )$ the Bochner space of $X$-valued $q$-integrable (or essentially bounded functions).
When $M = (a, b)$, we simply write $L^q (a, b; X )$ and $L^q (a, b)$.
Moreover, $f \in L^q_{loc} ([0, \infty); X )$ if and only if $f \in L^q (0, T ; X )$ for every $T > 0$.
 Furthermore, $L^q_{uloc} ([0, \infty); X )$ consists of all measurable $f : [0, \infty) \rightarrow X$ such that
 \begin{align*}
 \Vert f\Vert_{L^q_{uloc} ([0, \infty); X )}=\sup_{t\geq 0}\Vert f \Vert_{L^q(t,t+1;X)}<+\infty.
 \end{align*}
 In the case of a bounded interval $I=[0,T)$, we set $L^q_{uloc} (I; X ):=L^q (0,T; X) $, whereas in general we only have $L^q (0, \infty; X )\hookrightarrow L^q_{uloc} ([0, \infty); X )$. 
  For any interval $I=[0,T]$, for $T>0$, or $I=[0,\infty)$, 
 we denote by $BC(I;X)$ the Banach space of bounded continuous functions on $I$, equipped with the supremum norm. The space $BUC(I;X)$ is then its subspace of bounded and uniformly continuous functions. Moreover, we set $BC_w(I;X)$ to be the topological vector space of bounded weakly continuous functions $f:I\to X$. The function space 
 $C_0^\infty (0,T ; X)$ denotes the vector space of all $C^\infty$-functions $f : [0,T]\to X$ with compact support in $(0,T)$. In conclusion, $W^{1,p} (0, T ; X)$, $1 \leq p < \infty$, is the space of functions $f$ such that $\ddt f\in L^p(0,T;X)$ and $f\in L^p(0,T;X)$, where $\ddt$ denotes the vector-valued distributional derivative of $f$. Furthermore, replacing $L^p (0, T ; X)$ by $L^p_{uloc} (I; X)$ we define $W^{1,p}_{uloc} (I; X)$. We then set $H^1 (0, T ; X) =W^{1,2}(0,T;X)$, and $H^1_{uloc} (I ; X) =W^{1,2}_{uloc}(I;X)$.

% Next, we define the set
 %\begin{equation}
 %\mathcal{K}:=\left\{ \boldsymbol{\eta }\in \mathbf{H}^{1}(\Omega
 %):\sum_{i=1}^{N}\eta _{i}=1,\quad\eta _{i}\geq 0,\quad \forall i=1,\ldots
 %,N\right\} .  \label{K}
 %\end{equation}%
 For the sake of simplicity we will adopt the compact notation $\mathbf{v}%
 \geq k$, with $\mathbf{v}\in \mathbb{R}^{N}$ and $k\in \mathbb{R}$ to
 indicate the relation $v_{i}\geq k$ for any $i=1,\ldots ,N$, as well as we will write $\mathbf{z}=\mathbf{v}+k$ to indicate that $z_i=v_i+k$, $i=1,\ldots,N$. Moreover, given a matrix $\mathbf{A}\in \R^{N\times N}$ and a vector-valued function $\mathbf{g}\in \R^N$, we set
  $$(\mathbf{A}\nabla \mathbf{g})_{ij}:=\sum_{k=1}^N\mathbf{A}_{ik}{g}_{j,k}, \quad\text{for }i=1,\ldots, N,\quad j=1,\ldots, n,$$
 where ${g}_{j,k}$ is the $k$-th partial derivative of $g_j$.
 Concerning the free energy function, we recall that the
 Boltzmann-Gibbs entropy potential is given by 
 $$
 \Psi^1(\mathbf{u}):=\sum_{i=1}^N\theta u_i\ln u_i=\sum_{i=1}^N\psi(u_i),
 $$
where $\theta>0$ is the absolute temperature and the Boltzmann constant is assumed for simplicity equal to $1$. Thus
 \begin{equation*}
 \left(\Psip^1(\mathbf{u})\right) _{i}=\psi' (u_{i}):=\theta (\operatorname{ln} u_{i}+1),\quad \text{for }i=1,\ldots ,N.
 \end{equation*}%
 In order to include a large admissible class of entropy functionals, we assume
 
 \begin{itemize}
 	\item[(\textbf{E0})]  \label{E0}
 	$\psi \in C\left([0,1]\right) \cap C^{2}\left((0,1]\right)$  and $\psi ^{\prime \prime }(s)\geq \zeta >0,$ for all $s\in
 	(0,1];$
 	
 	\item[(\textbf{E1})] $\lim_{s\rightarrow 0^{+}}\psi^{\prime }\left( s\right)
 	=-\infty ;$
 \item[(\textbf{E2})] 
 \label{E2}
 $\psi$ is (real) analytic in $(0,1)$.
 \end{itemize}
 We also extend $\psi(s)=+\infty ,$ for any $s\in (-\infty ,0)$, and extend $\psi$ for all $s\in \lbrack 1,\infty )$ so that $\psi$ is still a $C^2$ function on $(0,+\infty)$. %Notice that
% assumption (\textbf{E1}) implies that there exists $\xi \in (0,1)$ such that
 %$\psi ^{\prime }(\xi )=0$. 
% Furthermore, there is no loss of generality in
 %assuming $\psi (\xi)=0$ since the potential may be defined up to a constant.
 Furthermore, there is no loss of generality in
assuming $\psi (0)=0$ since the potential may be defined up to a constant. These assumptions on $\psi$ are standard and have been adopted also in \cite{GGPS, GP}. 
% This also entails that $\psi (s)\geq 0,$ for all $s\in \lbrack 0,1]$.
 We also make the following assumptions
\begin{enumerate}
	\item[(\textbf{M0})] \label{M0} The mobility matrix $\mathbf{M}\colon \R^N\to\R^{N\times N}$ is symmetric, positive semidefinite and continuously differentiable. Moreover, there exists $C_\mathbf{M}$ such that $\vert\mathbf{M}(\mathbf{s})\vert\leq C_\mathbf{M}$ for any $\mathbf{s}\in \R^N$.  We also assume $\mathbf{M}$ to be nondegenerate, i.e., there exists a uniform constant $C_0>0$ such that
 \begin{align}
 \boldsymbol\xi^T\mathbf{M}(\mathbf{s})\boldsymbol\xi\geq C_0\Vert\boldsymbol\xi\Vert^2,\quad \forall \mathbf{s}\in \R^N,\forall \boldsymbol\xi\in T\Sigma.
 \label{nondeg}
 \end{align}
 Moreover we assume 
 \begin{align}
     \label{m}
 \sum_{j=1}^N \mathbf{M}_{ij}=0,\quad \text{ for any }i=1,\ldots,N.
 \end{align}
\RevA{This condition is standard for multi-component models (see, for instance, \cite{EL, GarckeElasticMisfit}) and it is the most natural assumption to be consistent with the constraint $\sum_{i=1}^N\varphi_i\equiv 1$. It corresponds to assuming that the sum of all the thermodynamic fluxes, defined, for $k=1,\ldots,N$, by $\mathbf J_k=-\sum_{l=1}^N\mathbf M_{kl}\nabla w_l$, is equal to zero.}	\item[(\textbf{M1})] 
\label{M1}The viscosity $\nu:\R^N\to \R$ is such that $\nu\in W^{1,\infty}(\R^N)$ and there exist $\nu_*,\nu^*>0$ such that
	\begin{align}
	0<\nu_*\leq \nu(\mathbf{s})\leq \nu^*, \quad \forall \mathbf{s}\in \R^N.
	\label{nu1}
	\end{align} 
	These conditions also imply that there exists $C_{\nu}>0$ such that $	\nu_{,\pphi}(\mathbf{s}):=\nabla_\mathbf{s}\nu(\mathbf{s})$ satisfies
	$$
	\vert \nu_{,\pphi}(\mathbf{s})\vert\leq C_\nu \quad \forall \mathbf{s}\in \R^N.
	$$
\end{enumerate}
 \begin{remark}
 	As the solution $\pphi$ will lie on the Gibbs simplex $\mathbf{G}$, we only need the functions $\mathbf{M}$ and $\nu$ on $G$. We then extend them to the whole of $\R^N$ such that Assumption (\textbf{M0}) is fulfilled.
 \end{remark}
In the sequel we will assume the matrix $\RevA{\boldsymbol\Gamma}$ (see \eqref{eq:NSCH1}) to be the $N\times N$ identity matrix $\mathbf{I}_N$. Actually, all the results still hold with matrices of the form $\RevA{\boldsymbol\Gamma}=\zeta \mathbf{I}_N$ with $\zeta>0$.
 
 \section{Preliminaries}\label{sec:Preliminaries}
 We now present some preliminary results concerning the following homogeneous Neumann elliptic problem with a logarithmic
 convex nonlinear term: Given
 $
 \mathbf{f}\in
 \widetilde{\mathbf{H}}_0,
 $ find $\pphi\in \mathbf{H}^2(\Omega)$ with $\Psi^1_{,\pphi}(\pphi)\in \mathbf{L}^2(\Omega)$ such that
 \begin{alignat}{2}
 \nonumber- \Delta\pphi+\mathbf{P}\Psi^1_{,\pphi}(\pphi)&=\mathbf f,&&\quad
 \text{ a.e. in }\Omega,\\
 \partial_\mathbf{n}\pphi&=0,&&\quad \text{ a.e. on }\partial\Omega,  \label{steady}\\
 \sum_{i=1}^N\varphi_i &= 1, &&\quad \text{ in } \Omega.\nonumber
 \end{alignat}
 
 In particular, we have  
 \begin{theorem}
 	\label{steaddy}
 	Let $\Omega\subset \R^n$ be a sufficiently smooth bounded domain. Assume $
 	\mathbf{f}\in
 	\widetilde{\mathbf{H}}_0$ and assumptions (\textbf{E0})-(\textbf{E1}).
 	\begin{enumerate}
 	\item Problem \eqref{steady} has a unique solution $\pphi$ such that  
 	\begin{align}
 	0<\pphi<1,\quad\text{ a.e. in }\Omega,
 	\label{quasi_sep}
 	\end{align}
 	and 
  \begin{align}
    \Vert \pphi\Vert_{\mathbf{H}^2(\Omega)}+\Vert \PsipA(\pphi)\Vert\leq C(1+\Vert \mathbf f\Vert).\label{H2b}
  \end{align}
  Moreover it also holds the following: for any sequence $(\mathbf{f}_k)_{k\in\N}\in \widetilde{\mathbf{H}}_0$, such that $\mathbf f_k\to \mathbf f$ as $n\to \infty$, we have
 	\begin{align}
 	\pphi_k\to \pphi\quad \text{ in   }\ \mathbf{H}^{2-s}(\Omega),\quad \text{for any } s\in\left(0,\frac 1 4\right),
 	\label{conv}
 	\end{align}
 	where $\pphi_k$ are the solutions to \eqref{steady} corresponding to $\mathbf f_k$, and $\pphi$ corresponds to $\mathbf f$.
 	
 	\item If we additionally assume that $\mathbf{f}\in \mathbf{L}^\infty(\Omega)$, then the unique solution $\pphi$ is also strictly
 	separated from the pure phases, i.e., there exists $0<{\delta }=\delta (\Vert \mathbf f\Vert_{\mathbf{L}^\infty(\Omega)})<\frac{1}{N}$, such that
 	\begin{equation}
 	\delta<\pphi(x)<1-(N-1)\delta  \label{prop}
 	\end{equation}%
 	for any $x\in \overline{\Omega }$. 
 	 
 	% \item If we assume, additionally to point (1), that
 	% \begin{align}
 	%  	\psi^{\prime\prime}(s)\leq e^{C\vert\psi^\prime(s)\vert+C},\qquad \forall s\in(0,1),
 	% \label{pot}
 	% \end{align}
 	% for some $C>0$ and $\mathbf f\in \widetilde{\mathbf{H}}_0\cap \mathbf{W}^{1,n}(\Omega)$, then also there exists
 	% $0<{\delta }=\delta (\Vert \mathbf f\Vert_{\mathbf{W}^{1,n}(\Omega)}^n,\Vert \pphi\Vert_{\mathbf{W}^{1,6}(\Omega)})<\frac{1}{N}$ such that the same \eqref{prop} holds. 
	
 \item 	Assume that the couple $(\mathbf f,\mathbf{m})\in \widetilde{\mathbf{V}}_0\times \Sigma$, with $\mathbf{m}\in(0,1)^N\cap \Sigma$, is such that there exists $\pphi\in \mathbf{G}$ almost everywhere, $\pphi\in \mathbf{H}^2(\Omega)$, $\overline{\pphi}=\mathbf{m}$, which solves \eqref{steady}.
 Then there exists $C=C(\min_{i=1,\ldots,N}{m}_i)>0$ such that 
 	 	\begin{align}
 	\Vert \pphi\Vert_{\mathbf{W}^{2,p}(\Omega)}+\Vert\Psi^1_{,\pphi}(\pphi)\Vert_{\mathbf{L}^p(\Omega)}+\left\vert\overline{\mathbf f}\right\vert\leq C(\min_{i=1,\ldots,N}{m}_i)(1+\Vert \nabla\mathbf f\Vert),
 		\label{g}
 	\end{align}
 for any $p\in[2,\infty)$ if $n=2$, $p\in[2,6]$ in $n=3$.
 	\end{enumerate}
 \end{theorem}
\begin{remark}
	\label{control} Notice that the assumption $\mathbf{m}\in(0,1)^N$, with $\mathbf{m}\in \Sigma$, of point (3) of the above theorem implies that there exists $\zeta>0$ such that $\zeta<\overline{\varphi}_{i}<1-\zeta$
	for any $i=1,\ldots,N$. Indeed, we have for any $i=1,\ldots,N$,
	\begin{equation*}
	0<\min_{j=1,\ldots,N}{m_j}\leq \overline{\varphi}_{i}=1-\sum_{j\not= i}{\overline{\varphi}}_{j}\leq 1-(N-1)\min_{j=1,\ldots,N}{m_{j}},
	\end{equation*}
	and thus we can choose, e.g.\ $\zeta=\min_{j=1,\ldots,N}{m_j}$ since $N\geq 2$.
\end{remark}
% \begin{remark}
% 	Notice that in case assumption \eqref{pot} holds we need simply $\mathbf f \in \widetilde{\mathbf{H}}_0\cap \mathbf{W}^{1,n}(\Omega)\not \hookrightarrow \widetilde{\mathbf{H}}_0\cap \mathbf{L}^\infty(\Omega)$ to gain the separation property. Nevertheless in the following we will not assume \eqref{pot}, since it will not be necessary. This is why we postponed the proof of part (3) in the Appendix.
% \end{remark}
\begin{remark}
We point out that in (3) the \textit{a priori} assumption that $(\mathbf f,\mathbf{m})\in \widetilde{\mathbf{V}}_0\cap \mathbf{H}^1(\Omega)\times \Sigma$ is such that there exists $\pphi$ satisfying \eqref{steady} with the fixed mean value $\overline{\pphi}=\mathbf{m}\in (0,1)^N\cap \Sigma$ is essential. Indeed, having a generic couple $(\mathbf f,\mathbf{m})\in \widetilde{\mathbf{V}}_0\times \Sigma$ does not ensure the existence of a solution to \eqref{steady}, since the problem would be overdetermined (see point (1): to find a unique $\pphi$ it is enough to fix $\mathbf f$ and not $\mathbf{m}$).
\end{remark}
 \textbf{Proof of Theorem \ref{steaddy}.}
 	\begin{enumerate}
 \item[\underline{Point (1)}] For the first part of the statement we can adapt the results already obtained in \cite[Theorem 8.1]{GGPS}. In particular, even though it concerns $\mathbf f\in \mathbf{L}^\infty(\Omega)\cap \widetilde{\mathbf{H}}_0$, we clearly see that, apart from the estimates related to the separation property (i.e., from \cite[(8.11)]{GGPS}), all the other results concerning well-posedness still hold for $\mathbf f\in \widetilde{\mathbf{H}}_0$. In particular we also obtain (see \cite[(8.9)]{GGPS}) that $\Psi_{,\pphi}(\pphi)\in \mathbf{L}^2(\Omega)$ and thus also \eqref{quasi_sep} consequently holds, together with \eqref{H2b}. To obtain the continuity result \eqref{conv} it is sufficient to consider the difference $\pphi_k-\pphi$ and multiply \eqref{steady}, written for this difference, by it, and then integrate over $\Omega$, to get, after integration by parts,
 	\begin{align*}
 	&\Vert\nabla(\pphi_k-\pphi)\Vert^2+C\Vert \pphi_k-\pphi\Vert^2\\&\leq \Vert\nabla(\pphi_k-\pphi)\Vert^2+(\Psi^1_{,\pphi}(\pphi_k)-\Psi^1_{,\pphi}(\pphi),\pphi_k-\pphi)_\Omega=(\mathbf f_k-\mathbf f,\pphi_k-\pphi)_\Omega\\&\leq \frac{1}{2C}\Vert \mathbf f_k-\mathbf f\Vert^2+\frac{C}{2}\Vert \pphi_k-\pphi\Vert^2,
 	\end{align*}
 	where we used the fact that $\pphi_k-\pphi\in T\Sigma$ for almost any $x\in\Omega$ and $\psi^{\prime\prime}\geq C>0$, i.e., $\psi$ is \textit{strictly} convex. This implies
 	$$
 	\Vert\pphi_k-\pphi \Vert_{\mathbf{H}^1(\Omega)}\leq C_0\Vert \mathbf f_k-\mathbf f\Vert,
 	$$ 
 	for some $C_0>0$. Therefore, since we also know (see \cite[(8.10)]{GGPS}) that 
 	$$
 	\Vert \pphi_k\Vert_{\mathbf{H}^2(\Omega)} \leq C(1+\Vert \mathbf f_k\Vert)\leq C, 
 	$$ 
 	for $k$ sufficiently large, being $(\mathbf f_k)_{k\in\N}$ convergent in $\mathbf{L}^2(\Omega)$, by a standard interpolation inequality between $\mathbf{H}^2(\Omega)$ and $\mathbf{L}^2(\Omega)$,  we easily deduce the validity of \eqref{conv}.
 	\item[\underline{Point (2)}] If we additionally assume $\mathbf f\in \mathbf{L}^\infty(\Omega)$, we have exactly the result in the statement of \cite[Thm 8.1]{GGPS}, so that \eqref{prop} holds.
 	%\item See Appendix \ref{app} for the proof.
 %	\setcounter{enumi}{3}
	\item[\underline{Point (3)}] To prove this result, we follow the argument in the proof of \cite[Thm.3.1]{GGPS}. In particular, we set
	\begin{equation*}
	\boldsymbol{\lambda  }:=%
	\overline{\mathbf{P}(%
		\PsipA(\pphi))},\quad \mathbf f_0:=\mathbf{f}-\boldsymbol{\lambda  }.
	\end{equation*}%
	We have, for all $\boldsymbol{\eta }\in \mathbf{H}%
	^{1}(\Omega )$ and for almost all $t\in (0,T)$,
	\begin{equation}
	(\mathbf f_0+\boldsymbol{\lambda },%
	\boldsymbol{\eta })_\Omega= (\nabla \pphi,\nabla
	\boldsymbol{\eta })_\Omega+(\mathbf{P}(\PsipA(\pphi)),\boldsymbol{%
		\eta })_\Omega.  \label{e3}
	\end{equation}%
	Then we exploit the convexity of $\Psi ^{1}$: for any $%
	\mathbf{k}\in \mathbf{G}$ (see \eqref{Gibbs}), 
	we have  $%
	\mathbf{k}-\pphi\in T\Sigma $ almost everywhere, so that we find
	\begin{align}
	C& \geq \int_{\Omega }\Psi ^{1}(\mathbf{k})\geq \int_{\Omega
	}\Psi ^{1}(\pphi)+\int_{\Omega }\Psi
	_{,\pphi}^{1}(\pphi)\cdot (\mathbf{k}-%
	\pphi)  \label{psi} \\
	& =\int_{\Omega }\Psi ^{1}(\pphi)+\int_{\Omega }\mathbf{P}\PsipA(\pphi%
	)\cdot (\mathbf{k}-\pphi),  \notag
	\end{align}%
	where we used 
	\begin{equation*}
	\int_{\Omega }\Psi ^{1}(\mathbf{k})\leq \max_{\mathbf{s}\in \lbrack
		0,1]}\Psi^1 (\mathbf{s})=C.
	\end{equation*}
	Moreover, we can choose $\boldsymbol{\eta }=\mathbf{k}-\pphi%
	$ in \eqref{e3} to deduce, on account of \eqref{psi}, that%
	\begin{align*}
	C& \geq \int_{\Omega }\Psi ^{1}(\mathbf{k}) \geq \int_{\Omega }\Psi ^{1}(\pphi)+ \Vert \nabla \pphi\Vert ^{2}+(\mathbf f%
	_{ 0},\mathbf{k}-\pphi)_\Omega+(\boldsymbol{\lambda
	},\mathbf{k}-\pphi)_\Omega.
	\end{align*}%
Observe that, since $\mathbf{k}\geq0$,
	\begin{equation*}
	\int_{\Omega }\sum_{i=1}^{N}k_{i}^{2}\leq \int_{\Omega }\left(
	\sum_{i=1}^{N}k_{i}\right) ^{2}=|\Omega |_n.
	\end{equation*}%
	Then by Cauchy-Schwarz's, Young's and Poincar\'{e}'s inequalities (all
	applied to $\mathbf f_{0}$), recalling that $0\leq \pphi\leq1$ and thus $\left\vert \int_\Omega \Psi^1(\pphi)\right\vert \leq C$,
	%for psi (cdelta), ok controllo dal basso
	%and the fact that $\Psi_\varepsilon^1(\mathbf{k})$, for $\varepsilon$ sufficiently small, is bounded when $\mathbf{k}\in \mathbf{G}$, %addirittura il k sarà separato, vedi sotto, quindi di sicuro controllo lo Psi indip da vareplsilon
	\begin{align}
	& (\boldsymbol{\lambda },\mathbf{k}-\pphi_{
	})_\Omega+ \Vert \nabla \pphi\Vert ^{2}  \notag \\
	& \leq -\int_{\Omega }\Psi ^{1}(\pphi)-(\mathbf f_{0},\mathbf{k}-\pphi%
	)_\Omega +C \notag \\
	& \leq C\left( 1+(\Vert \nabla \mathbf f\Vert)
	(1+\Vert \pphi\Vert )\right) \leq C(1+\Vert \nabla
	\mathbf f\Vert). \label{e4}
	\end{align}%
	Recalling Remark \ref{control}, there exists $\zeta=\min_{i=1,\ldots,N}{m_i}\in(0,\frac 1 N)$ such that, for all $%
	i=1,\ldots ,N$,
	\begin{equation*}
	0<\zeta<\overline{\varphi}_{i}<1-\zeta,
	\end{equation*}%
	with $\overline{\pphi}=\mathbf{m}$.
	Therefore we choose, for fixed $k,l=1,\ldots ,N$, $k\not=l$,
	\begin{equation*}
	\mathbf{k}=\mathbf{m}+\zeta\operatorname{sign}(\lambda _{k}-\lambda _{l})(\boldsymbol\eta_{k}-%
	\boldsymbol\eta_{l})\in \mathbf{G}
	\end{equation*}%
	in \eqref{e4}, where $\boldsymbol\eta_{j}$ is the vector with all zeros apart from the $j-th$ component, which is 1. Thus, from %
	\eqref{e4} we get that
	\begin{equation}
	|(\lambda _{k}-\lambda _{l})(t)|\leq \frac{C%
	}{\zeta|\Omega |_n}(1+\Vert \nabla \mathbf f\Vert).
	\label{s}
	\end{equation}%
Using the identity
	\begin{equation*}
	\boldsymbol{\lambda }=\frac{1}{N}\left( \sum_{l=1}^{N}(
	\lambda _{k}-\lambda _{l})\right) _{k=1,\ldots
		,N},
	\end{equation*}%
	we find
		\begin{equation}
	\vert \overline{\mathbf f}\vert=|\boldsymbol\lambda|\leq \frac{C%
	}{\zeta|\Omega |_n}(1+\Vert \nabla \mathbf f\Vert ),
	\label{s1}
	\end{equation}%
	where $\zeta$ only depends on $\min_{i=1,\ldots,N}{m}_i$ as desired. Notice that we are in the setting of point (1) of the statement of the theorem. In particular, we need to come back to the approximating scheme in $\varepsilon$ used in the proof of \cite[Thm.\ 8.1]{GGPS}: from \cite[(8.11)]{GGPS} we clearly see that, by passing to the limit in $\varepsilon$, there exists $C>0$ such that 
	$$
	\Vert \PsipA(\pphi)\Vert_{\mathbf{L}^p(\Omega)}\leq C\Vert \mathbf f\Vert_{\mathbf{L}^p(\Omega)}
	$$ 
	for any $p\geq 2$. Therefore, by Sobolev embeddings,
		$$
	\Vert \PsipA(\pphi)\Vert_{\mathbf{L}^p(\Omega)}\leq C\Vert \mathbf f\Vert_{\mathbf{H}^1(\Omega)},
	$$ 
	for any $p\geq2$ in the case $n=2$ and $p\in[2,6]$ for $n=3$. The control on the mean value of $\mathbf f$ in \eqref{s1} allows to deduce that $\Vert \mathbf f\Vert_{\mathbf{H}^1(\Omega)}\leq C(1+\Vert\nabla\mathbf f\Vert)$, so that
		$$
	\Vert \PsipA(\pphi)\Vert_{\mathbf{L}^p(\Omega)}\leq C(1+\Vert\nabla\mathbf f\Vert).
	$$ 
	 Then by elliptic regularity we easily deduce the corresponding $\mathbf{W}^{2,p}(\Omega)$ estimate on $\pphi$. This concludes the proof of \eqref{g}.
\end{enumerate}
\qed
 
 We define the total energy to be
 \begin{align}
 E_{tot}(\rho,\mathbf{v},\pphi):=\frac 1 2\int_\Omega \rho \vert\mathbf{v}\vert^2dx +\frac 1 2 \int_\Omega \vert \nabla \pphi\vert^2 dx +\int_\Omega \Psi(\pphi)dx,
 \label{total}
 \end{align}
% {\color{red}[For the time being I concentrate only on th cae $n=3$.]}
\RevA{where $\Psi(\mathbf{s}):=\Psi^1(\mathbf s)-\frac12 \mathbf s^T\mathbf A\mathbf s$, for any $\mathbf s\in [0,1]^N$, }and we introduce the notions of weak and strong solutions:
\begin{definition}[Weak solution]
	\label{weaks}
	Let $T \in (0, \infty]$ and set either $I = [0, \infty)$ if $T = \infty$ or $I = [0, T ]$ if $T < \infty$, $\mathbf{v}_0 \in\mathbf{ L}^2_\sigma(\Omega)$
	and $\pphi_0 \in \mathbf{H}^1 (\Omega)$ with $\pphi_0\in[0,1]^N$ and $\sum_{i=1}^N\varphi_{0,i}\equiv 1$ almost everywhere in $\Omega$. Assume also $\overline{\pphi}_0\in(0,1)^N$. We call the
	triple $(\mathbf{v},\pphi,\ww)$ a weak solution to \eqref{eq:NSCH1}-\eqref{eq:NSCH7} if 
	\begin{align*}
	&\mathbf{v} \in BC_w (I; \mathbf{H}_\sigma (\Omega)) \cap L ^2 (I; \mathbf{V}_\sigma) ,\\&
	\pphi\in BC_w (I; \mathbf{H}^1 (\Omega)) \cap {L}_{uloc}^2 (I; \mathbf{W}^{2,p} (\Omega))\cap H^1(I;(\mathbf{H}^1(\Omega))') ,\\&
	\sum_{i=1}^N\varphi_{i}\equiv 1 \quad \text{and}\quad\varphi_i>0\quad \forall i=1,\ldots,N \quad \text{a.e. in }\Omega\times I,
	\\& \PsipA(\pphi) \in \mathbf{L}^2_{uloc} (I; \mathbf{L}^p (\Omega)) ,\\&
	\ww\in \mathbf{L}^2_{uloc} (I; \widetilde{\mathbf{V}}_0),\quad \nabla\ww\in L^2(I;\mathbf{L}^2(\Omega)),
	\end{align*}
	where $p\in[2,\infty)$ if $n=2$ and $p\in[2,6]$ if $n=3$, and it satisfies, having set $Q_T:=\Omega\times(0,T)$,
	\begin{alignat}{2}
	\nonumber-(\rho\mathbf{v},\partial_t\boldsymbol\psi)_{Q_T}+(\operatorname{div}(\rho\mathbf{v}\otimes\mathbf{v}),\boldsymbol\psi)_{Q_T}&+(2\nu(\pphi)D\mathbf{v},D\boldsymbol\psi)_{Q_T}\\-((\mathbf{v}\otimes \mathbf{J}_\rho),\nabla\boldsymbol\psi)_{Q_T}&=((\nabla\pphi)^T\ww,\boldsymbol\psi)_{Q_T},&& \forall \boldsymbol\psi\in[{C}_{0,\sigma}^\infty(\Omega\times(0,T))]^n,\label{eq:weakAGG1}\\
	-(\pphi,\partial_t\boldsymbol\zeta)_{Q_T}-(\mathbf{v},\pphi^T\nabla\boldsymbol\zeta)_{Q_T}&=-(\mathbf{M}(\pphi)\nabla\ww,\nabla\boldsymbol\zeta)_{Q_T},&&\forall \boldsymbol\zeta\in {C}_0^\infty((0,T);[{C}^1(\overline{\Omega})]^N),\\
	\label{we} \ww&= -\Delta\pphi+\mathbf{P}(\Psi_{,\pphi}^1(\pphi)-\mathbf{A}\pphi),&&\text{ a.e. in }Q_T,\\
 \mathbf{J}_\rho&=-\nabla\ww^T{\mathbf{M}}(\pphi)\widetilde{\boldsymbol\rho},&&\text{ a.e. in }Q_T,\\
	 (\mathbf{v},\pphi)_{\vert t=0}&=(\mathbf{v}_0,\pphi_0),&&\text{ a.e. in }\Omega.
	\end{alignat}
	Moreover, having set $Q_{(s,t)}:=\Omega\times(s,t)$,
	\begin{align}
	&\nonumber E_{tot}(\rho(t),\mathbf{v}(t),\pphi(t))+\int_{Q_{(s,t)}}2\nu(\pphi)\vert D\mathbf{v}\vert^2d(x,\tau)+\int_{Q_{(s,t)}}\mathbf{M}(\pphi)\nabla\ww:\nabla\ww d(x,\tau)\\&\leq E_{tot}(\rho(s),\mathbf{v}(s),\pphi(s)),
	\label{energ}
	\end{align}
	for all $t\in[s,\infty)$ and almost all $s\in[0,\infty)$ has to hold (including $s=0$).
	\end{definition}

\begin{remark}
Notice that, recalling the definition of $\ww$,
 \begin{align*}
 \operatorname{div}(\nabla\pphi^T\nabla\pphi)&=\nabla\pphi^T\Delta\pphi+\frac 1 2 \nabla\left[\vert\nabla\pphi\vert^2\right]\\&=-\nabla\pphi^T\ww+\nabla\pphi^T\mathbf{P}\Psi_{,\pphi}(\pphi)+\frac 1 2 \nabla\left[\vert\nabla\pphi\vert^2\right]\\&
 =-\nabla\pphi^T\ww+\nabla \left[\frac 1 2\vert\nabla\pphi\vert^2+\Psi(\pphi)\right],
 \end{align*} 
 where $\nabla\pphi^T\mathbf{P}\Psi_{,\pphi}(\pphi)=\nabla\pphi^T\Psi_{,\pphi}(\pphi)$, since $\nabla\pphi^T\mathbf{e}=\mathbf{0}$, with $\mathbf{e}=(1,\ldots,1)$, due to the sum-to-one constraint of the components of $\pphi$. This means that (up to a redefinition of the pressure), we can substitute $-\operatorname{div}(\nabla\pphi^T\nabla\pphi)$ with $\nabla\pphi^T\ww$, as done in \eqref{eq:weakAGG1} of the definition above.
\end{remark}
 
We also define a strong solution as
\begin{definition}[Strong solution]
	\label{strong}
	Let $T >0$ and set either $I = [0, \infty)$ if $T = \infty$ or $I = [0, T ]$ if $T < \infty$, $\mathbf{v}_0 \in \mathbf{V}_\sigma$
	and $\pphi_0\in \mathbf{H}^2 (\Omega)$ with $\partial_\mathbf{n}\pphi_0=0$ on $\partial\Omega$ and such that $-\Delta\pphi_0+\mathbf{P}\Psip(\pphi_0)=:\ww_0\in \widetilde{\mathbf{V}}_0$, $\pphi_0\in[0,1]^N$ and $\sum_{i=1}^N\varphi_{0,i}\equiv 1$ almost everywhere in $\Omega$. Assume also $\overline{\pphi}\in(0,1)^N$. We call
	 $(\mathbf{v},\pphi,\ww,p)$ a strong solution to \eqref{eq:NSCH1}-\eqref{eq:NSCH7} if 
	\begin{align*}
	&\mathbf{v} \in L^\infty_{uloc} (I; \mathbf{V}_\sigma(\Omega)) \cap L ^2_{uloc} (I ; \mathbf{H}^2(\Omega))\cap H^1_{uloc}(I;\mathbf{H}_\sigma(\Omega)) ,\quad p\in L^2_{uloc}(I;H^1(\Omega)),\\&
	\pphi\in BUC(I; \mathbf{H}^1 (\Omega)) \cap  L^\infty(I;\mathbf{W}^{2,q}(\Omega))\cap H^1_{uloc}(I;\mathbf{H}^1(\Omega)) ,\\&
	\sum_{i=1}^N\varphi_{i}\equiv 1,\quad \text{and}\quad\varphi_i>0\quad \forall i=1,\ldots,N  \quad\text{a.e. in }\Omega\times I,\\ & \PsipA(\pphi) \in \mathbf{L}^\infty (I; \mathbf{L}^q (\Omega)) ,\quad
	\ww\in \mathbf{L}^\infty (I; \widetilde{\mathbf{V}}_0)\cap L^2_{uloc}(I;\mathbf{H}^3(\Omega)),
	\end{align*}
	with $q\in[2,\infty)$ when $n=2$ and $q=6$ when $n=3$, and it satisfies \eqref{eq:NSCH1}-\eqref{eq:NSCH7} almost everywhere in $\Omega\times I$.
	Moreover, the energy identity
	\begin{align}
	\nonumber &E_{tot}(\rho(t),\mathbf{v}(t),\pphi(t))+\int_{Q_{(s,t)}}2\nu(\pphi)\vert D\mathbf{v}\vert^2d(x,\tau)+\int_{Q_{(s,t)}}\mathbf{M}(\pphi)\nabla\ww:\nabla\ww d(x,\tau)\\& = E_{tot}(\rho(s),\mathbf{v}(s),\pphi(s)),
	\label{energ1}
	\end{align}
	has to hold for all $t\in[s,T)$ and almost all $s\in[0,T)$ (including $s=0$).
\end{definition} 
\begin{remark}
The validity of the energy identity is an immediate consequence of the regularity of strong solutions.
\end{remark}
\section{Main results}\label{sec:main}
Our main results are then the following:
\begin{theorem}[Existence of a weak solution]
	\label{weakk}
	Let $\Omega$ be a bounded domain in $\R^n$ with boundary of class $C^3$, $n = 2, 3$. \RevA{Assume (\textbf{E0})-(\textbf{E1}) and (\textbf{M0})-(\textbf{M1}) for $\Psi,\ \mathbf M,\ \nu$}.  Let \RevA{then } $\mathbf{v}_0 \in\mathbf{ L}^2_\sigma(\Omega)$
	and $\pphi_0\in \mathbf{H}^1 (\Omega)$ with $\pphi\in[0,1]^N$ and $\sum_{i=1}^N\varphi_i\equiv 1$ almost everywhere in $\Omega$. Assume also $\overline{\pphi}_0\in(0,1)^N$. Then there exists a global weak solution $(\mathbf{v},\pphi,\ww)$ of \eqref{eq:NSCH1}-\eqref{eq:NSCH7}, defined on $(0,\infty)$, in the sense of Definition \ref{weaks} with $T=\infty$.
\end{theorem}
Concerning strong solutions we have 
\begin{theorem}[Existence of a strong solution]
	\label{strong1}
	Let $\Omega$ be a bounded domain in $\R^n$ with boundary of class $C^3$, $n = 2, 3$. \RevA{Assume (\textbf{E0})-(\textbf{E1}) and (\textbf{M0})-(\textbf{M1}) for $\Psi,\ \mathbf M,\ \nu$, and assume that the mobility matrix $\mathbf{M}$ is constant.} Let \RevA{then }$\mathbf{v}_0 \in\mathbf{V}_\sigma$
	and $\pphi_0\in \mathbf{H}^2 (\Omega)$ with $\partial_\mathbf{n}\pphi_0=0$ on $\partial\Omega$ and such that $-\Delta\pphi_0+\mathbf{P}\Psip(\pphi_0)=:\ww_0\in \widetilde{\mathbf{V}}_0$, $\pphi_0\in[0,1]^N$ and $\sum_{i=1}^N\varphi_{0,i}\equiv 1$ almost everywhere in $\Omega$. Assume also $\overline{\pphi}_0\in(0,1)^N$. Then there exists a strong solution $(\mathbf{v},\pphi,\ww,p)$ of \eqref{eq:NSCH1}-\eqref{eq:NSCH7} in the sense of Definition \ref{strong}, which in the 2D case is globally defined on $(0,\infty)$, whereas in 3D it is defined on $(0,T)$, where $T>0$ can be bounded below by a constant depending only on suitable norms of the initial data.  
\end{theorem}
\begin{remark}
	Given a strong solution, the existence of pressure $p\in L^2_{uloc}(I;\mathbf{H}^1(\Omega))$ (where $I=(0,T)$ for $n=3$ and $I=[0,\infty)$ for $n=2$), $\overline{p}\equiv0$ can be retrieved in a classical way (see, e.g., \cite{Galdi}).
\end{remark}
\RevA{\begin{remark}
    We point out that we can show uniqueness of strong solutions, under the very general assumptions (\textbf{E0})-(\textbf{E1}) on the singular potential, only if we assume additionally that the phases are strictly separated from zero in the time interval of their existence, i.e., if a condition like \eqref{sep1} below holds. Indeed, in this setting one can use the weak-strong uniqueness result proven in Theorem \ref{ws} below to get the result.
\end{remark}}
By exploiting a similar proof as the one for the previous theorems, we can also prove a result concerning the multi-component Cahn-Hilliard equation with a divergence-free drift. In particular, we consider the problem 
\begin{alignat}{2}
\nonumber\partial_t\pphi+\RevA{(\nabla\pphi)\vv}-\operatorname{div}(\mathbf{M}\nabla\ww)&=0,&&\quad\text{ in } \Omega\times(0,\infty),\\
\ww&= -\Delta\pphi+\mathbf{P}(\Psi_{,\pphi}^1(\pphi)-\mathbf{A}\pphi),&&\quad\text{ in }\Omega\times(0,\infty), \label{Ch}\\ \nonumber
\partial_\mathbf{n}\pphi&=\partial_\mathbf{n}\ww=\mathbf0,&&\quad\text{ on }\partial\Omega\times(0,\infty).
\end{alignat}
We have the following 
\begin{theorem}
	\label{convective}
\RevB{Let $\Omega$ be a bounded domain in $\R^n$ with boundary of class $C^4$, $n = 2, 3$. \RevA{Assume (\textbf{E0})-(\textbf{E1}) and (\textbf{M0}), and set the mobility $\mathbf M$ to be constant. Then assume}  the initial condition $\pphi_0 \in \textbf{H}^2(\Omega)$ with $\partial_\textbf{n}\pphi_0=0$ on $\partial\Omega$ is such that $-\Delta\pphi_0+\textbf{P}\Psip(\pphi_0)=:\ww_0\in \widetilde{\textbf{V}}_0$, $\pphi_0\in[0,1]^N$, $\sum_{i=1}^N\varphi_{0,i}\equiv 1$ almost everywhere in $\Omega$ and $\overline{\pphi}_0\in(0,1)^N$. Fix then $I=(0,\infty)$ and assume that $\textbf{v} \in L^\infty(I;\textbf{H}_\sigma)\cap L^2(I; \textbf{V}_\sigma)$. Then, there exists a unique global (strong) solution to \eqref{Ch} such that}	\begin{align*}
&\pphi\in BUC(I; \mathbf{H}^1 (\Omega)) \cap  L^\infty(I;\mathbf{W}^{2,q}(\Omega))\cap H^1_{uloc}(I;\mathbf{H}^1(\Omega)) ,\\&
\sum_{i=1}^N\varphi_{i}\equiv 1, \quad \text{and}\quad\varphi_i>0\quad \forall i=1,\ldots,N  \quad \text{a.e. in }\Omega\times I,
\\& \PsipA(\pphi) \in \mathbf{L}^\infty (I; \mathbf{L}^q (\Omega)) ,\ \quad
\ww\in \mathbf{L}^\infty (I; \widetilde{\mathbf{V}}_0)\cap L^2_{uloc}(I;\mathbf{H}^3(\Omega)),
\end{align*}
with $q\in[2,\infty)$ when $n=2$ and $q=6$ when $n=3$, and it satisfies \eqref{Ch} in the almost everywhere sense. Moreover, it holds
\begin{align}
\Vert \nabla \ww\Vert_{L^\infty(0,\infty;\mathbf{L}^2(\Omega))}\leq C_1\left(\Vert \nabla\ww_0\Vert^2+\int_0^\infty\Vert \nabla\ww\Vert^2 ds+\int_0^\infty \Vert \nabla\textbf{v}\Vert^2 ds\right)^{\frac 1 2}\left(e^{C_1\int_0^\infty\Vert \nabla \textbf{v}\Vert^2ds}\right),
\label{w1}
\end{align}
together with
\begin{align}
\int_0^\infty\Vert\nabla\partial_t \pphi\Vert^2ds\leq C_1\left(\Vert \nabla\ww_0\Vert^2+\int_0^\infty\Vert \nabla\ww\Vert^2 ds+\int_0^\infty\Vert \nabla\textbf{v}\Vert^2 ds\right)\left(e^{C_1\int_0^\infty\Vert \nabla \textbf{v}\Vert^2 ds}\right),
\label{w3}
\end{align}
for some $C_1>0$.
We also have, additionally, for $k=1,2$,
\begin{align}
&\nonumber\int_0^\infty\Vert \nabla \ww\Vert_{\mathbf{H}^k(\Omega)}^2ds\\&\leq C_2\left(\Vert \nabla\ww_0\Vert^2+\int_0^\infty\Vert \nabla\ww\Vert^2 ds+\int_0^\infty \Vert \nabla\textbf{v}\Vert^2ds\right)\left(1+\int_0^\infty\Vert \nabla\textbf{v}\Vert^2ds\right)\left(e^{C_2\int_0^\infty\Vert \nabla \textbf{v}\Vert^2ds}\right),
\label{w4}
\end{align}
for some $C_2>0$.
In conclusion, in the case $\pphi$ is strictly separated, i.e., there exists $\delta\in(0,\tfrac 1 N)$ such that 
\begin{align}
\pphi(x,t)\geq \delta\quad \forall (x,t)\in\overline{\Omega\times I},
\label{separaz}
\end{align}
then it holds $\pphi\in L^\infty(I;\mathbf{H}^3(\Omega))\cap L^2_{uloc}(I;\mathbf{H}^4(\Omega))$.
\end{theorem}
\begin{remark}
\RevB{As it will be clear from the proof of Section \ref{sec:CH} below, we point out that uniqueness also holds for weak solutions (i.e., solutions with the regularity of Definition \ref{weaks}, concerning $\pphi$ and $\ww$), as long as, for instance, we assume $\textit{\textbf{v}}\in L^2(I;\mathbf{L}^3(\Omega)\cap \mathbf{L}^2_\sigma(\Omega))$.}
\end{remark}
Concerning the longtime behavior of weak solutions we have the following
\begin{theorem}[Regularity and asymptotic behavior of weak solutions]
	Let $\Omega$ be a bounded domain in $\R^n$ with boundary of class $C^4$, $n =
	2, 3$ and constant mobility matrix $\mathbf{M}$ \RevA{satisfying assumption (\textbf{M0}). Under assumptions (\textbf{E}0)-(\textbf{E}2) on the potential $\Psi$ and (\textbf{M1}) for $\nu$, }consider a global weak solution $(\mathbf{v},\pphi)$ given by Theorem \ref{weakk}.
	Then, the following results hold:
	\label{long}
	\begin{enumerate}
		\item  Global regularity of the concentration: for any $\tau > 0$, we have, letting $I=[\tau,+\infty)$,
			\begin{align*}
		&\pphi\in BUC(I; \mathbf{H}^1 (\Omega)) \cap  L^\infty(I;\mathbf{W}^{2,q}(\Omega))\cap H^1_{uloc}(I;\mathbf{H}^1(\Omega)) ,\\&
		\sum_{i=1}^N\varphi_{i}\equiv 1,\quad \text{and}\quad\varphi_i>0\quad \forall i=1,\ldots,N \quad \text{a.e. in }\Omega\times I,
		\\& \PsipA(\pphi) \in \mathbf{L}^\infty (I; \mathbf{L}^q (\Omega)) ,\quad
		\ww\in \mathbf{L}^\infty (I; \widetilde{\mathbf{V}}_0)\cap L^2_{uloc}(I;\mathbf{H}^3(\Omega)),
		\end{align*}
		with $q\in[2,\infty)$ when $n=2$ and $q=6$ when $n=3$, and it satisfies \eqref{Ch} in the almost everywhere sense. 
	
		\item Separation property: there exists $t^*>0$ and $\delta\in(0,\tfrac 1 N)$ such that 
		\begin{align}
		\pphi(x)\geq \delta\quad \forall (x,t)\in\overline{\Omega}\times[t^*,+\infty).
		\label{sep1}
  \end{align}
		\item Large time regularity for the velocity: there exists $t_R>0$ (possibly larger than $t^*$) such that 
		$$
		\mathbf{v}\in L^\infty(t_R,+\infty;\mathbf{V}_\sigma)\cap L^2(t_R;+\infty;\mathbf{W}_\sigma)\cap H^1(t_R;+\infty;\mathbf{H}_\sigma).
		$$
		\item Convergence to a stationary solution: $(\mathbf{v}(t),\pphi(t))\to (\mathbf{0},\pphi')$ in $\mathbf{H}_\sigma\times \mathbf{H}^{2r}(\Omega)$ as $t\to\infty$, for any $r\in(0,1)$,  where $\pphi'\in \mathbf{H}^2(\Omega)$ satisfies the stationary Cahn-Hilliard equation
\begin{alignat}{2}
		- \nonumber\Delta\pphi'+\mathbf{P}\Psi_{,\pphi}(\pphi')&=\overline{\mathbf{P}\Psi_{,\pphi}(\pphi')},&&\quad
		\text{ a.e. in }\Omega, \\
		\partial_\mathbf{n}\pphi'&=0,&&\quad \text{ a.e. on }\partial\Omega,		\label{steady1} \\
		\sum_{i=1}^N\varphi'_i &= 1,&& \quad \text{ in } \Omega,\nonumber
		\end{alignat}
		and $\overline{\pphi'}=\overline{\pphi}_0$.
		 
	\end{enumerate}
\end{theorem}

%{{\color{red}[TO DO: I would add a comment on existence of strong solutions only in the Cahn-Hilliard part, 
%which is possible in the approximating scheme also with $\alpha=0$. Maybe I would also add a Theorem about the existence of strong solutions only on advective Cahn-Hilliard equation, like Thm.2.4 in your work with Andrea Giorgini...]}}

In order to prove the above theorems, we perform a time discretization  scheme whose analysis is the object of the next section.
\section{Time discretization scheme}
\label{timediscr}
We introduce the sets:
\begin{align}
&\mathbf{H}_\mathbf{n}^2:=\{\mathbf f\in \mathbf{H}^2(\Omega)\cap\widetilde{\mathbf{H}}_0:\ \partial_\mathbf{n}\mathbf f=\mathbf{0}\text{ on }\partial\Omega\},
\\&\mathcal{Z}_1:=\{\pphi\in \mathbf{H}^2(\Omega):\ 0\leq \pphi,\ \sum_{i=1}^N\varphi_i=1, \overline{\pphi}\in(0,1)^N,\ \partial_\mathbf{n}\pphi=\mathbf{0}\text{ on } \partial\Omega \},\\&
\mathcal{Z}_2:=\{\pphi\in \mathcal{Z}_1:\ \Psi^1_{,\pphi}(\pphi)\in \mathbf{L}^2(\Omega)\}.
\end{align}
Notice that the conditions in $\mathcal{Z}_2$ also imply that $0< \pphi< 1$ almost everywhere in $\Omega$. Furthermore, for any $\pphi\in \mathcal{Z}_1$ it holds $\Psi(\pphi)\in \mathbf{L}^\infty(\Omega)$.
Setting $(\mathbf{v},\pphi,\ww):=(\mathbf{v}^{k+1},\pphi^{k+1},\ww^{k+1})$, we use the following scheme, for $h:=\frac 1 M$, $M\in\N$: let $\alpha\in[0,1].$ Given the couple $(\mathbf{v}^{k},\pphi^{k})$, $k\in \N_0:=\N\cup \{0\}$, with $\mathbf{v}^k\in \mathbf{H}_\sigma$ if $\alpha=0$, $\mathbf{v}^k\in \mathbf{W}_\sigma$ if $\alpha>0$, $\pphi^k\in \mathcal{Z}_1$, together with $\rho^k=\sum_{i=1}^N\widetilde{\rho}_i\varphi_i^k$, find $(\mathbf{v},\pphi,\ww)$ such that $\mathbf{v}\in \mathbf{V}_\sigma$, $\pphi\in\mathcal{Z}_2$ and $\ww\in \mathbf{H}_\mathbf{n}^2$, satisfying
  \begin{alignat}{2}
\nonumber\left(\dfrac{\rho \mathbf{v}-\rho^k\mathbf{v}^k}{h},\boldsymbol\psi\right)_{\Omega}+\alpha&\left(\frac{D \mathbf{v}-D\mathbf{v}^k}{h},\nabla\boldsymbol\psi\right)_{\Omega}+\left(
\operatorname{div}(\rho^k\mathbf{v}\otimes\mathbf{v})+\operatorname{div}(\mathbf{v}\otimes\mathbf{J}_\rho),\boldsymbol\psi\right)_{\Omega}
\\+\left(2\nu(\pphi^k)D\mathbf{v},D\boldsymbol\psi\right)_\Omega&\nonumber=\left((\nabla\pphi^k)^T\ww,\boldsymbol\psi\right)_\Omega,&&\quad \forall \boldsymbol\psi\in \mathbf{V}_\sigma,\\ \nonumber
\dfrac{\pphi-\pphi^k}{h}+(\nabla\pphi^k)\mathbf{v}&=\operatorname{div}(\mathbf{M}(\pphi^k)\nabla\ww)&&\quad\text{ in }\Omega,\\
\ww&= -\Delta\pphi+\mathbf{P}\left(\Psi_{,\pphi}^1(\pphi)-\frac  {\mathbf{A}\pphi}{2}-\frac{\mathbf{A}\pphi^k}{2}\right)\label{phi0}&&\quad\text{ in }\Omega,\\
\mathbf{J}_\rho&=-\nonumber\nabla\ww^T\widetilde{\mathbf{M}}(\pphi^k)&&\quad\text{ in }\Omega,\\
\partial_\mathbf{n}\pphi&=\mathbf{0}&&\quad \nonumber\text{ on }\partial\Omega,\\
\partial_\mathbf{n}\ww& =\mathbf{0}&&\quad \text{ on }\partial\Omega.\nonumber
\end{alignat}
where we set $\widetilde{\mathbf{M}}(\pphi):=\mathbf{M}(\pphi)\widetilde{\boldsymbol{\rho}}$. We also set, at least formally, $\ww_0:=-\Delta\pphi^0+\mathbf{P}\left(\Psi_{,\pphi}^1(\pphi^0)- {\mathbf{A}\pphi^0}\right).$
%{\color{red}[I removed the $\Gamma$ in front of the laplacian of $\pphi$, since it would generate cross terms related to $\mathbf{P}\Gamma\Delta\pphi$ and we cannot deal with them in the $L^2$ estimate of $\PsipA(\pphi)$...]}

%{\color{red}[We can then substitute the quadratic perturbation $\mathbf{A}\pphi$ with a more general concave $\Psi^2$...]}
\begin{remark}
	With respect to the scheme adopted in \cite[Sec.4]{ADG}, in this case we need to consider the (formal) addition of the term $\alpha\frac{\AA\mathbf{v}-\AA\mathbf{v}^k}{h}$ to gain extra regularity for the velocity $\mathbf{v}$. This is due to the fact that we aim at exploiting the same approximating scheme also to derive the existence of strong solutions.
\end{remark}
\begin{remark}
Notice that, by integrating \eqref{phi0}$_2$ over $\Omega$, we easily infer that 
\begin{align}
\overline{\pphi}\equiv \overline{\pphi}_k\in (0,1)^N\cap    \Sigma,\quad \forall k\in \N_0,
\label{conserv}
\end{align}
 i.e., the conservation of mass is preserved also at the discrete level.
 \end{remark}
\begin{remark}
	We observe that, when $\pphi\in[0,1]^N$ and $\sum_{i=1}^N\varphi_i=1$, it holds $\rho\geq \min_{i=1,\ldots,N}\widetilde{\rho}_i>0$.
	\label{minrho} 
\end{remark}
 \begin{remark}
 Observe that, by multiplying \eqref{phi0} by $\widetilde{\boldsymbol{\rho}}$ and recalling that $\rho^k=\sum_{i=1}^N\varphi_i^k\widetilde{\rho}_i=\widetilde{\boldsymbol{\rho}}\cdot \pphi^k$ for any $k$, we end up with 
 \begin{align}
 \frac{\rho-\rho^k}{h}+\nabla\rho^k\cdot\mathbf{v}+\operatorname{div}(\mathbf{J}_\rho)=0,
 \label{eqro}
 \end{align}
 so that we can also write, recalling $\operatorname{div}(\mathbf{v}\otimes\mathbf{J}_\rho)=\operatorname{div}(\mathbf{J}_\rho)\mathbf{v}+(\mathbf{J}_\rho\cdot\nabla)\mathbf{v}$,
 
 \begin{align}
(\operatorname{div}(\mathbf{v}\otimes\mathbf{J}_\rho),\boldsymbol\psi)_\Omega=((\mathbf{J}_\rho\cdot\nabla)\mathbf{v},\boldsymbol\psi)_\Omega+\frac{1}{2}(\operatorname{div}(\mathbf{J}_\rho)\mathbf{v},\boldsymbol\psi)_\Omega-\frac 1 2(\frac{\rho-\rho^k}{h},\mathbf{v}\cdot\boldsymbol\psi)_\Omega-\frac 1 2 (\nabla\rho^k\cdot \mathbf{v},\mathbf{v}\cdot\boldsymbol\psi)_\Omega.
 \label{Jrho}
 \end{align}
 \end{remark}
We thus obtain the following existence theorem for the discretized problem:
\begin{theorem}[Existence of a discretized solution]
	\label{disc}
We have the following results: assume $\Omega\subset\R^n$, $n=2,3$ a bounded domain with boundary of class $C^3$.
\begin{enumerate}
	\item Let the couple $(\mathbf{v}^{k},\pphi^{k})$, $k\in \N_0$, be given, with $\mathbf{v}^k\in \mathbf{H}_\sigma$ if $\alpha=0$, $\mathbf{v}^k\in \mathbf{W}_\sigma$ if $0<\alpha\leq 1$, $\pphi^k\in \mathcal{Z}_1$, together with $\rho_k=\sum_{i=1}^N\widetilde{\rho}_i\varphi_i^k$. Then there are some $(\mathbf{v},\pphi,\ww)$ such that $\mathbf{v}\in \mathbf{V}_\sigma$, $\pphi\in\mathcal{Z}_2$ and $\ww\in \mathbf{H}_\mathbf{n}^2$, satisfying \eqref{phi0}, together with the discrete energy estimate
	\begin{align}
	\nonumber&E_{tot}(\rho,\mathbf{v},\pphi)+\frac{\alpha }{2}\int_\Omega \left\vert D\mathbf{v}\right\vert^2dx  +h\int_\Omega \mathbf{M}(\pphi^k)\nabla \ww:\nabla\ww dx+\int_\Omega \rho^k\frac{\vert \mathbf{v}-\mathbf{v}^k\vert^2}2dx+\frac{\alpha }{2}\int_\Omega \left\vert D\mathbf{v}- D\mathbf{v}^k\right\vert^2dx\\&+2h\int_\Omega \nu(\pphi^k)\vert D\mathbf{v}\vert^2dx
	+\frac 1 2\int_\Omega\vert \nabla\pphi-\nabla \pphi^k\vert^2dx\leq E_{tot}(\rho^k,\mathbf{v}^k,\pphi^k)+\frac{\alpha }{2}\int_\Omega \left\vert D\mathbf{v}^k\right\vert^2dx.
	\label{ener}
	\end{align}
	Moreover, there exists $\delta_{k+1}\in(0,\frac 1 N)$, possibly depending on $k$, such that 
	\begin{align}
	\delta_{k+1}<\pphi(x)<1-(N-1)\delta_{k+1},\quad \forall x\in\overline{\Omega}.
	\label{sepk}
	\end{align}
	\item Assume additionally to (1) that $\pphi^k$ is such that $$\ww^k:=-\Delta\pphi^k+\mathbf{P}\left(\Psi_{,\pphi}^1(\pphi^k)-\frac  {\mathbf{A}\pphi^k}{2}-\frac{\mathbf{A}\pphi^{k-1}}{2}\right)\in \widetilde{\mathbf{V}}_0\cap \mathbf{L}^\infty(\Omega),$$ with $\pphi^{-1}=\pphi^0$ when $k=0$. Assume also that the mobility matrix $\mathbf{M}$ is constant and that $\overline{\pphi}^{k-1}=\overline{\pphi}^k=\overline{\pphi}^0$. Then it also holds
	\begin{align}
	&\nonumber\frac 1 {2} \left(\mathbf{M}\nabla \ww,\nabla\ww\right)_\Omega+\frac h 2 \left\Vert\frac{\nabla\pphi-\nabla\pphi^k}{h}\right\Vert^2\\&\leq \frac 1 {2} \left(\mathbf{M}\nabla \ww^k,\nabla\ww^k\right)_\Omega+hC(\min_{i=1,\ldots,N}\overline{\varphi}^0_i)\Vert \nabla\mathbf{v}\Vert\Vert \mathbf{v}\Vert(1+\Vert \nabla\ww^k\Vert^2+\Vert \nabla\ww\Vert^2)\nonumber\\&+hC(\min_{i=1,\ldots,N}\overline{\varphi}^0_i)\Vert \nabla\mathbf{v}\Vert^2(1+\Vert \nabla\ww^k\Vert^2)+hC\left(\left\Vert \frac{\pphi-\pphi^k}{h}\right\Vert^{2}_{\mathbf{V}_0'}+\left\Vert \frac{\pphi^k-\pphi^{k-1}}{h}\right\Vert^2_{\mathbf{V}_0'}\right),
	\label{h}
	\end{align}
	for some constants $C(\min_{i=1,\ldots,N}\overline{\varphi}^0_i),C>0$ independent of $h$, $k$ and $\alpha$.
	\item In conclusion, if, additionally to (2), $0<\alpha\leq 1$, and thus we assume $\mathbf{v}^k\in\mathbf{W}_\sigma$, then it also holds
	
		\begin{align}
\alpha\left\Vert\dfrac{ D \mathbf{v}-D \mathbf{v}^k}{h}\right\Vert^2\leq \frac{C}{\alpha}\left(\Vert \mathbf{w}\Vert_{\mathbf{H}^1(\Omega)}^2\Vert \nabla\pphi^k\Vert^2+\Vert\ww\Vert^2_{\mathbf{H}^2(\Omega)}\Vert \nabla \mathbf{v}\Vert^2+\Vert\nabla\mathbf{v}\Vert^4+\Vert\nabla\mathbf{v}\Vert^2\right),
	\label{c1}
	\end{align}
	\begin{align}
	&\nonumber\frac\alpha{2h}\left\Vert\AA\mathbf{v}\right\Vert^2\leq \frac\alpha{2h}\Vert\AA\mathbf{v}^k\Vert^2+C\left\Vert\AA\mathbf{v}\right\Vert^2\\&+C\left(\Vert \nabla\mathbf{v}\Vert^6+\Vert \nabla\ww\Vert^2\Vert\ww\Vert^2_{\mathbf{H}^2(\Omega)}\Vert \nabla \mathbf{v}\Vert^2+\Vert \nabla\pphi^k\Vert_{\mathbf{L}^6(\Omega)}^4\Vert \nabla\mathbf{v}\Vert^2\right. \nonumber\\&\left.+\Vert \mathbf{w}\Vert_{\mathbf{H}^3(\Omega)}^2\Vert \nabla\pphi^k\Vert^2+\left\Vert\dfrac{ D\mathbf{v}-D\mathbf{v}^k}{h}\right\Vert^2\right),
	\label{c2}
	\end{align}
	\begin{align}
	\nonumber&\alpha\left\Vert\dfrac{\AA\mathbf{v}-\AA\mathbf{v}^k}{h}\right\Vert^2\leq \frac{C}{\alpha}\left(\Vert \AA\mathbf{v}\Vert^2+\left\Vert\dfrac{D\mathbf{v}-D\mathbf{v}^k}{h}\right\Vert^2+\Vert \nabla\ww\Vert^2\Vert\ww\Vert^2_{\mathbf{H}^3(\Omega)}\Vert \nabla \mathbf{v}\Vert^2\right.\nonumber\\&\left.+\Vert \nabla\pphi^k\Vert_{\mathbf{L}^6(\Omega)}^4\Vert \nabla\mathbf{v}\Vert^2+\Vert\nabla\mathbf{v}\Vert^6+\Vert \mathbf{w}\Vert_{\mathbf{H}^2(\Omega)}^2\Vert \nabla\pphi^k\Vert^2\right),
	\label{c3}
	\end{align}
	for some constant $C>0$ independent of $h$, $k$ and $\alpha$.
\end{enumerate}
\end{theorem}
\begin{remark}
	Notice that the assumption at point (2) already implies, by Theorem \ref{steaddy} point (2), that there exists $\delta_k$ such that 
	\begin{equation}
	\delta_k<\pphi^k(x)<1-(N-1)\delta_k \label{prop2}
	\end{equation}%
	for any $x\in \overline{\Omega }$, i.e., $\pphi_k$ has to be strictly separated from pure phases.
	\label{separ}
\end{remark}
 \begin{remark}
	We point out that, assuming the matrix $\mathbf{M}$ constant, under the assumptions of point (3) by elliptic regularity for $\ww$ in \eqref{phi0}$_2$ we deduce, since $\partial\Omega$ is of class $C^4$, that $\ww\in \mathbf{H}^3(\Omega)$. Indeed, $\pphi,\pphi^k\in \mathbf{H}^2(\Omega)$ and
	thanks to Theorem \ref{steaddy} point (3) applied to the couple $(\mathbf f,\mathbf{m})=\left(\ww^k+\mathbf{P}\left(\frac  {\mathbf{A}\pphi^k}{2}+\frac{\mathbf{A}\pphi^{k-1}}{2}\right),\overline{\pphi}^k\right)$, we also get $\pphi^k\in \mathbf{W}^{2,p}(\Omega)$, $p\in[2,\infty)$ for $n=2$ and $p\in[2,6]$ for $n=3$, and thus $(\nabla\pphi^k)\mathbf{v}\in \mathbf{H}^1(\Omega)$.
	Moreover, it is easy to see that, by applying \cite[Thm. A.1]{GGW} to the problem (in weak formulation)
	\begin{align*}
	\left(\left(2\nu(\pphi^k)+\frac{\alpha}{h}\right)D\mathbf{v},D\boldsymbol\psi\right)_\Omega=(\mathbf g,\boldsymbol \psi)_\Omega,\quad \forall \boldsymbol\psi\in \mathbf{V}_\sigma,
	\end{align*}
	where 
	$$\mathbf g:=-\mathbb{P}\dfrac{\rho \mathbf{v}-\rho^k\mathbf{v}^k}{h}+\frac{\alpha}{h}\AA\mathbf{v}^k-
	\mathbb{P}\operatorname{div}(\rho^k\mathbf{v}\otimes\mathbf{v})-\mathbb{P}\operatorname{div}(\mathbf{v}\otimes\mathbf{J}_\rho)\in \mathbf{L}^{\frac 32}(\Omega),
$$
due to $\mathbf{v}^k\in \mathbf{W}_\sigma$, we easily deduce that $\mathbf{v}\in \mathbf{W}^{2,\frac 3 2}(\Omega)\cap \mathbf{V}_\sigma$, and thus by bootstrapping we end up with $\mathbf{v}\in \mathbf{W}_\sigma$.

	\label{H3}
\end{remark}
\subsection*{Proof of Theorem \ref{disc}}
	\begin{enumerate}
	\item[\underline{Point (1)}] First we show the \textit{a priori} estimate \eqref{ener} assuming to have a triple $(\mathbf{v},\pphi,\ww)$ such that $\mathbf{v}\in \mathbf{V}_\sigma$, $\pphi\in\mathcal{Z}_2$ and $\ww\in \mathbf{H}_\mathbf{n}^2$, satisfying \eqref{phi0}, given the couple $(\mathbf{v}^k,\pphi^k)$. 	Observe that property \eqref{sepk} is then clearly valid thanks to Theorem \ref{steaddy} point (2), being $\mathbf f=\ww+\mathbf{P}\mathbf{A}\frac{\pphi-\pphi^k}{2}\in \mathbf{L}^\infty(\Omega)$, by the embedding $\mathbf{H}^2_\mathbf{n}\hookrightarrow \mathbf{H}^2(\Omega)\hookrightarrow \mathbf{L}^\infty(\Omega)$. This means that we can work with the singular potential $\PsipA(\pphi)$ as if it were regular, without the necessity of making lengthy regularization procedures.
	
	 Let us first observe that 
	\begin{align}
	0=&\int_\Omega\operatorname{div}\left(\mathbf{J}_\rho\frac{\vert \mathbf{v}\vert^2}2\right)dx=\int_\Omega \operatorname{div}(\mathbf{J}_\rho)\frac {\vert \mathbf{v}\vert^2}2dx+\int_\Omega((\mathbf{J}_\rho\cdot\nabla)\mathbf{v})\cdot\mathbf{v}dx
	\label{v1}
	\end{align}
	and 
		\begin{align}
	0=&\int_\Omega \operatorname{div}\left(\rho^k\mathbf{v}\frac{\vert\mathbf{v}\vert^2}{2}\right)dx=\int_\Omega (\nabla\rho^k\cdot \mathbf{v})\frac{\vert\mathbf{v}\vert^2}{2}dx+\int_\Omega \rho^k((\mathbf{v}\cdot\nabla)\mathbf{v})\cdot \mathbf{v}dx,
	\label{v2}
	\end{align}
		recalling that $\operatorname{div}\mathbf{v}=0$,	so that 
	$$
	\int_\Omega(\operatorname{div}(\rho^k\mathbf{v}\otimes\mathbf{v})\cdot \mathbf{v}-(\nabla\rho^k\cdot \mathbf{v})\frac{\vert\mathbf{v}\vert^2}{2})dx=0.
	$$
	By exploiting the basic identity 
 \begin{align}
 \mathbf{a}\cdot(\mathbf{a}-\mathbf{b})=\frac{\vert \mathbf{a}\vert^2}2-\frac{\vert \mathbf{b}\vert^2}2+\frac{\vert \mathbf{a}-\mathbf{b}\vert^2}2,
 \label{identity1}
 \end{align}
 we also obtain (see also \cite{ADG})
	\begin{align}
	&\nonumber(\rho\mathbf{v}-\rho^k\mathbf{v}^k)\cdot\mathbf{v}=(\rho-\rho^k)\vert\mathbf{v}\vert^2+\rho^k(\mathbf{v}-\mathbf{v}^k)\cdot\mathbf{v}\\&=\rho\frac{\vert\mathbf{v}\vert^2}2-\rho^k\frac{\vert\mathbf{v}^k\vert^2}2+(\rho-\rho^k)\frac{\vert\mathbf{v}\vert^2}2+\rho^k\frac{\vert \mathbf{v}-\mathbf{v}^k\vert^2}{2}.
	\label{diffv}
	\end{align}
	Analogously,
	\begin{align*}
	\frac{D(\mathbf{v}-\mathbf{v}^k)}{h}: D \mathbf{v}=\frac{1}{2h}\vert D\mathbf{v}\vert^2-\frac{1}{2h}\vert D\mathbf{v}^k\vert^2+\frac{1}{2h}\vert D\mathbf{v}-D\mathbf{v}^k\vert^2.
	\end{align*}

	Putting everything together, from \eqref{phi0}$_1$ with $\boldsymbol\psi=\mathbf{v}$ we deduce 
	\begin{align}
	&\nonumber 0=\int_\Omega\left(\frac{\rho\vert\mathbf{v}\vert^2-\rho^k\vert \mathbf{v}_k\vert^2}{2h}\right)dx+\alpha\int_\Omega\left(\frac{1}{2h}\vert D\mathbf{v}\vert^2-\frac{1}{2h}\vert D\mathbf{v}^k\vert^2+\frac{1}{2h}\vert D\mathbf{v}-D\mathbf{v}^k\vert^2\right)dx \\&+\int_\Omega \frac{\rho^k\vert\mathbf{v}-\mathbf{v}^k\vert^2}{2h}dx+\int_\Omega 2\eta(\pphi^k)\vert D\mathbf{v}\vert^2dx-\int_\Omega (\nabla\pphi^k)^T\ww\cdot\mathbf{v}dx.
	\label{vel} 
	\end{align} 
	Then, testing \eqref{phi0}$_2$ with $\ww$ we get 
	\begin{align}
	0=\int_\Omega \frac{\pphi-\pphi^k}h\ww dx +\int_\Omega (\nabla\pphi^k)\mathbf{v}\cdot\ww dx+\int_\Omega \nabla\ww :\mathbf{M}(\pphi^k)\nabla\ww dx,
	\label{pphi}
	\end{align}
	and testing \eqref{phi0}$_3$ by $\frac{\pphi-\pphi^k}h$ we instead deduce, since clearly $\frac{\pphi-\pphi^k}h\in T\Sigma$ for a.e.\ $x\in\Omega$,
	$$
	\int_\Omega \frac{\pphi-\pphi^k}h\cdot \ww dx=\int_\Omega\nabla\pphi:\nabla\frac{\pphi-\pphi^k}{h} dx+\int_\Omega \PsipA(\pphi)\cdot \frac{\pphi-\pphi^k}h dx-\int_\Omega \mathbf{A}\frac{\pphi+\pphi^k}2\cdot\frac{\pphi-\pphi^k}h dx.
	$$
	Now, by monotonicity of $\Psi^1$,
	$$
	(\PsipA(\pphi),\frac{\pphi-\pphi^k}h)\geq \int_\Omega \dfrac{\Psi^1(\pphi)-\Psi^1(\pphi^k)}{h}dx.
	$$
	Moreover, we have, using that is $\mathbf{A}$ symmetric, 
	$$
	-\int_\Omega \mathbf{A}\frac{\pphi+\pphi^k}2\cdot\frac{\pphi-\pphi^k}h dx=-\int_\Omega \frac{1}{2h}\pphi^T\mathbf{A}\pphi+\int_\Omega \frac{1}{2h}(\pphi^k)^T\mathbf{A}\pphi^k.
	$$
	Again by the same identity \eqref{identity1} as before we also have
	\begin{align*}
	\int_\Omega\nabla\pphi:\nabla\frac{\pphi-\pphi^k}{h} dx= \frac{1}{2h}\int_\Omega \vert \nabla\pphi\vert^2dx-\frac{1}{2h}\int_\Omega \vert \nabla\pphi^k\vert^2dx+\frac{1}{2h}\int_\Omega \vert \nabla(\pphi-\pphi^k)\vert^2dx.
	\end{align*}
	Therefore, summing up \eqref{vel} with \eqref{pphi}, recalling that $\int_\Omega (\nabla\pphi^k)\mathbf{v}\cdot\ww dx=\int_\Omega (\nabla\pphi^k)^T\ww\cdot\mathbf{v}dx$, we infer, multiplying everything by $h$,
	\begin{align}
\nonumber&E_{tot}(\rho,\mathbf{v},\pphi)+\frac{\alpha }{2}\int_\Omega \left\vert D\mathbf{v}\right\vert^2dx  +h\int_\Omega \mathbf{M}(\pphi^k)\nabla \ww:\nabla\ww dx+\int_\Omega \rho^k\frac{\vert \mathbf{v}-\mathbf{v}^k\vert^2}2dx+\frac{\alpha }{2}\int_\Omega \left\vert D\mathbf{v}-D\mathbf{v}^k\right\vert^2dx\\&+2h\int_\Omega \nu(\pphi^k)\vert D\mathbf{v}\vert^2dx
+\frac 1 2\int_\Omega\vert \nabla\pphi-\nabla \pphi^k\vert^2dx\leq E_{tot}(\rho^k,\mathbf{v}^k,\pphi^k)+\frac{\alpha }{2}\int_\Omega \left\vert D\mathbf{v}^k\right\vert^2dx,
\label{ener1}
\end{align}
	which is exactly \eqref{sepk}. In order to show existence of weak solutions (we perform the proof in the 3D case, the two dimensional one being analogous) we want to use the Leray-Schauder's principle as in \cite{ADG}. In particular, we define the operators $\mathcal{L}^k,\mathcal{F}^k:X\to Y$, where
	\begin{align*}
&	X:=\mathbf{V}_\sigma\times \mathcal{Z}_2\times \mathbf{H}_\mathbf{n}^2,\\&
Y:=\mathbf{V}_\sigma^\prime\times \widetilde{\mathbf{H}}_0\times \widetilde{\mathbf{H}}_0,
	\end{align*}
and, for any triple $(\mathbf{u},\pphi,\ww)\in X$,
\begin{align}
\mathcal{L}^k(\mathbf{u},\pphi,\ww):=\begin{bmatrix}
L^k_\alpha(\mathbf{u})\\
-\operatorname{div}(\mathbf{M}(\pphi^k)\nabla\ww)+\int_\Omega\ww dx\\
-\Delta\pphi+\mathbf{P}\PsipA(\pphi)
\end{bmatrix}
\label{L}
\end{align} 
with 
$$
<L^k_\alpha(\mathbf{v}),\boldsymbol\psi>:=\int_\Omega \left(\nu(\pphi^k)+\frac \alpha h \right)D\mathbf{v}:D\boldsymbol\psi dx\quad \forall \boldsymbol\psi\in \mathbf{V}_\sigma,
$$
and the second and third line regarded pointwise, and 
\begin{align}
\label{F}
\mathcal{F}^k(\mathbf{u},\pphi,\ww):=\begin{bmatrix}
-\frac{\rho\mathbf{v}-\rho^k\mathbf{v}^k}{h}+\frac{\alpha}{h}\AA\mathbf{v}^k-\operatorname{div}(\rho^k\mathbf{v}\otimes\mathbf{v})+(\nabla\pphi^k)^T\ww\\-\left(\operatorname{div}(\mathbf{J}_\rho)-\frac{\rho-\rho^k}h-\mathbf{v}\cdot\nabla\rho^k\right)\frac {\mathbf{v}}2-(\mathbf{J}_\rho\cdot\nabla)\mathbf{v}\\ \\ 
-\frac{\pphi-\pphi^k}h-(\nabla\pphi^k)\mathbf{v}+\int_\Omega\ww dx\\ \\ 
\ww+\mathbf{P}\mathbf{A}\frac{\pphi+\pphi^k}{2}
\end{bmatrix},
\end{align}
where we exploited \eqref{Jrho} in the first line. Then $\mathbf{z} = (\mathbf{u},\pphi,\ww) \in X$ is a weak solution of the time discrete problem \eqref{phi0} if and only if
	\begin{align}
	\label{ll}
	\mathcal{L}^k(\mathbf{u},\pphi,\ww)-\mathcal{F}^k(\mathbf{u},\pphi,\ww)=0.
	\end{align}
From standard theory of partial differential equations (e.g., the Lax-Milgram Lemma) we get the invertibility of the linear continuous operator $$L_k^\alpha:\mathbf{V}_\sigma\to \mathbf{V}_\sigma^\prime.$$
Now we consider for given $\mathbf f \in \widetilde{\mathbf{H}}_0$ the elliptic boundary value problem
\begin{alignat*}{2}
-\operatorname{div}(\mathbf{M}(\pphi^k)\nabla\ww)+\int_\Omega\ww dx&=\mathbf f&&\quad\text{in }\Omega,\\
\partial_\mathbf{n}\ww&=\mathbf{0}&&\quad \text{on }\partial\Omega.
\end{alignat*}
By the Lax-Milgram Lemma there exists a unique solution $\ww\in \widetilde{\mathbf{V}}_0$. Indeed, due the properties of the matrix $\mathbf{M}$, we can rewrite the problem in weak form as: find $\ww \in \widetilde{\mathbf{V}}_0$ such that 
\begin{align}
\label{weak}
\int_\Omega \mathbf{M}(\pphi^k)\nabla\ww:\nabla \boldsymbol\psi +\int_\Omega \boldsymbol\psi dx \cdot \int_\Omega \ww dx=\int_\Omega \mathbf f\cdot \boldsymbol\psi dx\quad \forall \boldsymbol\psi\in \widetilde{\mathbf{V}}_0,
\end{align}
and by \eqref{nondeg} we clearly see that the bilinear form on $\widetilde{\mathbf{V}}_0\times \widetilde{\mathbf{V}}_0$ on the left-hand side is continuous and coercive ($\vert \mathbf{M}(\mathbf{s})\vert $ is bounded above by some constant, see assumption (\textbf{M0})).  Since clearly $\mathbf f\in\widetilde{\mathbf{V}}_0^\prime$, then by Lax\RevA{-}Milgram Lemma there exists a unique $\ww\in \widetilde{\mathbf{V}}_0$ satisfying the weak problem \eqref{weak} and such that 
\begin{align}
\Vert \ww\Vert_{\mathbf{H}^1(\Omega)}\leq C_k\Vert\mathbf f\Vert.
\label{f}
\end{align}
We now want to show by a bootstrapping argument that actually we have $\ww\in \mathbf{H}^2(\Omega)\cap \widetilde{\mathbf{H}}_0$. Indeed, $\ww$ formally solves 
\begin{alignat}{2}
-\Delta \ww&=(\mathbf{M}(\pphi^k)+\mathbf{e}\otimes \mathbf{e})^{-1}\left(\mathbf{P}(\nabla\mathbf{M}(\pphi^k):\nabla\ww)-\int_\Omega \ww dx+\mathbf f\right) &&\quad \text{ in }\Omega,\label{p}\\
\partial_\mathbf{n}\ww&=0&&\quad \text{ on }\partial\Omega,
\end{alignat}
where $\mathbf{e}=(1,\ldots,1)^T$. \RevA{Note that here we have adopted the notation $$(\nabla\mathbf{M}(\pphi^k):\nabla\ww)_i:= \sum_{k=1}^N \nabla \mathbf M_{ik}\cdot \nabla w_k,$$ for any component $i=1,\ldots,N$.} Clearly, for any $x\in\Omega$, if at the same point a function $\mathbf g$ belongs to $T\Sigma$, then also $\widehat{\mathbf g}:=(\mathbf{M}(\pphi^k(x))+\mathbf{e}\otimes \mathbf{e})^{-1}\mathbf{g}\in T\Sigma$, so that the problem above is well posed. Indeed, since $\mathbf{M}(\cdot)$ is positive definite on $T\Sigma$, we have that there exists $\widetilde{\mathbf g}\in T\Sigma$ such that $\mathbf g=\mathbf{M}(\pphi^k(x))\widetilde{\mathbf g}=(\mathbf{M}(\pphi^k(x))+\mathbf{e}\otimes \mathbf{e})\widetilde{\mathbf g}$. This implies that $\widehat{\mathbf g}=\widetilde{\mathbf g}$ and thus $ \widehat{\mathbf g}\in T\Sigma$ as well. Observe that, due to $\pphi^k\in \mathbf{H}^2(\Omega)$, by the embedding $\mathbf{H}^2(\Omega)\hookrightarrow\mathbf{W}^{1,6}(\Omega)$, recalling that $\mathbf{M}^\prime(\pphi^k)$ is bounded by (\textbf{M}0) and $\pphi^k\in[0,1]^N$ by assumption, 
$$
\Vert \mathbf{P}(\nabla\mathbf{M}(\pphi^k):\nabla\ww)\Vert_{\mathbf{L}^{\frac 3 2}(\Omega)}\leq C\Vert \nabla\mathbf{M}(\pphi^k):\nabla\ww\Vert_{\mathbf{L}^{\frac 3 2}(\Omega)}\leq C\Vert \pphi^k\Vert_{\mathbf{W}^{1,6}(\Omega)}\Vert \nabla\ww\Vert\leq C_k\Vert \mathbf f\Vert.
$$
By elliptic regularity we thus obtain from \eqref{f} and \eqref{p} that $\Vert \ww\Vert_{ \mathbf{W}^{2,\frac 3 2}(\Omega)}\leq C_k\Vert \mathbf f\Vert$. Observe now that $\mathbf{W}^{2,\frac 3 2}(\Omega)\hookrightarrow \mathbf{W}^{1,3}(\Omega)$, so that we can write, recalling the previous estimates, 
$$
\Vert \mathbf{P}(\nabla\mathbf{M}(\pphi^k):\nabla\ww)\Vert\leq \Vert \nabla\mathbf{M}(\pphi^k):\nabla\ww\Vert\leq C\Vert \pphi^k\Vert_{\mathbf{W}^{1,6}(\Omega)}\Vert \ww\Vert_{\mathbf{W}^{1,3}(\Omega)}\leq C_k\Vert \mathbf f\Vert,
$$
 implying, again by elliptic regularity from \eqref{p}, $\ww \in \mathbf{H}^2(\Omega)$ and thus $\ww\in \mathbf{H}^2_n$ as desired, since $\ww$ also solves \eqref{weak}. We also deduce the inequality 
 \begin{align}
 \Vert \ww\Vert_{\mathbf{H}^2(\Omega)}\leq C_k\Vert \mathbf f\Vert.
 \label{ineq}
 \end{align}
 Clearly this last inequality, as the problem \eqref{weak} is linear, ensures that if $\mathbf f_m\to \mathbf f$ in $\mathbf{L}^2(\Omega)$ as $m\to \infty$ then the corresponding $\ww_m$ converge to $\ww$ in $\mathbf{H}^2(\Omega)$.
 
 Concerning the third row of $\mathcal{L}^k$, we easily see by Theorem \ref{steaddy} point (1) that this operator is invertible from $\widetilde{\mathbf{H}}_0$ to $\mathcal{Z}_2$, and it is also continuous from $\widetilde{\mathbf{H}}_0$ to $\mathbf{H}^{2-{s}}(\Omega)$ with $s\in(0,\tfrac 1 4)$. Altogether we obtain that $\mathcal{L}^k$ is invertible from $Y$ to $X$, but since $X$ is not a Banach space due to $\mathcal{Z}_2$ ($\mathcal{Z}_2$ is not even a vector space), we need to redefine new spaces: in particular, let us set, for some $s\in(0,\tfrac1 4)$,
	\begin{align*}
&	\widetilde{X}:=\mathbf{V}_\sigma\times \mathbf{H}^{2-s}(\Omega)\times \mathbf{H}_\mathbf{n}^2,\ \quad
\widetilde{Y}:=\mathbf{L}^{\frac 3 2}(\Omega)\times \mathbf{W}^{1,\frac 3 2}(\Omega)\times \widetilde{\mathbf{V}}_0.
\end{align*}     
Now, this means that $(\mathcal{L}^k)^{-1}$ is continuous from $Y$ to $\widetilde{X}$. We also observe that $\widetilde{Y}\hookrightarrow\hookrightarrow Y$, so that the restriction $(\mathcal{L}^k)^{-1}_{\vert \widetilde{Y}}:\widetilde{Y}\to (\mathcal{L}^k)^{-1}(\widetilde{Y})\hookrightarrow\hookrightarrow\widetilde{X}$ is a compact operator, since it is the composition of a compact and a bounded linear operator. By exactly the same argument as in \cite[Sec.4.2]{ADG}, with the only adaptations to consider vectors $\pphi,\pphi^k,\ww$ instead of scalars (recall that $\mathbf{M}(\pphi^k)$ is bounded), setting the operator $A$ as to be the identity and recalling $\mathbf{v}^k\in \mathbf{W}_\sigma$, we deduce that the operator $\mathcal{F}^k_{\vert \widetilde{X}}:\widetilde{X}\to \widetilde{Y}$ is continuous and maps bounded sets into bounded sets. In order to apply a fixed point theorem, we note that if there exists $\mathbf h\in \widetilde{Y}\hookrightarrow\hookrightarrow Y$ solving
$$
\mathbf h -\mathcal{F}^k\circ(\mathcal{L}^k)^{-1}(\mathbf h)=0,
$$  
then there exists $\mathbf{z}\in {X}$ solution to \eqref{ll}, just by setting $\mathbf{z}=(\mathcal{L}^k)^{-1}(\mathbf h)\in X$. We now define the operator $\mathcal{K}^k:=
 \mathcal{F}^k\circ(\mathcal{L}^k)^{-1}: \widetilde{Y}\to \widetilde{Y}$, which is continuous and compact on $\widetilde{Y}$, being the composition of a continuous operator and a continuous compact operator. To apply Leray-Schauder's fixed point Theorem we are left to show that there exists $R$ such that, if $\mathbf h\in \widetilde{Y}$ and $0\leq \lambda\leq 1$ satisfy $\mathbf h=\lambda\mathcal{K}^k(\mathbf h)$, then
 $$
 \Vert \mathbf h\Vert_{\widetilde{Y}}\leq R.
 $$  
 Let us now define $\mathbf{z}_\lambda=(\mathbf{v},\pphi,\ww):=(\mathcal{L}^k)^{-1}\mathbf h$, so that $\mathcal{L}^k(\mathbf{z})=\lambda\mathcal{F}^k(\mathbf{z})$, i.e., $\mathbf{z}_\lambda=(\mathbf{v},\pphi,\ww)$ solves
\begin{alignat}{2}
&\nonumber\qquad\left(2\nu(\pphi)D\mathbf{v},D\boldsymbol\psi\right)_\Omega+\frac\alpha h(D \mathbf{v}, D\boldsymbol\psi)_\Omega\\&\RevA{\nonumber=\int_\Omega\lambda\left[-\frac{\rho\mathbf{v}-\rho^k\mathbf{v}^k}{h}+\frac \alpha h\AA\mathbf{v}^k-\operatorname{div}(\rho^k\mathbf{v}\otimes\mathbf{v})+(\nabla\pphi^k)^T\ww\right.}\\ \nonumber & \left.-\left(\operatorname{div}(\mathbf{J}_\rho)-\frac{\rho-\rho^k}h-\mathbf{v}\cdot\nabla\rho^k\right)\frac {\mathbf{v}}2-(\mathbf{J}_\rho\cdot\nabla)\mathbf{v}\right]\cdot \boldsymbol\psi dx\quad \forall \boldsymbol\psi\in C^\infty_{0,\sigma}(\Omega),\\ \label{phia}
-\operatorname{div}(\mathbf{M}(\pphi^k)\nabla\ww)+\int_\Omega\ww dx&=\lambda\left[-\frac{\pphi-\pphi^k}h-(\nabla\pphi^k)\mathbf{v}+\int_\Omega\ww dx\right],\\ \nonumber
-\Delta\pphi+\mathbf{P}\PsipA(\pphi)&=\lambda\left[\ww+\mathbf{P}\mathbf{A}\frac{\pphi+\pphi^k}{2}\right],\\ \nonumber
\partial_\mathbf{n}\pphi&=0\quad\qquad\qquad\qquad\qquad\qquad\qquad\qquad\qquad\qquad\text{on }\partial\Omega,\\ \nonumber
\partial_\mathbf{n}\ww&=0\quad\qquad\qquad\qquad\qquad\qquad\qquad\qquad\qquad\qquad\text{on }\partial\Omega.
\end{alignat}
 Therefore, by a similar method as for obtaining the energy estimate \eqref{ener}, we infer 
 	\begin{align}
 \nonumber&\lambda\int_\Omega\frac{\rho}{2}\vert \mathbf{v}\vert^2dx+\frac{1}{2}\int_\Omega\vert \nabla\pphi\vert^2-\lambda\int_\Omega\frac{\pphi^T\mathbf{A}\pphi}{2}+\int_\Omega\Psi^1(\pphi)dx +h\int_\Omega \mathbf{M}(\pphi^k)\nabla \ww:\nabla\ww dx\\&\quad+{\alpha}\int_\Omega\vert D\mathbf{v}\vert^2dx+\frac{\alpha\lambda}{2}\int_\Omega\vert D\mathbf{v}- D\mathbf{v}^k\vert^2dx\nonumber+\lambda\int_\Omega \rho^k\frac{\vert \mathbf{v}-\mathbf{v}^k\vert^2}2+(1-\lambda)h\left(\int_\Omega \ww dx\right)^2\\&\quad +2h\int_\Omega \nu(\pphi^k)\vert D\mathbf{v}\vert^2dx
 +\frac 1 2\int_\Omega\vert \nabla\pphi-\nabla \pphi^k\vert^2dx\nonumber\\&\leq \frac{\alpha\lambda}{2}\int_\Omega\vert D\mathbf{v}^k\vert^2dx+\lambda\int_\Omega\frac{\rho^k}{2}\vert \mathbf{v}^k\vert^2dx+\frac{1}{2}\int_\Omega\vert \nabla\pphi^k\vert^2-\lambda\int_\Omega\frac{{(\pphi^k)}^T\mathbf{A}\pphi^k}{2}+\int_\Omega\Psi^1(\pphi^k)dx
 \label{ener2}
 \end{align}
 Note that, clearly, since $\mathbf{z}_\lambda\in {X}$, it holds $\pphi\in[0,1]^N$ and thus $\rho\geq 0$ by Remark \ref{minrho}. \RevA{Also, we have $\left\vert \int_\Omega \Psi^1(\pphi)dx\right\vert\leq C$}. Moreover, $\pphi^k\in[0,1]^N$ and thus $\left\vert \int_\Omega \Psi^1(\pphi^k)dx\right\vert\leq C$. Recall also that $\mathbf{v}^k\in \mathbf{W}_\sigma$. Being $\lambda\in[0,1]$, this entails that we have the following controls
  	\begin{align}
& \nonumber{\alpha}\int_\Omega\vert D\mathbf{v}\vert^2dx+\frac{1}{2h}\int_\Omega\vert \nabla\pphi\vert^2dx+\int_\Omega \mathbf{M}(\pphi^k)\nabla \ww:\nabla\ww dx\\&+(1-\lambda)\left(\int_\Omega \ww dx\right)^2+2\int_\Omega \nu(\pphi^k)\vert D\mathbf{v}\vert^2dx\leq \frac{C_k}{h}.
 \label{ener2bis}
 \end{align}
 Since $h$ is fixed at this level. In the case $\lambda\in\left[0, \beta\right]$, for some $\beta\in(0,1)$ suitably chosen, it is easy to see that it holds %(recall that $\ww\in \mathbf{H}^2_n$)  
 \begin{align*}
 \Vert \ww\Vert_{\widetilde{\mathbf{V}}_0}+\Vert \mathbf{v}\Vert_{\mathbf{V}_\sigma}\leq {C_k(\beta)},
 \end{align*}
 where we include in $C_k$ also the dependence on $h$.
 Now, from \eqref{phia}$_3$, exploiting Theorem \ref{steady} point (3), with $\mathbf f=\lambda\left[\ww+\mathbf{P}\mathbf{A}\frac{\pphi+\pphi^k}{2}\right]\in \widetilde{\mathbf{H}}_0$, we easily see from \eqref{H2b} that 
 $$
 \Vert \pphi\Vert_{\mathbf{H}^2(\Omega)}\leq C_k(\beta).
 $$ 
 In the case $\lambda\in\left (\beta,1\right]$ instead, we observe that, since by \eqref{ener2} $\sqrt{1-\lambda}\left\vert \int_\Omega \ww dx\right\vert dx\leq \frac{C_k}{\sqrt{h}}$, we have
 $$\left\vert \frac{h(\lambda-1)}{\lambda}\int_\Omega \ww dx\right\vert \leq {\sqrt{h}C_k}\frac{\sqrt{1-\beta}}{\beta},\quad \forall \lambda\in(\beta,1],$$ 
 and thus we can fix $\beta=\beta_k$ independent of $\lambda$, such that
 $$ {\sqrt{h}C_k}\frac{\sqrt{1-\beta}}{\beta}\leq \frac 1 2\min_{i=1,\ldots,N}\overline{\pphi^k_i},\quad \forall \lambda\in(\beta,1].$$
 Observe now that, from \eqref{phia}$_2$, by multiplying the equation by $\frac{1}{\vert\Omega\vert_n}$ and then integrating over $\Omega$, we end up, after integration by parts,
  $$\overline{\pphi}=\overline{\pphi}^k+\frac{h(\lambda-1)}{\lambda}\int_\Omega \ww dx.$$ 
  Being $\overline{\pphi}^k\in(0,1)^N$, it thus holds $\overline{\pphi}\in (\tilde{\delta}_k,1-\tilde{\delta}_k),$ for any $\lambda\in(\beta,1]$, with $\tilde{\delta}_k>0$ independent of $\lambda$.
Moreover, from \eqref{ener2bis} we get 
 $$
 \Vert \pphi\Vert_{\mathbf{H}^1(\Omega)}+\Vert \mathbf{v}\Vert_{\mathbf{V}_\sigma}+\Vert \nabla \ww\Vert\leq C_k,
 $$
 where we incorporated the dependence on $h$ in $C_k$.
 In order to control $\overline{\ww}$, we exploit Theorem \ref{steady} point (3) with $\mathbf f=\lambda\left[\ww+\mathbf{P}\mathbf{A}\frac{\pphi+\pphi^k}{2}\right]$ and $\mathbf{m}=\overline{\pphi}$. Indeed, we get, since now $\overline{\pphi}\in(\tilde{\delta}_k,1-\tilde{\delta}_k)$ uniformly in $\lambda$, 
 $$
 \Vert \pphi\Vert_{\mathbf{H}^2(\Omega)}+\vert \overline\ww\vert\leq C_k(1+\lambda\Vert \nabla\ww\Vert+\lambda\Vert \nabla\pphi\Vert+\lambda\Vert \nabla\pphi^k\Vert)\leq C_k,
 $$
 thanks to \eqref{ener2bis} and the fact that $\pphi^k\in [0,1]^N$. This entails that also 
 $$
 \Vert \ww\Vert_{\widetilde{\mathbf{V}}_0}\leq C_k.
 $$
 To sum up, we have proven that there exists a constant $C_k$ such that
   	\begin{align}
\Vert \pphi\Vert_{\mathbf{H}^2(\Omega)}+\Vert \ww\Vert_{\widetilde{\mathbf{V}}_0}+\Vert \mathbf{v}\Vert_{\mathbf{V}_\sigma}\leq C_k \quad \forall \lambda\in[0,1].
\label{p1}
\end{align}
 Furthermore, since $\ww$ solves \eqref{phia}$_2$, by \eqref{ineq} we deduce that, additionally,
 $$
\Vert\ww\Vert_{\mathbf{H}^2(\Omega)}\leq C_k(\Vert \pphi\Vert+\Vert \pphi^k\Vert+\Vert \nabla\pphi^k\Vert_{\mathbf{L}^6(\Omega)}\Vert \mathbf{v}\Vert_{\mathbf{L}^3(\Omega)}+\Vert \ww\Vert)\leq C_k,  
 $$  
 exploiting \eqref{p1}, $\pphi^k\in \mathbf{H}^2(\Omega)$ and the embeddings $\mathbf{H}^2(\Omega)\hookrightarrow\mathbf{W}^{1,6}(\Omega)$ and $\mathbf{V}_\sigma\hookrightarrow \mathbf{L}^3(\Omega)$.
Therefore, we have thus shown that there exists a constant $C_k>0$ such that 
$$
\Vert \mathbf{z}_\lambda\Vert_{\widetilde{X}}\leq C_, \quad \forall {\lambda\in[0,1].}
$$ 
Since then $\mathcal{F}^k$ maps bounded sets of $\widetilde{X}$ into bounded sets of $\widetilde{Y}$, it also holds
$$
\Vert\mathcal{K}^k(\mathbf h)\Vert_{\widetilde{Y}}\leq C_k,
$$ 
so that clearly, since $\lambda\in[0,1]$, there exists $R>0$ such that 
$$
\Vert\mathbf h\Vert_{\widetilde{Y}}=\lambda\Vert\mathcal{K}^k(\mathbf h)\Vert_{\widetilde{Y}}\leq C_k, \quad \forall \RevA{\lambda\in[0,1]},
$$
as needed. This allows to apply Leray-Schauder's Theorem and conclude that there exists a fixed point $\mathbf f\in \widetilde{Y}$ of $\mathcal{K}_k$, so that $\mathbf{z}=(\mathcal{L}^k)^{-1}(\mathbf f)\in X$ is a solution to \eqref{ll}, i.e., a weak solution to \eqref{phi0} satisfying the energy estimate \eqref{ener} and, consequently, \eqref{sepk}. Indeed, as already noticed, since $\ww\in \mathbf{H}_\mathbf{n}^2\hookrightarrow \mathbf{L}^\infty(\Omega)$, this comes directly from Theorem \ref{steaddy} point (2).
\item[\underline{Point (2)}] Under the extra assumption on $\ww_k$, we can actually show a more refined estimate with respect to \eqref{ener}. First recall property \eqref{sepk}, which ensures that, considering $k\in\N$, each $\pphi$ is strictly separated (depending on $k$) from the pure phases. Moreover, also at $k=0$ we have the same result, thanks to Remark \ref{separ}. This means that for any $k\in\N_0$ we can deal with the terms regarding the singular potential $\PsipA$ as if it were a regular one, since we always deal with values belonging to the interior of $(0,1)^N$, where $\PsipA$ is smooth by assumption. Let us then consider \eqref{phi0}$_2$ and test it by $\frac{\ww-\ww^k}{h}\in \mathbf{H}^1(\Omega)$. We get (recall that we assumed $\mathbf{M}$ to be a constant) 
\begin{align}
\left(\frac {\pphi-\pphi^k}h,\frac{\ww-\ww^k}{h}\right)_\Omega+\left((\nabla\pphi^k)\mathbf{v},\frac{\ww-\ww^k}{h}\right)_\Omega+\left(\mathbf{M}\nabla\ww,\frac{\nabla\ww-\nabla\ww^k}{h}\right)_\Omega=0.
\label{wd}
\end{align} 
Moreover, we also observe that we have, having set $\pphi_{-1}=\pphi^0$ if $k=0$,
$$
\frac{\ww-\ww^k}{h}=-\frac{\Delta\pphi-\Delta\pphi^k}{h}+\mathbf{P}\frac{\PsipA(\pphi)-\PsipA(\pphi^k)}{h}-\mathbf{P}\frac{\mathbf{A}(\pphi-\pphi^{k-1})}{h}.
$$
Observe now that, by the identity \eqref{identity1},
\begin{align*}
&\left(\mathbf{M}\nabla\ww,\frac{\nabla\ww-\nabla\ww^k}{h}\right)_\Omega\\&=\frac 1 {2h} \left(\mathbf{M}\nabla \ww,\nabla\ww\right)_\Omega-\frac 1 {2h} \left(\mathbf{M}\nabla \ww^k,\nabla\ww^k\right)_\Omega+\frac 1 {2h} \left(\mathbf{M}(\nabla\ww-\nabla\ww^k),\nabla\ww-\nabla\ww^k\right)_\Omega.
\end{align*}
Then we have, since $\pphi-\pphi^k\in T\Sigma$ for any $x\in\Omega$,
\begin{align*}
&\left(\frac {\pphi-\pphi^k}h,\frac{\ww-\ww^k}{h}\right)_\Omega\\&=\left\Vert \frac{\nabla\pphi-\nabla\pphi^k}{h}\right\Vert^2+\left(\frac {\pphi-\pphi^k}h,\frac{\PsipA(\pphi)-\PsipA(\pphi^k)}{h}\right)_\Omega-\left(\frac {\pphi-\pphi^k}h,\frac{\mathbf{A}(\pphi-\pphi^{k-1})}{h}\right)_\Omega.
\end{align*}
Futhermore, by assumption (\textbf{E0}),
\begin{align*}
&\left(\frac {\pphi-\pphi^k}h,\frac{\PsipA(\pphi)-\PsipA(\pphi^k)}{h}\right)_\Omega\\&=\frac 1 {h^2}\sum_{i=1}^N\int_\Omega \left(\vert \varphi_i-\varphi^k_i\vert^2\int_0^1\left(\psi^{\prime\prime}(s\varphi_i+(1-s)\varphi_i^k)\right)ds\right)dx\geq 0.
\end{align*}
We now introduce the operator $L: \mathbf{V}_0\to \mathbf{V}_0'$ such that, for any $\mathbf{f}\in \mathbf{V}_0$,
$$
<L\mathbf{f},\boldsymbol\psi>=\int_\Omega \nabla\mathbf{f}:\nabla\boldsymbol\psi dx\quad \forall \boldsymbol\psi\in \mathbf{V}_0,
$$
which is invertible by the Lax-Milgram Lemma thanks to Poincar\'{e}'s inequality. Since $\pphi-\pphi^{k-1}\in \mathbf{V}_0$ and $\pphi-\pphi^k\in \mathbf{V}_0$, we have
\begin{align*}
-\left(\frac {\pphi-\pphi^k}h,\frac{\mathbf{A}(\pphi-\pphi^{k-1})}{h}\right)_\Omega&=-\left(\frac {\pphi-\pphi^k}h,\frac{\mathbf{A}(\pphi-\pphi^{k})}{h}\right)_\Omega-\left(\frac {\pphi-\pphi^k}h,\frac{\mathbf{A}(\pphi^k-\pphi^{k-1})}{h}\right)_\Omega\\&=-\left(\nabla \left(\frac {\pphi-\pphi^k}h\right),\nabla L^{-1}\left(\frac{\mathbf{A}(\pphi-\pphi^{k})}{h}\right)\right)_\Omega\\&\quad -\left(\nabla \left(\frac {\pphi-\pphi^k}h\right),\nabla L^{-1}\left(\frac{\mathbf{A}(\pphi^k-\pphi^{k-1})}{h}\right)\right)_\Omega.
\end{align*}
Recalling that $\Vert\nabla L^{-1}(\cdot)\Vert$ is an equivalent norm on $\mathbf{V}_0'$, by Cauchy-Schwarz and Young's inequalities we get 
$$
-\left(\frac {\pphi-\pphi^k}h,\frac{\mathbf{A}(\pphi-\pphi^{k-1})}{h}\right)_\Omega\leq \frac 1 8 \left\Vert \frac{\nabla\pphi-\nabla\pphi^k}{h}\right\Vert^2+C\left(\left\Vert \frac{\pphi-\pphi^k}{h}\right\Vert^2_{\mathbf{V}_0'}+\left\Vert \frac{\pphi^k-\pphi^{k-1}}{h}\right\Vert^2_{\mathbf{V}_0'}\right).
$$
Clearly, in the case $k=0$, the third term in the right-hand side does not appear. We are then left with one last term:
\begin{align*}
\left((\nabla\pphi^k)\mathbf{v},\frac{\ww-\ww^k}{h}\right)_\Omega=\left((\nabla\pphi^k)\mathbf{v},-\frac{\Delta\pphi-\Delta\pphi^k}{h}+\frac{\PsipA(\pphi)-\PsipA(\pphi^k)}{h}-\frac{\mathbf{A}(\pphi-\pphi^{k-1})}{h}\right)_\Omega.
\end{align*}
Note that the projector $\mathbf{P}$ does not come into play, since $(\nabla\pphi^k)\mathbf{v}\in T\Sigma$ almost everywhere in $\Omega$. This can be immediately seen from equation \eqref{phi0}$_2$, but also directly, since $(\nabla\pphi^k)\mathbf{v}\cdot \mathbf{e}=\mathbf{v}\cdot (\nabla\pphi^k)^T\mathbf{e}\equiv \mathbf{0}$ and due to $\sum_{i=1}^N\pphi^k_i\equiv1$.
Now, integrating by parts, recalling the embeddings, for $n=2,3$, $\mathbf{H}^1(\Omega)\hookrightarrow \mathbf{L}^4(\Omega)$, $\mathbf{W}^{2,4}(\Omega)\hookrightarrow \mathbf{W}^{1,\infty}(\Omega)$,
\begin{align*}
&\left((\nabla\pphi^k)\mathbf{v},-\frac{\Delta\pphi-\Delta\pphi^k}{h}\right)_\Omega\\&=\left((D^2\pphi^k)\cdot \mathbf{v},\frac{\nabla\pphi-\nabla\pphi^k}{h}\right)_\Omega+\left(\nabla\pphi^k\nabla\mathbf{v},\dfrac{\nabla\pphi-\nabla\pphi^k}{h}\right)_\Omega\\&\leq \Vert\pphi^k\Vert_{\mathbf{W}^{2,4}(\Omega)}\Vert \mathbf{v}\Vert_{\mathbf{L}^4(\Omega)}\left\Vert \frac{\nabla\pphi-\nabla\pphi^k}{h}\right\Vert+\Vert\pphi^k\Vert_{\mathbf{W}^{1,\infty}(\Omega)}\Vert \nabla\mathbf{v}\Vert\left\Vert \frac{\nabla\pphi-\nabla\pphi^k}{h}\right\Vert\\&\leq \frac 1 8\left\Vert \frac{\nabla\pphi-\nabla\pphi^k}{h}\right\Vert^2+C\Vert\nabla\mathbf{v}\Vert^2\Vert\pphi^k\Vert_{\mathbf{W}^{2,4}(\Omega)}^2.
\end{align*}
Then, similarly as before, by means of the operator $L$ we have by H\"{o}lder's and Young's inequalities,
\begin{align*}
&\left((\nabla\pphi^k)\mathbf{v},-\frac{\mathbf{A}(\pphi-\pphi^{k-1})}{h}\right)_\Omega\leq C\Vert\pphi^k\Vert_{\mathbf{W}^{1,4}(\Omega)}\Vert \mathbf{v}\Vert_{\mathbf{L}^4(\Omega)}\left(\left\Vert \frac{\pphi-\pphi^k}{h}\right\Vert+\left\Vert \frac{\pphi^k-\pphi^{k-1}}{h}\right\Vert\right)\\& \leq C\Vert\pphi^k\Vert_{\mathbf{W}^{1,4}(\Omega)}\Vert \mathbf{v}\Vert_{\mathbf{L}^4(\Omega)}\\&\RevA{\quad\times\left(\left\Vert \frac{\nabla\pphi-\nabla\pphi^k}{h}\right\Vert^{\frac 1 2}\left\Vert \frac{\pphi-\pphi^k}{h}\right\Vert^{\frac 1 2}_{\mathbf{V}_0'}+\left\Vert \frac{\nabla\pphi^k-\nabla\pphi^{k-1}}{h}\right\Vert^{\frac 1 2}\left\Vert \frac{\pphi^k-\pphi^{k-1}}{h}\right\Vert^{\frac 1 2}_{\mathbf{V}_0'}\right)}\\&\leq \frac 1 8\left\Vert \frac{\nabla\pphi-\nabla\pphi^k}{h}\right\Vert^2+C\left(\Vert\nabla\mathbf{v}\Vert^2\Vert\pphi^k\Vert_{\mathbf{W}^{1,4}(\Omega)}^2+\left\Vert \frac{\pphi-\pphi^k}{h}\right\Vert^{2}_{\mathbf{V}_0'}+\left\Vert \frac{\pphi^k-\pphi^{k-1}}{h}\right\Vert^2_{\mathbf{V}_0'}\right).
\end{align*}
In conclusion, by integration by parts and standard inequalities, recalling that $\operatorname{div}\mathbf{v}=0$,
\begin{align}
&\nonumber\left((\nabla\pphi^k)\mathbf{v},\frac{\PsipA(\pphi)-\PsipA(\pphi^k)}{h}\right)_\Omega\\&\RevA{=\nonumber\sum_{i=1}^N\frac 1 h\int_\Omega\left( \nabla\varphi^k_{i}\cdot\mathbf{v}\int_0^1 \psi''(s\varphi_i+(1-s)\varphi^k_i)(\varphi_i-\varphi^k_i)ds\right)dx}\\&\nonumber=\sum_{i=1}^N\frac 1 h\int_\Omega\left(\int_0^1 (s\nabla\varphi_{i}+(1-s)\nabla\varphi_i^k)\cdot\mathbf{v} \psi''(s\varphi_i+(1-s)\varphi^k_i)(\varphi_i-\varphi^k_i)ds\right)dx\\&\nonumber-\sum_{i=1}^N\frac 1 h\int_\Omega\left(\int_0^1 s(\nabla\varphi_i-\nabla\varphi^k_i)\cdot\mathbf{v} \left(\partial_s\psi'(s\varphi_i+(1-s)\varphi^k_i)\right)ds\right)dx\\&\nonumber
=\sum_{i=1}^N\frac 1 h\int_\Omega\left( \int_0^1\mathbf{v}\cdot \left(\nabla\psi'(s\varphi_i+(1-s)\varphi^k_i)\right)(\varphi_i-\varphi^k_i)ds\right)dx\\&\nonumber
-\sum_{i=1}^N\frac{1}h\left[\int_\Omega s(\nabla\varphi_i-\nabla\varphi^k_i)\cdot\mathbf{v}\psi'(s\varphi_i+(1-s)\varphi^k_i)dx\right]_{0}^1\\&\nonumber+\sum_{i=1}^N\frac{1}h\int_\Omega\left(\int_0^1 (\nabla\varphi_i-\nabla\varphi^k_i)\cdot\mathbf{v}\psi'(s\varphi_i+(1-s)\varphi^k_i)ds\right)dx\\&\nonumber
=\sum_{i=1}^N\frac 1 h\int_\Omega\left( \int_0^1\mathbf{v}\cdot \left(\nabla\psi'(s\varphi_i+(1-s)\varphi^k_i)\right)(\varphi_i-\varphi^k_i)ds\right)dx\\&\nonumber
-\sum_{i=1}^N\frac{1}h\int_\Omega (\nabla\varphi_i-\nabla\varphi^k_i)\cdot\mathbf{v}\psi'(\varphi_i)dx\nonumber\\&\nonumber\RevA{-\sum_{i=1}^N\frac 1 h\int_\Omega\left( \int_0^1\mathbf{v}\cdot \left(\nabla\psi'(s\varphi_i+(1-s)\varphi^k_i)\right)(\varphi_i-\varphi^k_i)ds\right)dx}\\&\nonumber=-\sum_{i=1}^N\frac{1}h\int_\Omega (\nabla\varphi_i-\nabla\varphi^k_i)\cdot\mathbf{v}\psi'(\varphi_i)dx\leq C \left\Vert \frac{\nabla\pphi-\nabla\pphi^k}{h}\right\Vert\Vert \mathbf{v}\Vert_{\mathbf{L}^3(\Omega)}\Vert \PsipA(\pphi)\Vert_{\mathbf{L}^6(\Omega)}\\&\leq C\left\Vert\frac{\nabla\pphi-\nabla\pphi^k}{h}\right\Vert\Vert \nabla\mathbf{v}\Vert^\frac 1 2\Vert \mathbf{v}\Vert^\frac 1 2\Vert \PsipA(\pphi)\Vert_{\mathbf{L}^6(\Omega)}\nonumber\\&\RevA{\leq \frac 1 8 \left\Vert\frac{\nabla\pphi-\nabla\pphi^k}{h}\right\Vert^2+C\Vert \nabla\mathbf{v}\Vert\Vert \mathbf{v}\Vert\Vert \PsipA(\pphi)\Vert_{\mathbf{L}^6(\Omega)}^2
},
\label{parts}
\end{align}
where we exploited also Sobolev-Gagliardo-Nirenberg's inequalities in 3D. It is then clear that in the 2D case it is enough to use Sobolev-Gagliardo-Nirenberg's inequalities for the $\mathbf{L}^4(\Omega)$ norm, leading to an analogous result. Now we observe that, since \eqref{phi0}$_3$ holds and by the assumption on $\ww_k$, recalling that it holds $\overline{\pphi}=\overline{\pphi}^k=\overline{\pphi}^0$, we have by Theorem \ref{steaddy} point (3), with $\mathbf f=\mathbf{P}(\ww+\mathbf{A}\frac{\pphi-\pphi^k}{2})$ and $\mathbf{m}=\overline{\pphi}^0$, that 
$$
\Vert \PsipA(\pphi)\Vert_{\mathbf{L}^6(\Omega)}\leq C(\min_{i=1,\ldots,N}\overline{\varphi}^0_i)(1+\Vert\nabla\ww\Vert),\quad \Vert \pphi^k\Vert_{\mathbf{W}^{2,4}(\Omega)}\leq C(\min_{i=1,\ldots,N}\overline{\varphi}^0_i)(1+\Vert\nabla\ww^k\Vert).
$$
Putting everything together in \eqref{wd} we end up with 
\begin{align*}
&\frac 1 {2h} \left(\mathbf{M}\nabla \ww,\nabla\ww\right)+\frac 1 2 \left\Vert\frac{\nabla\pphi-\nabla\pphi^k}{h}\right\Vert^2\\&\leq \frac 1 {2h} \left(\mathbf{M}\nabla \ww^k,\nabla\ww^k\right)+C(\min_{i=1,\ldots,N}\overline{\varphi}^0_i)\Vert \nabla\mathbf{v}\Vert\Vert \mathbf{v}\Vert(1+\Vert \nabla\ww^k\Vert^2+\Vert \nabla\ww\Vert^2)\\&+C(\min_{i=1,\ldots,N}\overline{\varphi}^0_i)\Vert \nabla\mathbf{v}\Vert^2(1+\Vert \nabla\ww^k\Vert^2)+C\left(\left\Vert \frac{\pphi-\pphi^k}{h}\right\Vert^{2}_{\mathbf{V}_0'}+\left\Vert \frac{\pphi^k-\pphi^{k-1}}{h}\right\Vert^2_{\mathbf{V}_0'}\right),
\end{align*}
 which is exactly \eqref{h}.
 \item[\underline{Point (3)}]
 
 Assume now that we have $\alpha>0$ and $\mathbf{v}^k\in \mathbf{W}_\sigma$, together with the other assumptions at point (2). We can thus test \eqref{phi0}$_1$ by $\boldsymbol\psi=\frac{\mathbf{\mathbf{v}}-\mathbf{v}^k}{h}\in \mathbf{V}_\sigma$, to get 
 \begin{align}
 &\nonumber\left(\dfrac{\rho \mathbf{v}-\rho^k\mathbf{v}^k}{h},\frac{\mathbf{v}-\mathbf{v}^k}{h}\right)_\Omega+\alpha\left\Vert\frac{D\mathbf{v}-D\mathbf{v}^k}{h}\right\Vert^2+\left(\operatorname{div}(\rho^k\mathbf{v}\otimes\mathbf{v}),\frac{\mathbf{v}-\mathbf{v}^k}{h} \right)_\Omega\\&+\left(\operatorname{div}(\mathbf{v}\otimes\mathbf{J}_\rho),\frac{\mathbf{v}-\mathbf{v}^k}{h}\right)_\Omega+\left(2\nu(\pphi^k)D\mathbf{v},D\frac{\mathbf{v}-\mathbf{v}^k}{h}\right)_\Omega=\left((\nabla\pphi^k)^T\ww,\frac{\mathbf{v}-\mathbf{v}^k}{h}\right)_\Omega.
 \label{o}
 \end{align}    
 Now we have 
 \begin{align}
 \left(\dfrac{\rho \mathbf{v}-\rho^k\mathbf{v}^k}{h},\frac{\mathbf{v}-\mathbf{v}^k}{h}\right)_\Omega=\left(\rho^k\dfrac{ \mathbf{v}-\mathbf{v}^k}{h},\frac{\mathbf{v}-\mathbf{v}^k}{h}\right)_\Omega+\left(\dfrac{ \rho-\rho^k}{h}\mathbf{v},\frac{\mathbf{v}-\mathbf{v}^k}{h}\right)_\Omega.
 \label{dec}
 \end{align}
 Recalling Remark \ref{minrho}, we have
 $$
 \left(\rho^k\dfrac{ \mathbf{v}-\mathbf{v}^k}{h},\frac{\mathbf{v}-\mathbf{v}^k}{h}\right)_\Omega\geq C_\rho\left\Vert\dfrac{ \mathbf{v}-\mathbf{v}^k}{h}\right\Vert^2,
 $$
 where $C_\rho=\min_{i=1,\ldots,N}\widetilde{\rho}_i>0$.
 Then we have, recalling \eqref{eqro}, 
 \begin{align*}
 &\left(\operatorname{div}(\rho^k\mathbf{v}\otimes\mathbf{v}),\frac{\mathbf{v}-\mathbf{v}^k}{h} \right)_\Omega+\left(\operatorname{div}(\mathbf{v}\otimes\mathbf{J}_\rho),\frac{\mathbf{v}-\mathbf{v}^k}{h}\right)_\Omega+\left(\dfrac{ \rho-\rho^k}{h}\mathbf{v},\frac{\mathbf{v}-\mathbf{v}^k}{h}\right)_\Omega\\&=\left(\rho^k(\mathbf{v}\cdot\nabla)\mathbf{v},\frac{\mathbf{v}-\mathbf{v}^k}{h}\right)_\Omega+\left((\mathbf{J}_\rho\cdot\nabla)\mathbf{v},\frac{\mathbf{v}-\mathbf{v}^k}{h} \right)_\Omega. 
 \end{align*}
 Now, by H\"{o}lder's,  Poincar\'{e}'s and Young's inequalities, we get, recalling that $\rho^k$ is bounded uniformly in $k$ due to $\pphi^k\in[0,1]^N$,
 \begin{align*}
 &\left\vert\left(\rho^k(\mathbf{v}\cdot\nabla)\mathbf{v},\frac{\mathbf{v}-\mathbf{v}^k}{h}\right)_\Omega\right\vert\leq C\Vert \mathbf{v}\Vert_{\mathbf{L}^4(\Omega)}\Vert \nabla\mathbf{v}\Vert\left\Vert\dfrac{ \mathbf{v}-\mathbf{v}^k}{h}\right\Vert_{\mathbf{L}^4(\Omega)}\leq \frac{\alpha}8\left\Vert\dfrac{D\mathbf{v}-D\mathbf{v}^k}{h}\right\Vert^2+\frac{C}{\alpha}\Vert\nabla\mathbf{v}\Vert^4.
 \end{align*} 
 Then, concerning the other terms, recalling the definition of $\mathbf{J}_\rho$, we get
 \begin{align*}
 &\left\vert\left((\mathbf{J}_\rho\cdot\nabla)\mathbf{v},\frac{\mathbf{v}-\mathbf{v}^k}{h} \right)_\Omega\right\vert\leq C\Vert \nabla\ww\Vert_{\mathbf{L}^4(\Omega)}\Vert \nabla\mathbf{v}\Vert\left\Vert\dfrac{ \mathbf{v}-\mathbf{v}^k}{h}\right\Vert_{\mathbf{L}^4(\Omega)}\\&\leq C\Vert \ww\Vert_{\mathbf{H}^2(\Omega)}\Vert \nabla\mathbf{v}\Vert\left\Vert\dfrac{D\mathbf{v}-D\mathbf{v}^k}{h}\right\Vert\leq\frac {\alpha} 8\left\Vert\dfrac{D\mathbf{v}-D\mathbf{v}^k}{h}\right\Vert^2+\frac{C}{\alpha}\Vert\ww\Vert^2_{\mathbf{H}^2(\Omega)}\Vert \nabla \mathbf{v}\Vert^2.
 \end{align*}
 We then continue the analysis of \eqref{o}: we have, recalling that $\nu(\pphi^k)$ is bounded uniformly in $k$,
 \begin{align*}
 &\left(2\nu(\pphi^k)D\mathbf{v},D\frac{\mathbf{v}-\mathbf{v}^k}{h}\right)_\Omega\leq \frac{C}{\alpha}\Vert \nabla\mathbf{v}\Vert^2+ \frac {\alpha} 8\left\Vert\dfrac{D\mathbf{v}-D\mathbf{v}^k}{h}\right\Vert^2.
 \end{align*}
 We are only left with one term: standard inequalities entail
 \begin{align*}
 &\left\vert\left((\nabla\pphi^k)^T\ww,\frac{\mathbf{v}-\mathbf{v}^k}{h}\right)_\Omega\right\vert\leq \Vert \nabla\pphi^k\Vert\Vert\ww\Vert_{\mathbf{L}^4(\Omega)}\left\Vert\dfrac{ \mathbf{v}-\mathbf{v}^k}{h}\right\Vert_{\mathbf{L}^4(\Omega)}\\&\leq C\Vert \nabla\pphi^k\Vert\Vert\ww\Vert_{\mathbf{H}^1(\Omega)}\left\Vert\dfrac{D\mathbf{v}-D\mathbf{v}^k}{h}\right\Vert\leq \frac {\alpha}{8}\left\Vert\dfrac{D\mathbf{v}-D\mathbf{v}^k}{h}\right\Vert^2+\frac C\alpha\Vert \mathbf{w}\Vert_{\mathbf{H}^1(\Omega)}^2\Vert \nabla\pphi^k\Vert^2.
 \end{align*}
 To sum up, coming back to \eqref{o}, in the end we obtain 
 \begin{align*}
 \alpha\left\Vert\dfrac{D\mathbf{v}-D\mathbf{v}^k}{h}\right\Vert^2\leq \frac{C}{\alpha}\left(\Vert \mathbf{w}\Vert_{\mathbf{H}^1(\Omega)}^2\Vert \nabla\pphi^k\Vert^2+\Vert\ww\Vert^2_{\mathbf{H}^2(\Omega)}\Vert \nabla \mathbf{v}\Vert^2+\Vert\nabla\mathbf{v}\Vert^4+\Vert\nabla\mathbf{v}\Vert^2\right),
 \end{align*}
 which is exactly \eqref{c1}.

 Recalling Remark \ref{H3}, \eqref{phi0}$_1$ also holds in almost everywhere sense in $\Omega$. We can thus test \eqref{phi0}$_1$ by $\boldsymbol\psi=\AA\mathbf{v}\in \mathbf{H}_\sigma$, to get 
 \begin{align}
 &\nonumber\left(\dfrac{\rho \mathbf{v}-\rho^k\mathbf{v}^k}{h},\AA\mathbf{v}\right)_\Omega+\alpha\left(\dfrac{\AA\mathbf{v}-\AA\mathbf{v}^k}{h},\AA\mathbf{v}\right)_\Omega+\left(\operatorname{div}(\rho^k\mathbf{v}\otimes\mathbf{v}),\AA\mathbf{v} \right)_\Omega+\left(\operatorname{div}(\mathbf{v}\otimes\mathbf{J}_\rho),\AA\mathbf{v}\right)_\Omega\\&-\left(2\operatorname{div}(\nu(\pphi^k)D\mathbf{v}),\AA\mathbf{v}\right)_\Omega=\left((\nabla\pphi^k)^T\ww,\AA\mathbf{v}\right)_\Omega.
 \label{o2}
 \end{align}    
 Now we have 
 \begin{align}
 \left(\dfrac{\rho \mathbf{v}-\rho^k\mathbf{v}^k}{h},\AA\mathbf{v}\right)_\Omega=\left(\rho^k\dfrac{ \mathbf{v}-\mathbf{v}^k}{h},\AA\mathbf{v}\right)_\Omega+\left(\dfrac{ \rho-\rho^k}{h}\mathbf{v},\AA\mathbf{v}\right)_\Omega.
 \label{dec3}
 \end{align}
 Then we have, by the Poincar\'{e}, Young and Korn inequalities,
 $$
\left\vert \left(\rho^k\dfrac{ \mathbf{v}-\mathbf{v}^k}{h},\AA\mathbf{v}\right)_\Omega\right\vert\leq C\left\Vert\dfrac{ D\mathbf{v}-D\mathbf{v}^k}{h}\right\Vert^2+\frac{1}{2}\Vert \AA\mathbf{v}\Vert^2,
 $$
 where $C>0$ clearly depends on $\max_{i=1,\ldots,N}\widetilde{\rho}_i$.
 Then, recalling \eqref{Jrho}, 
 \begin{align*}
 &\left(\operatorname{div}(\rho^k\mathbf{v}\otimes\mathbf{v}),\AA\mathbf{v} \right)_\Omega+\left(\operatorname{div}(\mathbf{v}\otimes\mathbf{J}_\rho),\AA\mathbf{v}\right)_\Omega+\left(\dfrac{ \rho-\rho^k}{h}\mathbf{v},\AA\mathbf{v}\right)_\Omega\\&=\left(\rho^k(\mathbf{v}\cdot\nabla)\mathbf{v},\AA\mathbf{v}\right)_\Omega+\left((\mathbf{J}_\rho\cdot\nabla)\mathbf{v},\AA\mathbf{v} \right)_\Omega. 
 \end{align*}
 Now, by H\"{o}lder's,  Poincar\'{e}'s, Young's and Agmon's inequalities (in 3D and analogously in 2D), we get, recalling that $\rho^k$ is bounded uniformly in $k$ due to $\pphi^k\in[0,1]^N$,
 \begin{align*}
 &\left\vert\left(\rho^k(\mathbf{v}\cdot\nabla)\mathbf{v},\AA\mathbf{v}\right)_\Omega\right\vert\leq C\Vert \mathbf{v}\Vert_{\mathbf{L}^\infty(\Omega)}\Vert \nabla\mathbf{v}\Vert\left\Vert\AA\mathbf{v}\right\Vert\\&\leq C\Vert \nabla\mathbf{v}\Vert^{\frac 3 2}\Vert\left\Vert\AA\mathbf{v}\right\Vert^\frac 3 2\leq \frac{1}2\left\Vert\AA\mathbf{v}\right\Vert^2+C\Vert\nabla\mathbf{v}\Vert^6.
 \end{align*} 
 Then, concerning the other terms, recalling the definition of $\mathbf{J}_\rho$, we get in 3D (the 2D case is analogous, using $\mathbf{L}^\infty(\Omega)$ norm instead of $\mathbf{L}^6(\Omega)$, together with Agmon's inequality)
 \begin{align*}
 &\left\vert\left((\mathbf{J}_\rho\cdot\nabla)\mathbf{v},\AA\mathbf{v} \right)_\Omega\right\vert\leq C\Vert \nabla\ww\Vert_{\mathbf{L}^6(\Omega)}\Vert \nabla\mathbf{v}\Vert_{\mathbf{L}^3(\Omega)}\Vert \AA\mathbf{v}\Vert\\&\leq C\Vert \ww\Vert_{\mathbf{H}^3(\Omega)}^\frac 1 2 \Vert \nabla\ww\Vert^\frac 1 2\Vert \nabla\mathbf{v}\Vert^\frac 1 2\Vert \AA\mathbf{v}\Vert^\frac 3 2 \\&\leq\frac {1} 2\left\Vert\AA\mathbf{v}\right\Vert^2+C\Vert \nabla\ww\Vert^2\Vert\ww\Vert^2_{\mathbf{H}^3(\Omega)}\Vert \nabla \mathbf{v}\Vert^2,
 \end{align*}
  where we used the 3D inequality
 \begin{align}
 \Vert \mathbf f\Vert_{\mathbf{L}^6(\Omega)}\leq C\Vert \mathbf f\Vert^{\frac 1 2}\Vert \mathbf f\Vert^{\frac 1 2}_{\mathbf{H}^2(\Omega)}.
 \label{pp}
 \end{align}
 We then continue the analysis of \eqref{o2}: in 3D (and analogously in 2D, by using $\mathbf{L}^4(\Omega)$ norm instead) we have, recalling that $\nu(\pphi^k)$ and $\nu_{,\pphi}$  are bounded uniformly in $k$,
 \begin{align}
 \nonumber&\left\vert\left(2\operatorname{div}(\nu(\pphi^k)D\mathbf{v}),\AA\mathbf{v}\right)_\Omega\right\vert\leq\left\vert\left(2\nu(\pphi^k)\Delta\mathbf{v},\AA\mathbf{v}\right)_\Omega\right\vert+\left\vert\left(2 D\mathbf{v}(\nabla\pphi^k)^T\nu_{,\pphi}(\pphi^k),\AA\mathbf{v}\right)_\Omega\right\vert
\\&\leq  C\Vert \AA\mathbf{v}\Vert^2+\Vert \nabla\pphi^k\Vert_{\mathbf{L}^6(\Omega)}\Vert D\mathbf{v}\Vert_{\mathbf{L}^3(\Omega)}\Vert \AA\mathbf{v}\Vert\nonumber\\&
 \leq C\Vert \AA\mathbf{v}\Vert^2+\Vert \nabla\pphi^k\Vert_{\mathbf{L}^6(\Omega)}\Vert \nabla\mathbf{v}\Vert^{\frac 1 2}\Vert \AA\mathbf{v}\Vert^\frac 3 2\nonumber
\\&\leq C\Vert \AA\mathbf{v}\Vert^2+C\Vert \nabla\pphi^k\Vert_{\mathbf{L}^6(\Omega)}^4\Vert \nabla\mathbf{v}\Vert^2.\label{div}
 \end{align}
 We are only left with one term: standard inequalities entail
 \begin{align*}
 &\left\vert\left((\nabla\pphi^k)^T\ww,\AA\mathbf{v}\right)_\Omega\right\vert\leq \Vert \nabla\pphi^k\Vert\Vert\ww\Vert_{\mathbf{L}^\infty(\Omega)}\left\Vert\AA\mathbf{v}\right\Vert\\&\leq C\Vert \nabla\pphi^k\Vert\Vert\ww\Vert_{\mathbf{H}^2(\Omega)}\left\Vert\AA\mathbf{v}\right\Vert\leq C\left\Vert\AA\mathbf{v}\right\Vert^2+ C\Vert \mathbf{w}\Vert_{\mathbf{H}^2(\Omega)}^2\Vert \nabla\pphi^k\Vert^2.
 \end{align*}
 In conclusion, as usual we notice that
 $$
 \alpha\left(\frac{\AA\mathbf{v}-\AA\mathbf{v}^k}{h},\AA\mathbf{v}\right)_\Omega=\frac{\alpha}{2h}\Vert \AA\mathbf{v}\Vert^2-\frac{\alpha}{2h}\Vert \AA\mathbf{v}^k\Vert^2+\frac{\alpha}{2h}\Vert \AA\mathbf{v}-\AA\mathbf{v}^k\Vert^2.
 $$
 To sum up, coming back to \eqref{o2}, in the end we obtain 
 \begin{align*}
 &\frac\alpha{2h}\left\Vert\AA\mathbf{v}\right\Vert^2\leq \frac\alpha{2h}\Vert\AA\mathbf{v}^k\Vert^2+C\left\Vert\AA\mathbf{v}\right\Vert^2+ C\left\Vert\dfrac{ D\mathbf{v}-D\mathbf{v}^k}{h}\right\Vert^2\\&+C\left(\Vert \nabla\mathbf{v}\Vert^6+\Vert \nabla\ww\Vert^2\Vert\ww\Vert^2_{\mathbf{H}^2(\Omega)}\Vert \nabla \mathbf{v}\Vert^2+\Vert \nabla\pphi^k\Vert_{\mathbf{L}^6(\Omega)}^4\Vert \nabla\mathbf{v}\Vert^2+\Vert \mathbf{w}\Vert_{\mathbf{H}^3(\Omega)}^2\Vert \nabla\pphi^k\Vert^2\right)_\Omega,
 \end{align*}
 which is exactly \eqref{c2}.

Furthermore, we multiply \eqref{phi0}$_1$ (written in almost everywhere formulation) by $\frac{\AA\mathbf{v}-\AA\mathbf{v}^{k}}{h}\in \mathbf{H}_\sigma$, to get, after integrating over $\Omega$,
 \begin{align}
 &\nonumber\left(\dfrac{\rho \mathbf{v}-\rho^k\mathbf{v}^k}{h},\frac{\AA\mathbf{v}-\AA\mathbf{v}^k}{h}\right)_\Omega+\alpha\left\Vert\frac{\AA\mathbf{v}-\AA\mathbf{v}^k}{h}\right\Vert^2+\left(\operatorname{div}(\rho^k\mathbf{v}\otimes\mathbf{v}),\frac{\AA\mathbf{v}-\AA\mathbf{v}^k}{h} \right)_\Omega\\&+\left(\operatorname{div}(\mathbf{v}\otimes\mathbf{J}_\rho),\frac{\AA\mathbf{v}-\AA\mathbf{v}^k}{h}\right)_\Omega-\left(2\operatorname{div}(\nu(\pphi^k)D\mathbf{v}),\frac{\AA\mathbf{v}-\AA\mathbf{v}^k}{h}\right)_\Omega\nonumber\\&=\left((\nabla\pphi^k)^T\ww,\frac{\AA\mathbf{v}-\AA\mathbf{v}^k}{h}\right)_\Omega.
 \label{o3}
 \end{align}    
 Now we have 
 \begin{align}
 \left(\dfrac{\rho \mathbf{v}-\rho^k\mathbf{v}^k}{h},\frac{\AA\mathbf{v}-\AA\mathbf{v}^k}{h}\right)_\Omega=\left(\rho^k\dfrac{ \mathbf{v}-\mathbf{v}^k}{h},\frac{\AA\mathbf{v}-\AA\mathbf{v}^k}{h}\right)_\Omega+\left(\dfrac{ \rho-\rho^k}{h}\mathbf{v},\frac{\AA\mathbf{v}-\AA\mathbf{v}^k}{h}\right)_\Omega.
 \label{dec2}
 \end{align} 
 We thus have, since $\rho^k$ is bounded, by Korn's inequality,
$$
\left\vert\left(\rho^k\dfrac{ \mathbf{v}-\mathbf{v}^k}{h},\frac{\AA\mathbf{v}-\AA\mathbf{v}^k}{h}\right)_\Omega\right\vert\leq C\left\Vert\dfrac{\mathbf{v}-\mathbf{v}^k}{h}\right\Vert\left\Vert\frac{\AA\mathbf{v}-\AA\mathbf{v}^k}{h}\right\Vert\leq \frac \alpha 8 \left\Vert\frac{\AA\mathbf{v}-\AA\mathbf{v}^k}{h}\right\Vert^2+\frac{C}\alpha\left\Vert\dfrac{D\mathbf{v}-D\mathbf{v}^k}{h}\right\Vert^2.
$$
 Then we have, recalling \eqref{eqro}, 
 \begin{align*}
&\left(\operatorname{div}(\rho^k\mathbf{v}\otimes\mathbf{v}),\dfrac{\AA\mathbf{v}-\AA\mathbf{v}^k}{h} \right)_\Omega+\left(\operatorname{div}(\mathbf{v}\otimes\mathbf{J}_\rho),\dfrac{\AA\mathbf{v}-\AA\mathbf{v}^k}{h}\right)_\Omega+\left(\dfrac{ \rho-\rho^k}{h}\mathbf{v},\dfrac{\AA\mathbf{v}-\AA\mathbf{v}^k}{h}\right)_\Omega\\&=\left(\rho^k(\mathbf{v}\cdot\nabla)\mathbf{v},\dfrac{\AA\mathbf{v}-\AA\mathbf{v}^k}{h}\right)_\Omega+\left((\mathbf{J}_\rho\cdot\nabla)\mathbf{v},\dfrac{\AA\mathbf{v}-\AA\mathbf{v}^k}{h} \right)_\Omega. 
 \end{align*}
Now, by H\"{o}lder's,  Poincar\'{e}'s and Young's inequalities, together with 3D Sobolev-Gagliardo-Nirenberg's inequality (and analogously in 2D), we get, recalling that $\rho^k$ is bounded uniformly in $k$, being $\pphi^k\in[0,1]^N$,
\begin{align*}
&\left\vert\left(\rho^k(\mathbf{v}\cdot\nabla)\mathbf{v},\dfrac{\AA\mathbf{v}-\AA\mathbf{v}^k}{h}\right)_\Omega\right\vert\leq C\Vert \mathbf{v}\Vert_{\mathbf{L}^6(\Omega)}\Vert \nabla\mathbf{v}\Vert_{\mathbf{L}^3(\Omega)}\left\Vert\dfrac{ \AA\mathbf{v}-\AA\mathbf{v}^k}{h}\right\Vert
\\&\leq C\Vert \nabla\mathbf{v}\Vert^\frac 3 2\Vert \AA\mathbf{v}\Vert^\frac 1 2\left\Vert\dfrac{ \AA\mathbf{v}-\AA\mathbf{v}^k}{h}\right\Vert\leq \frac{\alpha}8\left\Vert\dfrac{\AA\mathbf{v}-\AA\mathbf{v}^k}{h}\right\Vert^2+\frac{C}{\alpha}\Vert\nabla\mathbf{v}\Vert^6+C\Vert \AA\mathbf{v}\Vert^2.
\end{align*} 
Then, concerning the other terms, recalling the definition of $\mathbf{J}_\rho$, we get (in 3D and analogously in 2D) thanks to \eqref{pp}
\begin{align*}
&\left\vert\left((\mathbf{J}_\rho\cdot\nabla)\mathbf{v},\dfrac{\AA\mathbf{v}-\AA\mathbf{v}^k}{h} \right)_\Omega\right\vert\leq C\Vert \nabla\ww\Vert_{\mathbf{L}^6(\Omega)}\Vert \nabla\mathbf{v}\Vert_{\mathbf{L}^3(\Omega)}\left\Vert\dfrac{\AA\mathbf{v}-\AA\mathbf{v}^k}{h}\right\Vert\\&\leq C\Vert \nabla\ww\Vert^{\frac 12}\Vert \ww\Vert_{\mathbf{H}^3(\Omega)}^\frac 1 2 \Vert \nabla\mathbf{v}\Vert^\frac 1 2 \Vert \AA\mathbf{v}\Vert^\frac 1 2\left\Vert\dfrac{\AA\mathbf{v}-\AA\mathbf{v}^k}{h}\right\Vert\\&\leq\frac {\alpha} 8\left\Vert\dfrac{\AA\mathbf{v}-\AA\mathbf{v}^k}{h}\right\Vert^2+C\Vert \AA\mathbf{v}\Vert^2+\frac{C}{\alpha}\Vert \nabla\ww\Vert^2\Vert\ww\Vert^2_{\mathbf{H}^3(\Omega)}\Vert \nabla \mathbf{v}\Vert^2.
\end{align*}
Continuing the analysis of \eqref{o3}: we have, recalling that $\nu(\pphi^k)$ and $\nu_{,\pphi}$  are bounded uniformly in $k$, similarly to \eqref{div},
\begin{align*}
&\left\vert\left(2\operatorname{div}(\nu(\pphi^k)D\mathbf{v}),\dfrac{\AA\mathbf{v}-\AA\mathbf{v}^k}{h}\right)_\Omega\right\vert\\&\leq \left\vert\left(2\nu(\pphi^k)\Delta\mathbf{v},\dfrac{\AA\mathbf{v}-\AA\mathbf{v}^k}{h}\right)_\Omega\right\vert+\left\vert\left(2 D\mathbf{v}(\nabla\pphi^k)^T\nu_{,\pphi}(\pphi^k),\dfrac{\AA\mathbf{v}-\AA\mathbf{v}^k}{h}\right)_\Omega\right\vert
\\&\leq  \frac{\alpha}{8}\left\Vert \dfrac{\AA\mathbf{v}-\AA\mathbf{v}^k}{h}\right\Vert^2+\frac{C}{\alpha}\Vert \AA\mathbf{v}\Vert^2+\Vert \nabla\pphi^k\Vert_{\mathbf{L}^6(\Omega)}\Vert D\mathbf{v}\Vert_{\mathbf{L}^3(\Omega)}\left\Vert \dfrac{\AA\mathbf{v}-\AA\mathbf{v}^k}{h}\right\Vert\\&\leq 
\frac{\alpha}{8}\left\Vert \dfrac{\AA\mathbf{v}-\AA\mathbf{v}^k}{h}\right\Vert^2+\frac{C}{\alpha}\Vert \AA\mathbf{v}\Vert^2+\Vert \nabla\pphi^k\Vert_{\mathbf{L}^6(\Omega)}\Vert \nabla\mathbf{v}\Vert^\frac 1 2 \Vert \AA\mathbf{v}\Vert^\frac1 2\left\Vert \dfrac{\AA\mathbf{v}-\AA\mathbf{v}^k}{h}\right\Vert
\\&\leq  \frac{\alpha}{4}\left\Vert \dfrac{\AA\mathbf{v}-\AA\mathbf{v}^k}{h}\right\Vert^2+\frac C\alpha\Vert \AA\mathbf{v}\Vert^2+\frac C\alpha\Vert \nabla\pphi^k\Vert_{\mathbf{L}^6(\Omega)}^4\Vert \nabla\mathbf{v}\Vert^2.
\end{align*}
We are only left with one term: standard inequalities and embeddings entail
\begin{align*}
&\left\vert\left((\nabla\pphi^k)^T\ww,\dfrac{\AA\mathbf{v}-\AA\mathbf{v}^k}{h}\right)_\Omega\right\vert\leq \Vert \nabla\pphi^k\Vert\Vert\ww\Vert_{\mathbf{L}^\infty(\Omega)}\left\Vert\dfrac{ \AA\mathbf{v}-\AA\mathbf{v}^k}{h}\right\Vert\\&\leq C\Vert \nabla\pphi^k\Vert\Vert\ww\Vert_{\mathbf{H}^2(\Omega)}\left\Vert\dfrac{ \AA\mathbf{v}-\AA\mathbf{v}^k}{h}\right\Vert\leq \frac {\alpha}{8}\left\Vert\dfrac{ \AA\mathbf{v}-\AA\mathbf{v}^k}{h}\right\Vert^2+\frac C\alpha\Vert \mathbf{w}\Vert_{\mathbf{H}^2(\Omega)}^2\Vert \nabla\pphi^k\Vert^2.
\end{align*}

To sum up, coming back to \eqref{o3}, in the end we obtain 
\begin{align*}
\alpha\left\Vert\dfrac{\AA\mathbf{v}-\AA\mathbf{v}^k}{h}\right\Vert^2&\leq \frac{C}{\alpha}\left(\Vert \AA\mathbf{v}\Vert^2+\left\Vert\dfrac{D\mathbf{v}-D\mathbf{v}^k}{h}\right\Vert^2+\Vert \nabla\ww\Vert^2\Vert\ww\Vert^2_{\mathbf{H}^3(\Omega)}\Vert \nabla \mathbf{v}\Vert^2\right.\\&\left.\quad +\Vert \nabla\pphi^k\Vert_{\mathbf{L}^6(\Omega)}^4\Vert \nabla\mathbf{v}\Vert^2+\Vert\nabla\mathbf{v}\Vert^6+\Vert \mathbf{w}\Vert_{\mathbf{H}^2(\Omega)}^2\Vert \nabla\pphi^k\Vert^2\right),
\end{align*}
which is exactly \eqref{c3}. This concludes the proof.
\end{enumerate}
\section{Proof of Theorem \ref{weakk}}\label{sec:Weak}
	In order to prove the theorem we want to make use of the time discretization scheme \eqref{phi0} with $\alpha=0$. To this aim, we need to first approximate the initial datum $\pphi_0$ so that it belongs to $\mathbf{H}^2(\Omega)$. %[{\color{red} Here I propose a different approximating problem just for simplicity of solution, since the constraint on the sum of the components is somehow non classical and needs some justifications, in my opinion}] 
 For example, one can define $\pphi_0^M$ as solution of
	\begin{align}
	\nonumber-\frac 1M\Delta\pphi_0^M+\pphi_0^M&=\pphi_0,&&\quad \text{ in }\Omega,\\ \label{a}
	\sum_{i=1}^N\varphi_{0,i}^M&=1,&&\quad \text{ in }\Omega,\\
	\partial_\mathbf{n}\pphi_0^M&=\mathbf{0},&&\quad\text{ on }\partial\Omega.\nonumber
	\end{align}
	We can easily find first a unique $\widetilde{\pphi}_0^M\in \widetilde{\mathbf{V}}_0$ by Lax-Milgram Theorem such that, in a weak sense,
		\begin{align*}
	-\frac 1M\Delta\widetilde{\pphi}_0^M+\widetilde{\pphi}_0^M&=\pphi_0-\frac 1 N \mathbf{e},&&\quad \text{ in }\Omega,\\
	\sum_{i=1}^N\widetilde{\varphi}_{0,i}^M&=0,&&\quad \text{ in }\Omega,\\
	\partial_\mathbf{n}\widetilde{\pphi}_0^M&=\mathbf{0},&&\quad\text{ on }\partial\Omega,
	\end{align*} 
 with $\mathbf{e}=(1,\ldots,1)$.
	Notice that also $\widetilde{\pphi}_0^M\in \mathbf{H}^2(\Omega)$ by elliptic regularity.
	We then obtain $\pphi_0^M=\widetilde{\pphi}_0^M+\frac 1 N \mathbf{e}\in  \mathbf{H}^2(\Omega)$, satisfying \eqref{a} and, by testing the equation by $-\Delta\pphi_0^M$,
\begin{align*}
	\Vert\nabla \pphi_0^M\Vert^2\leq \Vert \nabla\pphi_0\Vert\Vert \nabla\pphi_0^M\Vert.
	\end{align*}
%whereas by testing it by $-\frac 1 N\Delta\pphi_0^M$ we infer
%\begin{align*}
%\frac 1 {N^2}\Vert \Delta\pphi^M_0\Vert^2\leq \Vert \pphi_0\Vert^2.
%\end{align*}
Moreover, by testing the first  equation in \eqref{a} with $\pphi_0^M$, we obtain
$$
\Vert \pphi_0^M\Vert^2\leq \Vert \pphi_0^M\Vert\Vert \pphi_0\Vert,
$$
 so that in the end we deduce
	\begin{align}
	\Vert\nabla \pphi_0^M\Vert\leq \Vert \nabla\pphi_0\Vert,\quad \Vert \pphi_0^M\Vert\leq \Vert\pphi_0\Vert,
	\label{H1contr}
	\end{align}
	and thus also
	$$
	\Vert\pphi_0- \pphi_0^M\Vert\leq \frac 1 M \Vert \nabla\pphi_0^M\Vert\Vert \nabla(\pphi_0^M-\pphi_0)\Vert\leq \frac2 M\Vert \nabla\pphi_0\Vert^2\to 0\quad \text{as }M\to \infty.
	$$
	By these two estimates, since $\mathbf{H}^1(\Omega)$ is a Hilbert space, thus reflexive, such that $\mathbf{H}^1(\Omega)\hookrightarrow\hookrightarrow \mathbf{L}^2(\Omega)$, we easily end up with 
	$$
	\pphi_0^M\rightharpoonup \pphi_0\quad\text{in }\mathbf{H}^1(\Omega) \text{ as }M\to \infty.
	$$
	Moreover, thanks to \eqref{H1contr}, we have 
	$$
	\limsup_{M\to\infty}\Vert\nabla \pphi_0^M\Vert\leq \Vert \nabla\pphi_0\Vert,
	$$
	with $\nabla\pphi_0^M\rightharpoonup \nabla\pphi_0$ in $\mathbf{L}^2(\Omega)$. This implies that $\nabla\pphi_0^M\to \nabla\pphi_0$ in $\mathbf{L}^2(\Omega)$, and thus also 
		$$
	\pphi_0^M\to \pphi_0\quad\text{in }\mathbf{H}^1(\Omega) \text{ as }M\to \infty.
	$$
	In conclusion, we notice that, clearly, $\overline{\pphi}_0^M=\overline{\pphi}_0$ for any $M$. By testing the equation \eqref{a}$_1$ with $-(\pphi^M_0)^-$ (the negative part taken componentwise) we get
 $$
 \frac 1 M\Vert \nabla(\pphi^M_0)^-\Vert^2+\Vert (\pphi^M_0)^-\Vert^2= -\int_\Omega\pphi_0\cdot(\pphi_0^M)^- dx\leq 0,
 $$
 recalling $\pphi_0\in[0,1]^N$ almost everywhere in $\Omega$, and thus $\pphi_0^M\geq 0$. Since it holds \eqref{a}$_2$, this also implies that $\pphi_0^M\leq 1$. Therefore $\pphi_0^M\in[0,1]^N$ in $\overline{\Omega}$, recalling the embedding $\mathbf{H}^2(\Omega)\hookrightarrow \mathbf{C}(\overline{\Omega})$. Therefore we can actually consider the time discretization scheme \eqref{phi0}, setting $\mathbf{v}^0=\mathbf{v}_0$ and $\pphi^0=\pphi_0^M$, together with $\alpha=0$. We now introduce the functions $f^M$ on $[-h,\infty)$ through $f^M(t):=f^{k+1}$ for any $t\in[(k-1)h,k h)$, with $k\in \N_0$ and $f\in\{\mathbf{v},\pphi\}$. In the case of $\ww^M$, we cannot make the definition starting from $t=-h$, since the initial $\ww_0$ is not defined, but nevertheless we can introduce $\ww^M$ on $[0,\infty)$ as
$\ww^M(t):=\ww^{k+1}$ for any $t\in[kh,(k+1)h)$ for $k\in \N_0$. We also set $\rho^M:=\widetilde{\boldsymbol\rho}\cdot \pphi^M$. We define, for any $t\geq 0$, 
\begin{align*}
&(\Delta_h^+f)(t):=f(t+h)-f(t),\quad (\Delta_h^-f)(t):=f(t)-f(t-h),\\&
\partial_{t,h}^+f(t):=\frac 1 h (\Delta_h^+f)(t),\quad \partial_{t,h}^-f(t):=\frac 1 h (\Delta_h^-f)(t),\\&
f_h:=f(t-h).
\end{align*}
From now we can follow the very same arguments of \cite[Sec.5]{ADG} to pass to the limit in $N$. We only give a sketch of the proof, and we refer \cite{ADG} for the complete details. For an arbitrary $\boldsymbol\psi\in \mathbf{C}^\infty_{0,\sigma}(\Omega\times(0,T))$ we set $\widetilde{\boldsymbol\psi}:=\int_{kh}^{(k+1)h}\boldsymbol\psi dt$ and use it as a test function in \eqref{phi0}$_1$, summing over $k\in\N_0$, to get 
\begin{align}
&\nonumber\int_0^\infty\int_\Omega \partial_{t,h}^-(\rho^M\mathbf{v}^M)\cdot \boldsymbol\psi dx dt+\int_0^\infty\int_\Omega \operatorname{div}(\rho_h^M\mathbf{v}^M\otimes \mathbf{v}^M)\cdot \boldsymbol\psi dxdt \\&\nonumber
+\int_0^\infty\int_\Omega 2\nu(\pphi_h^M)D\mathbf{v}^M: D\boldsymbol\psi dx dt -\int_0^\infty\int_\Omega(\mathbf{v}^M\otimes \mathbf{J}_\rho^M):D\boldsymbol\psi dxdt\\&
=\int_0^\infty\int_\Omega(\nabla\pphi^M_h)^T\nabla\ww^Mdxdt
\label{velo}
\end{align}
for all $\boldsymbol\psi\in \mathbf{C}^\infty_{0,\sigma}(\Omega\times(0,T))$. We point out that, given $\boldsymbol\psi\in \mathbf{C}^\infty_{0,\sigma}(\Omega\times(0,T))$ there exists $\overline{h}=\overline{h}(\operatorname{supp}(\boldsymbol\psi))>0$ sufficiently small so that 
\begin{align*}
\int_0^\infty\int_\Omega \partial_{t,h}^-(\rho^M\mathbf{v}^M)\cdot \boldsymbol\psi dx dt=-\int_0^\infty\int_\Omega (\rho^M\mathbf{v}^M)\cdot \partial_{t,h}^+\boldsymbol\psi dx dt,\quad \forall h\leq \overline{h}.
\end{align*}
Analogously we get
\begin{align}
\int_0^\infty\int_\Omega\partial_{t,h}^-\pphi^M\cdot \boldsymbol\zeta dxdt-\int_0^\infty\int_\Omega \mathbf{v}^M\cdot (\pphi_h^M)^T\nabla\boldsymbol\zeta dx dt=\int_0^\infty\int_\Omega\mathbf{M}(\pphi_h^M)\nabla\ww^M:\nabla\boldsymbol\zeta dx dt,
\label{dtp}
\end{align}
 for all $\boldsymbol\zeta\in C_0((0,\infty);\mathbf{C}^1(\overline{\Omega}))$, and 
 \begin{align}
 \ww^M=-\Delta\pphi^M+\mathbf{P}\left(\PsipA(\pphi^M)-\mathbf{A}\frac{\pphi^M+\pphi^M_h}{2}\right)
 \label{ww}
 \end{align}
 holds pointwise in $\Omega\times(0,\infty)$ almost everywhere. Now we observe that $\int_{kh}^{(k+1)h}f^Mdt=h f^{k+1}$, so that \eqref{ener} can be rewritten, for any $k\in \N_0$, as 
	\begin{align}
\nonumber&E_{tot}(\rho^M(kh),\mathbf{v}^M(kh),\pphi^M(kh))+\int_{kh}^{(k+1)h}\int_\Omega \mathbf{M}(\pphi^M_h)\nabla \ww^M:\nabla\ww^M dxd\tau\\&+\frac{1}{h}\int_{kh}^{(k+1)h}\int_\Omega \rho^M\frac{\vert \mathbf{v}^M-\mathbf{v}^M_h\vert^2}2dxd\tau+2\int_{kh}^{(k+1)h}\int_\Omega \nu(\pphi_h^M)\vert D\mathbf{v}^M\vert^2dxd\tau
\nonumber\\&+\frac 1 {2h}\int_{kh}^{(k+1)h}\int_\Omega\vert \nabla\pphi^M-\nabla \pphi^M_h\vert^2dxd\tau\leq E_{tot}(\rho^M((k-1)h),\mathbf{v}^M((k-1)h),\pphi^M((k-1)h)).
\label{enerbis}
\end{align}
We now sum up these inequalities over $k=0,\ldots,Lh$, for any $L>0$, $L\in \N$, ending up with 
 	\begin{align}
 \nonumber&E_{tot}(\rho^M((L-1)h),\mathbf{v}^M((L-1)h),\pphi^M((L-1)h))+\int_{0}^{Lh}\int_\Omega \mathbf{M}(\pphi^M_h)\nabla \ww^M:\nabla\ww^M dxd\tau\\&+\frac{1}{h}\int_{0}^{Lh}\int_\Omega \rho^M\frac{\vert \mathbf{v}^M-\mathbf{v}^M_h\vert^2}2dxd\tau+2\int_{0}^{Lh}\int_\Omega \nu(\pphi_h^M)\vert D\mathbf{v}^M\vert^2dxd\tau
 \nonumber\\&+\frac 1 {2h}\int_{0}^{Lh}\int_\Omega\vert \nabla\pphi^M-\nabla \pphi^M_h\vert^2dxd\leq E_{tot}(\rho_0^M,\mathbf{v}_0,\pphi^M_0)\leq C,
 \label{enerbis2}
 \end{align}
 uniformly in $M$, thanks to \eqref{H1contr}. We have set $\rho_0^M:=\widetilde{\boldsymbol{{\rho}}}\cdot \pphi_0^M>0$. This means, since $\rho^M\geq \min_{i=1,\ldots,N}\widetilde{\rho}_i$, that we have obtained the following uniform bounds (recall the conservation of mass, so that $\overline{\pphi}^M\equiv\overline{\pphi}_0^M=\overline{\pphi}_0$)
 \begin{align}
 &\nonumber\Vert \mathbf{v}^M\Vert_{L^\infty(0,\infty;\mathbf{H}_\sigma)}+ \Vert \mathbf{v}^M\Vert_{L^2(0,\infty;\mathbf{V}_\sigma)}+\Vert \pphi^M\Vert_{L^\infty(0,\infty;\mathbf{H}^1(\Omega))}\\&+ \frac{1}{\sqrt{h}}\Vert\pphi^M-\pphi^M_h\Vert_{L^2(0,\infty;\mathbf{H}^1(\Omega))}+\Vert \nabla\ww^M\Vert_{L^2(0,\infty;\mathbf{L}^2(\Omega))}\leq C.
 \label{unif}
 \end{align}
 Moreover, by Theorem \ref{steaddy} point (3) applied to \eqref{ww}, with $\mathbf f=\ww^M+\mathbf{P}\left(\mathbf{A}\frac{\pphi^M+\pphi^M_h}{2}\right)$ and $\mathbf{m}=\overline{\pphi}_0$, we get, for any $T>0$,
 \begin{align}
 \Vert \pphi^M\Vert_{L^2(0,T;\mathbf{W}^{2,p}(\Omega))}+\Vert \PsipA(\pphi^M)\Vert_{L^2(0,T;\mathbf{L}^p(\Omega))}+\left\Vert\overline{\ww}^M\right\Vert_{L^2(0,T)}\leq C(T),
 \label{wwm}
 \end{align}
 for $p\in[2,\infty)$ when $n=2$ and $p\in[2,6]$ if $n=3$.
 This means, by compactness arguments, that there exists a triple $(\mathbf{v},\pphi,\ww)$ such that, up to nonrelabeled subsequences,
 \begin{alignat*}{2}
& \mathbf{v}^M\rightharpoonup^* \mathbf{v}&&\quad\text{ in }L^\infty(0,\infty;\mathbf{H}_\sigma),\\&
\mathbf{v}^M\rightharpoonup \mathbf{v}&&\quad\text{ in }L^2(0,\infty;\mathbf{V}_\sigma),\\&
\pphi^M\rightharpoonup^* \pphi&&\quad\text{ in }L^\infty(0,\infty;\mathbf{H}^1(\Omega)),\\&
\pphi^M\rightharpoonup \pphi&&\quad\text{ in }L^2(0,T;\mathbf{W}^{2,p}(\Omega)),\quad \forall T>0,\\&
\ww^M\rightharpoonup\ww&&\quad\text{ in }L^2(0,T;\mathbf{H}^1(\Omega)),\quad\quad \forall T>0,\\&
\nabla\ww^M\rightharpoonup \nabla\ww &&\quad \text{ in }L^2(0,\infty;\mathbf{L}^2(\Omega)),
 \end{alignat*}
 with $p\in[2,\infty)$ if $n=2$, $p\in[2,6]$ if $n=3$.
We now define $\widetilde{\pphi}^M(t):=\frac 1 h\int_{t-h}^t\pphi^M(s)ds$, which is the piecewise linear interpolant of $\pphi^M(kh)$, $k\in \N_0$. Then we have $\partial_{t}\widetilde{\pphi}^M=\partial_{t,h}^-\pphi^M$.
Thus, it is easy to see from \eqref{dtp} and \eqref{unif} that 
\begin{align}
\Vert \partial_t\widetilde{\pphi}^M\Vert_{L^2(0,\infty;(\mathbf{H}^1(\Omega))')}\leq C,
\label{dt1}
\end{align}
uniformly in $M$. Since we have
\begin{align}
\Vert \widetilde{\pphi}^M-\pphi^M\Vert_{(\mathbf{H}^1(\Omega))'}\leq 2h\Vert \partial_t\widetilde{\pphi}^M\Vert_{(\mathbf{H}^1(\Omega))'},
\label{ident}
\end{align}
we immediately infer
\begin{align}
\widetilde{\pphi}^M- \pphi^M\to 0\quad \text{in }L^2(0,\infty;(\mathbf{H}^1(\Omega))'),\quad \text{as }M\to \infty.
\label{cov}
\end{align}
Due to the fact that, clearly, from \eqref{unif} we also have $\widetilde{\pphi}^M$ uniformly bounded in $L^\infty(0;\infty;\mathbf{H}^1(\Omega))$, by Aubin-Lions Theorem we deduce, thanks to \eqref{dt1},
$$
\widetilde{\pphi}^M\to \widetilde{\pphi},\quad \text{ in }L^2(0,T;\mathbf{L}^2(\Omega)),
$$
 for any $0<T<\infty$, implying also, up to another subsequence, that the convergence takes place pointwise almost everywhere in $\Omega\times(0,\infty)$. Therefore, thanks to \eqref{cov} we can deduce that $\widetilde{\pphi}\equiv \pphi$, entailing then
 \begin{align*}
& \partial_t \widetilde{\pphi}^M\equiv\partial_{t,h}^-\pphi^M\rightharpoonup \partial_t\pphi\quad \text{ in }L^2(0,\infty; (\mathbf{H}^1(\Omega))'),\\&
{\pphi}^M\to \pphi\quad \text{ in }L^2(0,T; \mathbf{H}^1(\Omega)),\quad \forall 0<T<\infty,\\&
{\pphi}^M\to \pphi\quad \text{a.e. in }\Omega\times(0,\infty).
 \end{align*}
 Where the last two convergences are obtained by interpolation from the bound \eqref{wwm} and the convergence \eqref{cov}.
 Notice that, by the last convergence we also deduce that $\sum_{i=1}^N\varphi_i\equiv 1$ almost everywhere in $\Omega\times(0,\infty)$, and also that $\pphi\in[0,1]^N$ almost everywhere in $\Omega\times(0,\infty)$. Furthermore, thanks to \eqref{unif}, we also have that $\pphi^M-\pphi^M_h\to 0$ in $L^2(0,\infty;\mathbf{H}^1(\Omega))$, and thus also 
 $$
 \pphi_h^M\to \pphi
 $$
 in $L^2(0,T;\mathbf{H}^1(\Omega))$ for any $T>0$, and, up to subsequences, also pointwise almost everywhere in $\Omega\times(0,\infty)$.
Moreover, since $\pphi\in H^1_{uloc}([0,\infty);(\mathbf{H}^1(\Omega))')\cap L^2_{uloc}([0,\infty);\mathbf{H}^1(\Omega))\hookrightarrow BUC([0,\infty);\mathbf{L}^2(\Omega))$ and since $\pphi\in L^\infty(0,\infty;\mathbf{H}^1(\Omega))$, we have by \cite[Lemma 2.1]{ADG} that
$$
\pphi\in BC_w([0,\infty);\mathbf{H}^1(\Omega))).
$$
Finally, again by Aubin-Lions Lemma, since $\widetilde{\pphi}^M$ is also uniformly bounded in the spaces $H^1(0,\infty;(\mathbf{H}^1(\Omega))')\cap L^\infty(0,\infty;\mathbf{H}^1(\Omega))$, it holds, up to subsequences, $$\widetilde{\pphi}^M\to \pphi\quad \text{ in }BC([0,T];\mathbf{L}^2(\Omega)),\quad \forall 0<T<\infty,$$
so that
$$
\widetilde{\pphi}^M(0)\to \pphi(0)\quad \text{in }\mathbf{L}^2(\Omega).
$$
Since also $$\widetilde{\pphi}^M(0)=\frac 1 h\int_{-h}^0\pphi^M(s)ds=\pphi_0^M\to \pphi_0\quad \text{ in }\mathbf{H}^1(\Omega),$$
we deduce ${\pphi}(0)=\pphi_0$ almost everywhere in $\Omega$.
The same arguments hold for $\rho^M$, since it is a linear combination of the components of $\pphi^M$. To show the convergence in \eqref{ww}, we need the following observations. First, we have, for any $k\in\N$, that $\ww^k\in \mathbf{H}^2(\Omega)\hookrightarrow \mathbf{L}^\infty(\Omega)$, so that, by Theorem \ref{steaddy} Point (2) there exists $\delta_k>0$, depending on $k$, such that
$$
0<\delta_k<\pphi^k<1-\delta_k\quad \text{ in }\overline{\Omega}.
$$ 
Therefore, we are rigorously allowed to observe, since $\psi^\prime$ is monotone increasing, that, for any $k\in\N$ and for any $i=1,\ldots,N$,
\begin{align*}
-\psi^\prime(\eta)h\left\vert \{x\in\Omega:\ \varphi_i^k\leq \eta\}\right\vert_n\leq -h\int_{ \{x\in\Omega:\ \varphi_i^k\leq \eta\}}\psi^\prime(\varphi^k_i)dx,\quad \forall \eta\in\left(0,\frac 1 2\right)
\end{align*}
i.e., for any $k\in\N_0$,
\begin{align*}
&-\psi^\prime(\eta)\int_{kh}^{(k+1)h}\left\vert \{x\in\Omega:\ \varphi_i^M(t)\leq \eta\}\right\vert_n dt\\&\RevA{\leq -\int_{kh}^{(k+1)h}\left(\int_{ \{x\in\Omega:\ \varphi_i^M(t)\leq \eta\}}\psi^\prime(\varphi^M_i(t))dx\right)dt,\quad \forall \eta\in\left(0,\frac 1 2\right).}
\end{align*}
Summing up for $k=0,\ldots,L$ such that $(L+1)h<T$, by \eqref{wwm} we conclude
\begin{align*}
-\psi^\prime(\eta)\left\vert\{(x,t)\in \Omega\times(0,(L+1)h): \ \varphi_i^M(x,t)\leq \eta\}\right\vert_{n+1}\leq \Vert \psi^\prime(\varphi^M_i)\Vert_{L^2(0,T;\mathbf{L}^2(\Omega))}\leq C(T)\quad \forall   \eta\in\left(0,\frac 1 2\right).
\end{align*}
Since for any $L>0$ there exists $T>0$ such that $(L+1)h<T$, the result above holds for any $L>0$ and thus for any $T>0$. 
Since $\pphi^M\to \pphi$ almost everywhere in $\Omega\times(0,\infty)$, we can apply Fatou's Lemma to deduce, for any $T>0$, 
\begin{align*}
&\left\vert\{(x,t)\in \Omega\times(0,T): \ \varphi_i(x,t)\leq \eta\}\right\vert_{n+1}\\&\leq \liminf_{N\to\infty}\left\vert\{(x,t)\in \Omega\times(0,T): \ \varphi_i^M(x,t)\leq \eta\}\right\vert_{n+1}\leq \frac{C(T)}{-\psi^\prime(\eta)}\quad \forall   \eta\in\left(0,\frac 1 2\right),
\end{align*}
and then, since, by assumption (\textbf{E1}), $-\psi^\prime(s)\to +\infty$ as $s\to 0$, letting $\eta\to0$ we end up with 
\begin{align}
0<	\pphi<1,\quad \text{a.e. in }\Omega\times(0,\infty),
\label{sepa}
\end{align}
by repeating the same argument for any component $i=1,\ldots,N$. Thanks to this result we deduce that 
$$
\PsipA(\pphi^M)\to \PsipA(\pphi)\quad \text{a.e. in }\Omega\times(0,\infty), 
$$
so that, since also $\PsipA(\pphi^M)$ is uniformly bounded in $L^2(0,T;\mathbf{L}^2(\Omega))$ for any $0<T<\infty$, we infer by the generalized Lebesgue Theorem that
$$
\PsipA(\pphi^M)\rightharpoonup \PsipA(\pphi) \quad \text{ in }L^2(0,T;\mathbf{L}^2(\Omega)), \quad \forall 0<T<\infty,
$$
which is enough to pass to the limit as $M\to \infty$ in the weak formulation of \eqref{ww}. 

To conclude the passage to the limit, we need to gain some strong convergence on the velocity sequence $\mathbf{v}^M$. Since this argument is identical to the one of \cite[Sec.5.1]{ADG}, we only refer the main steps. First, we define $\widetilde{\rho^M\mathbf{v}^M}(t):=\int_{t-h}^h\rho^M\mathbf{v}^Mds$ and then show, by standard arguments exploiting \eqref{unif} and \eqref{wwm}, that 
$$
\Vert \partial_t\mathbb{P}\widetilde{\rho^M\mathbf{v}^M}\Vert_{L^{\frac 8 7}(0,T;\mathbf{W}^{-1,4}(\Omega))}\leq C(T),\quad \forall 0<T<\infty. 
$$
 Then, since also $\mathbb{P}\widetilde{\rho^M\mathbf{v}^M}$ is uniformly bounded in $L^2(0,T;\mathbf{H}^1(\Omega))$, by Aubin-Lions Lemma we deduce 
$$
\mathbb{P}\widetilde{\rho^M\mathbf{v}^M}\to \mathbf{w}\quad \text{ in }L^2(0,T;\mathbf{L}^2(\Omega))\quad \forall 0<T<\infty.
$$
Since also $\mathbb{P}\widetilde{\rho^M\mathbf{v}^M}$ is uniformly bounded in $L^\infty(0,\infty;\mathbf{L}^2(\Omega))$, we have $\mathbf{w}\in L^\infty(0,\infty;\mathbf{L}^2(\Omega))$. Now, the projector $\mathbb{P}: L^2(0,T;\mathbf{L}^2(\Omega))\to L^2(0,T;\mathbf{H}_\sigma)$ is weakly continuous, so that, by the weak convergence (coming from the previous convergence results) $\widetilde{\rho^M\mathbf{v}^M}\rightharpoonup \rho\mathbf{v}$ in $L^2(0,T;\mathbf{L}^2(\Omega))$, we infer $\mathbf{w}=\mathbb{P}\rho\mathbf{v}$. Again by the fact that
$$
\Vert \mathbb{P}\widetilde{\rho^M\mathbf{v}^M}-\mathbb{P}{\rho^M\mathbf{v}^M}\Vert_{L^{\frac 8 7}(0,T;\mathbf{W}^{-1,4}(\Omega))}\leq 2h\Vert \partial_t\mathbb{P}\widetilde{\rho^M\mathbf{v}^M}\Vert_{L^{\frac 8 7}(0,T;\mathbf{W}^{-1,4}(\Omega))},
$$
we also infer that also 
$$
\mathbb{P}{\rho^M\mathbf{v}^M}\to \mathbb{P}\rho\mathbf{v}\quad \text{ in }L^2(0,T;\mathbf{L}^2(\Omega)),\quad \forall 0<T<\infty.
$$
 Now a standard argument (see \cite[Sec.5.1]{ADG}) allows to deduce
$$
(\rho^{M})^{\frac 1 2}\mathbf{v}^M\to \rho^\frac 1 2\mathbf{v},
$$ 
from which, recalling that, thanks to \eqref{sepa}, $\rho\geq \min_{i=1,\ldots,N}\widetilde{\rho}_i>0$, and 
$$
\rho^M\to\rho\quad \text{ a.e. in }\Omega\times(0,\infty),
$$ 
we easily deduce
$$
\mathbf{v}^M=(\rho^M)^{-\frac 1 2 }(\rho^{M})^{\frac 1 2}\mathbf{v}^M\to \mathbf{v}\quad \text{ in }L^2(0,T;\mathbf{L}^2(\Omega))\quad \forall 0<T<\infty,
$$
which also implies, up to another subsequence, the almost everywhere convergence on $\Omega\times(0,\infty)$. We can now pass to the limit as $M\to \infty$ in \eqref{velo} and \eqref{dtp}. To conclude the proof, one should show that $\mathbf{v}\in BC_w([0,T];\mathbf{L}^2(\Omega))$ for any $0<T<\infty$ and that $\mathbf{v}(0)=\mathbf{v}_0$, together with the validity of the energy estimate \eqref{energ}. Since the arguments needed are the same, we refer to \cite[Sec.5.2-5.3]{ADG} for the detailed proof. The proof is concluded.

\section{Proof of Theorem \ref{strong1}}
\label{str}
\subsection{Existence}
	The proof of the theorem relies again on the time discretization scheme and on Theorem \ref{disc}. In particular, we need to make a specific choice for the initial value $\pphi^0$ of the scheme, in order to be able to define $\ww^0\in \mathbf{L}^\infty(\Omega)\cap \widetilde{\mathbf{V}}_0$. Let us consider, for some $T>0$, the Cahn-Hilliard equation
	\begin{alignat}{2}
	\nonumber\partial_t\widetilde{\pphi}-\Delta \widetilde{\ww}&=0,&&\quad\text{in }\Omega\times(0,T),\\ \nonumber
	\widetilde{\ww}&=-\Delta{\widetilde{\pphi}}+\mathbf{P}{\Psi}_{,\pphi}(\widetilde{\pphi}),&&\quad\text{in }\Omega\times(0,T),\\ \nonumber
	\partial_\mathbf{n}\widetilde{\pphi} &=0,&&\quad \text{on }\partial\Omega\times(0,T),\\
	\partial_\mathbf{n}\widetilde{\ww} &=0,&&\quad \text{on }\partial\Omega\times(0,T),	\label{a1}\\
	\widetilde{\pphi}(0)&=\pphi_0,&&\quad \text{in }\Omega, \nonumber\\ \nonumber
	\sum_{i=1}^N\widetilde{\varphi}_i&\equiv 1,&&\quad\text{in }\Omega\times(0,T).
	\end{alignat} 
	By a similar approximating scheme as in the proof of \cite[Thm.3.1]{GGPS} (see also \cite[Remark 3.10]{GGPS}) we can easily prove (by testing \eqref{a1}$_1$ first by $\widetilde{\ww}$ and then by $\partial_t\widetilde{\ww}$) that there exists a unique solution $(\widetilde{\pphi},\widetilde{\ww})$ such that
	\begin{align*}
	&\widetilde{\pphi}\in C([0,T];\mathbf{H}^1(\Omega))\cap H^1(0,T;\mathbf{H}^1(\Omega)),\\&
	\widetilde{\ww}\in L^\infty(0,T;\mathbf{H}^1(\Omega))\cap L^2(0,T;\mathbf{H}^3(\Omega)),
	\end{align*}
	so that 
	\begin{align}
	\Vert \nabla\widetilde{\ww}\Vert_{L^\infty(0,T;\mathbf{L}^2(\Omega))}\leq C\Vert\nabla\ww_0\Vert.
	\label{p2}
	\end{align}
	This means that there exists a sequence $t_M\to 0$ as $M\to \infty$ satisfying 
	\begin{align}
	\sum_{i=1}^N\widetilde{\varphi}_i(t_M)\equiv 1,\quad \{\widetilde{\ww}(t_M)\}_M\subset \mathbf{H}^2(\Omega),\quad \Vert \nabla\widetilde{\ww}(t_M)\Vert\leq C\Vert\nabla\ww_0\Vert,\quad \forall M\in\N,
	\label{conv1}
	\end{align}
	and, by the regularity of $\widetilde{\pphi}$,
	\begin{align}
	\widetilde{\pphi}(t_M)\to \pphi_0\quad \text{ in }\mathbf{H}^1(\Omega).
	\label{conv2}
	\end{align}
	Moreover, it also holds $\overline{\widetilde{\pphi}}\equiv\overline{\pphi}_0$. We then define $\ww_0^M:=\widetilde{\ww}(t_M)$ and $\pphi_0^M:=\widetilde{\pphi}(t_M)$. We additionally consider a $\mathbf{W}_\sigma$ approximating sequence $\{\mathbf{v}_0^\alpha\}_\alpha$, for $\alpha>0$ of $\mathbf{v}_0$ in the following way $\mathbf{v}_0^\alpha:=(\mathbb{I}+\sqrt{\alpha}\AA)^{-1}\mathbf{v}_0\in \mathbf{W}_\sigma$, where $\mathbb{I}$ is the identity operator. Notice that $\mathbf{v}_0^\alpha$ satisfies
	$$
	\mathbf{v}_0^\alpha+\sqrt{\alpha}\AA \mathbf{v}_0^\alpha=\mathbf{v}_0.
	$$
	By testing against $\AA\mathbf{v}_0^\alpha$ we easily obtain 
	\begin{align}
	\Vert\nabla\mathbf{v}_0^\alpha\Vert\leq \Vert \nabla\mathbf{v}_0\Vert,\quad \alpha\Vert\AA\mathbf{v}_0^\alpha\Vert^2\leq C\Vert \mathbf{v}_0\Vert^2. 
	\label{appr}
	\end{align}
	Therefore, we have 
	$$
	\Vert \mathbf{v}_0^\alpha-\mathbf{v}_0\Vert^2\leq \sqrt{\alpha}\Vert \nabla\mathbf{v}_0^\alpha\Vert\Vert \nabla\mathbf{v}_0^\alpha-\nabla\mathbf{v}_0\Vert\leq C\sqrt{\alpha}\to 0\quad\text{as }\alpha\to 0.	$$
	This implies, together with \eqref{appr},
	\begin{align}
	 \mathbf{v}_0^\alpha\to \mathbf{v}_0\quad \text{in }\mathbf{V}_\sigma\quad \text{as }\alpha\to 0.
	\label{con}
	\end{align}
	We can thus consider the discretization scheme, with $\pphi^0=\pphi_0^M$, $\ww^0=\ww_0^M$, $\mathbf{v}^0=\mathbf{v}_0^\alpha$.  Clearly by construction it holds
	$$
	-\Delta\pphi^0+\mathbf{P}\Psip(\pphi^0)=\ww^0\quad \text{in } \Omega,
	$$
	for any $M$. By this observation, Theorem \ref{disc} holds for any $k\in\N_0$: all the assumptions are verified with the choice $\pphi^0=\pphi_0^M$, $\ww^0=\ww_0^M$ and $\mathbf{v}^0=\mathbf{v}_0^\alpha$ with $\alpha\in(0,1]$. This means that \eqref{ener} holds. By repeating the same arguments as in the proof of Theorem \ref{weakk}, having introduced the functions $(\mathbf{v}^M,\pphi^M,\ww^M)$, this time all defined on $[-h,+\infty)$, we obtain, as in \eqref{enerbis2},
	 	\begin{align}
	\nonumber&E_{tot}(\rho^M((L-1)h),\mathbf{v}^M((L-1)h),\pphi^M(L-1))+\frac{\alpha }{2}\int_\Omega \left\vert\nabla\mathbf{v}^M((L-1)h)\right\vert^2dx\\&+\frac{\alpha }{2h}\int_{0}^{Lh}\int_\Omega \left\vert\nabla\mathbf{v}^M-\nabla \mathbf{v}^M_h\right\vert^2dx+\int_{0}^{Lh}\int_\Omega \mathbf{M}(\pphi^M_h)\nabla \ww^M:\nabla\ww^M dxd\tau\nonumber\\&+\frac{1}{h}\int_{0}^{Lh}\int_\Omega \rho^M\frac{\vert \mathbf{v}^M-\mathbf{v}^M_h\vert^2}2dxd\tau+2\int_{0}^{Lh}\int_\Omega \nu(\pphi_h^M)\vert D\mathbf{v}^M\vert^2dxd\tau
	+\frac 1 {2h}\int_{0}^{Lh}\int_\Omega\vert \nabla\pphi^M-\nabla \pphi^M_h\vert^2dxd\tau\nonumber\\&\leq E_{tot}(\rho_0^M,\mathbf{v}_0^\alpha,\pphi^M_0)+\frac{\alpha }{2}\int_\Omega \left\vert\nabla\mathbf{v}_0^\alpha\right\vert^2dx\leq C,
	\label{enerbis2bis}
	\end{align}
where the last bound, uniform in $M,L$ and $\alpha$ (for $M$ sufficiently large), is due to \eqref{conv2} and \eqref{appr}. 	
This ensures (see \eqref{unif}) (recall the conservation of mass, so that $\overline{\pphi}^M\equiv\overline{\pphi}_0^M=\overline{\pphi}_0$)
\begin{align}
&\sqrt{\alpha}\Vert\mathbf{v}^M\Vert_{{L}^{\infty}(0,\infty;\mathbf{V}_\sigma)}+\Vert \mathbf{v}^M\Vert_{L^\infty(0,\infty;\mathbf{H}_\sigma)}+ \Vert \mathbf{v}^M\Vert_{L^2(0,\infty;\mathbf{V}_\sigma)}\nonumber\\&+\Vert \pphi^M\Vert_{L^\infty(0,\infty;\mathbf{H}^1(\Omega))}+\Vert \nabla\ww^M\Vert_{L^2(0,\infty;\mathbf{L}^2(\Omega))}\leq C,
\label{unif1}
\end{align}
together with
 \begin{align}
\Vert \pphi^M\Vert_{L^2(0,T;\mathbf{H}^2(\Omega))}+\Vert \PsipA(\pphi^M)\Vert_{L^2(0,T;\mathbf{L}^2(\Omega))}+\left\Vert\overline{\ww}^M\right\Vert_{L^2(0,T)}\leq C(T),
\label{wwm1}
\end{align}	
uniformly in $M$, by an application of Theorem \ref{steaddy} (see also \eqref{wwm}). Having defined the linear interpolation $\widetilde{\pphi}^M$, we also deduce (see \eqref{dt1})
\begin{align}
\Vert \partial_t\widetilde{\pphi}^M\Vert_{L^2(0,\infty;(\mathbf{H}^1(\Omega))')}\leq C.
\label{dt3}
\end{align}
Now we can consider \eqref{h}, which is valid for any $k\in\N_0$, having assumed the convention $\pphi^{-1}=\pphi^0$. This can be rewritten as, for any $k\in\N$,
\begin{align*}
&\nonumber\frac 1 {2} \left(\mathbf{M}\nabla \ww^M(kh),\nabla\ww^M(kh)\right)+\frac 1 2\int_{kh}^{(k+1)h} \left\Vert\nabla\partial_t\widetilde{\pphi}^M\right\Vert^2ds\\&\leq \frac 1 {2} \left(\mathbf{M}\nabla \ww^M((k-1)h),\nabla\ww^M((k-1)h)\right)\\&\RevA{+C(\min_{i=1,\ldots,N}\overline{\varphi}^0_i)\int_{kh}^{(k+1)h}\Vert \nabla\mathbf{v}^M\Vert\Vert\mathbf{v}^M\Vert(1+\Vert \nabla\ww^M_h\Vert^2+\Vert \nabla\ww^M\Vert^2)ds\nonumber}\\&+C(\min_{i=1,\ldots,N}\overline{\varphi}^0_i)\int_{kh}^{(k+1)h}\Vert \nabla\mathbf{v}^M\Vert^2(1+\Vert\nabla\ww^M_h\Vert^2)ds+C\int_{kh}^{(k+1)h}\left(\left\Vert \partial_t\widetilde{\pphi}^M\right\Vert^{2}_{\mathbf{V}_0'}+\left\Vert \partial_t\widetilde{\pphi}^M_h\right\Vert^2_{\mathbf{V}_0'}\right)ds,
\end{align*}	
where we used $\partial_t\widetilde{\pphi}^M(t)=0$ when $t\in[-h,0)$. This implies, by summing this inequality over $k=0,\ldots,L$, for any fixed $L>0$, 
\begin{align}
&\nonumber\frac 1 {2} \left(\mathbf{M}\nabla \ww^M((L-1)h),\nabla\ww^M((L-1)h)\right)+\frac 1 2\int_{0}^{Lh} \left\Vert\nabla\partial_t\widetilde{\pphi}^M\right\Vert^2ds\\&\leq \frac 1 {2} \left(\mathbf{M}\nabla \ww^M_0,\nabla\ww^M_0\right)+C(\min_{i=1,\ldots,N}\overline{\varphi}^0_i)\int_{0}^{Lh}\Vert \nabla\mathbf{v}^M\Vert\Vert\mathbf{v}^M\Vert(1+\Vert \nabla\ww^M\Vert^2)ds\nonumber\\&\nonumber +C(\min_{i=1,\ldots,N}\overline{\varphi}^0_i)\int_{0}^{Lh}(1+\Vert\nabla\mathbf{v}^M\Vert^2)\Vert \nabla\ww^M_h\Vert^2ds+C\int_{0}^{Lh}\left(\left\Vert \partial_t\widetilde{\pphi}^M\right\Vert^{2}_{\mathbf{V}_0'}+\left\Vert \partial_t\widetilde{\pphi}^M_h\right\Vert^2_{\mathbf{V}_0'}\right)ds\\&\leq \frac 1 {2} \left(\mathbf{M}\nabla \ww^M_0,\nabla\ww^M_0\right)+C(\min_{i=1,\ldots,N}\overline{\varphi}^0_i)\int_{0}^{Lh}\Vert \nabla\mathbf{v}^M\Vert\Vert \nabla\ww^M\Vert^2ds\nonumber\\&+C(\min_{i=1,\ldots,N}\overline{\varphi}^0_i)\int_{-h}^{(L-1)h}\Vert \nabla\mathbf{v}^M(h+s)\Vert^2\Vert \nabla\ww^M(s)\Vert^2ds+C_0,
\label{d1}
\end{align}
where in the last step we performed a change of variables and exploited also \eqref{unif1} and \eqref{dt3}, being $C_0>0$ uniform in $M,L,\alpha$. Observe now that, by the definition of $\ww^M$ and Cauchy-Schwarz inequality, 
\begin{align*}
&\int_{(L-1)h}^{Lh}\Vert \nabla\mathbf{v}^M\Vert\Vert \nabla\ww^M\Vert^2ds=\Vert \nabla\ww^M((L-1)h)\Vert^2\int_{(L-1)h}^{Lh}\Vert \nabla\mathbf{v}^M\Vert ds\\&\leq \sqrt{h}\Vert \nabla\ww^M((L-1)h)\Vert^2\left(\int_{(L-1)h}^{Lh}\Vert \nabla\mathbf{v}^M\Vert^2 ds\right)^\frac 1 2\leq C_1\sqrt{h}\left(\mathbf{M}\nabla \ww^M((L-1)h),\nabla\ww^M(
(L-1)h)\right),
\end{align*}
for $C_1>0$ uniform in $M,L,\alpha$, where we exploited once more \eqref{unif1} and property \eqref{nondeg} of $\mathbf{M}$. Recalling now that $\ww^M(t)=\ww^M((L-1)h)$ for any $t\in[(L-1)h,Lh)$, we deduce 
\begin{align}
&\nonumber \left(\frac 1 {2}-C_1\sqrt{h}\right)\left(\mathbf{M}\nabla \ww^M(t),\nabla\ww^M(t)\right)+\frac 1 2\int_{0}^{t} \left\Vert\nabla\partial_t\widetilde{\pphi}^M\right\Vert^2ds\\&\leq \frac 1 {2} \left(\mathbf{M}\nabla \ww^M_0,\nabla\ww^M_0\right)+C_0+C(\min_{i=1,\ldots,N}\overline{\varphi}^0_i)\int_{-h}^{t}\Vert \nabla\mathbf{v}^M(h+s)\Vert^2\Vert \nabla\ww^M(s)\Vert^2ds,
\label{uloc}
\end{align}
 for any $t\in[(L-1)h,Lh)$ and for any $L>0$, implying that this inequality holds for any $t\geq0$, since neither $C_1$ depends on $L$. We can now choose, e.g., $h<\overline{h}$, with $\overline{h}:=	\frac 1 {16C_1^2}$, to apply Gronwall's Lemma and obtain 
 \begin{align*}
 &\sup_{t\geq 0}\left(\mathbf{M}\nabla \ww^M(t),\nabla\ww^M(t)\right)\leq C(\Vert \ww_0^M\Vert_{\mathbf{H}^1(\Omega)})e^{C(\min_{i=1,\ldots,N}\overline{\varphi}^0_i)\int_{-h}^\infty\Vert\nabla\mathbf{v}^M\Vert^2 ds}\\&\leq C\left(\Vert \ww_0^M\Vert_{\mathbf{H}^1(\Omega)},E_{tot}(\rho_0^M,\mathbf{v}_0^\alpha,\pphi_0^M),\min_{i=1,\ldots,N}\overline{\varphi}^0_i\right)\leq C,
 \end{align*} 
 independently of $N,\alpha$, by exploiting \eqref{unif1}. Another application of Theorem \ref{steaddy} point (3), with the choice $\mathbf f=\ww^M+\mathbf{P}\left(\mathbf{A}\frac{\pphi^M+\pphi^M_h}{2}\right)$ and $\mathbf{m}=\overline{\pphi}_0$, then gives
  \begin{align}
 \Vert \pphi^M\Vert_{L^\infty(0,\infty;\mathbf{H}^2(\Omega))}+\Vert \PsipA(\pphi^M)\Vert_{L^\infty(0,\infty;\mathbf{L}^2(\Omega))}+\left\Vert\overline{\ww}^M\right\Vert_{L^\infty(0,\infty)}\leq C,
 \label{wwm2}
 \end{align}	
 uniformly in $N,\alpha$, so that we can obtain
  \begin{align}
 \Vert\ww^M\Vert_{L^\infty(0,\infty;\mathbf{H}^1(\Omega))}\leq C.
 \label{ext}
 \end{align}
 By a similar argument as in \eqref{d1}, but on any interval $[t,t+1]$, $t\geq0$, we also get
 \begin{align}
 \Vert\ww^M\Vert_{L^\infty(0,\infty;\mathbf{H}^1(\Omega))}+\Vert \partial_t\widetilde{\pphi}^M\Vert_{L^2_{uloc}([0,\infty);\mathbf{H}^1(\Omega))}\leq C,
 \label{ext1}
 \end{align}
 uniformly in $N,\alpha$. By standard comparison in \eqref{phi0}$_2$, it now easy to see that 
 \begin{align}
 \Vert \ww^M\Vert_{L^2_{uloc}([0,\infty);\mathbf{H}^3(\Omega))}\leq C,
 \label{regmu}
 \end{align}
 uniformly in $N,\alpha$. In conclusion, we study the relation \eqref{c1}, which is valid for any $k\in \N_0$ thanks to Remark \ref{H3} and to the choice of $\mathbf{v}_0^\alpha\in \mathbf{W}_\sigma$. By the same arguments as before, we have, for any $T>0$,
 	\begin{align}
 &\alpha\int_0^{Lh}\left\Vert D \partial_t\widetilde{\mathbf{v}}^M\right\Vert^2\nonumber\\&\leq \frac{C}{\alpha}\int_0^{Lh}\left(\Vert \mathbf{w}^{M}\Vert_{\mathbf{H}^1(\Omega)}^2\Vert \nabla\pphi^M_h\Vert^2+\Vert\ww^M\Vert^2_{\mathbf{H}^2(\Omega)}\Vert \nabla \mathbf{v}^M\Vert^2+\Vert\nabla\mathbf{v}^M\Vert^4+\Vert\nabla\mathbf{v}^M\Vert^2\right)ds\leq C(L,\alpha),
 \label{c1b}
 \end{align}
 where $\widetilde{\mathbf{v}}^M(t):=\frac 1 h \int_{t-h}^t\mathbf{v}^Mds$, and thanks to \eqref{unif1}, \eqref{wwm1} and \eqref{regmu}, entailing that
 \begin{align}
 \left\Vert \partial_t\widetilde{\mathbf{v}}^M\right\Vert_{L^2(0,T;\mathbf{H}^1(\Omega))}\leq C(\alpha,T),\quad \forall T>0,
 \label{la}
 \end{align}
 with a constant $C(\alpha,T)$ which depends on $\alpha>0$ and $T$. Now we can exploit \eqref{c2} to get, summing over $k\in \N_0$ up to $k=L-1$,
 	\begin{align}
 &\nonumber\frac\alpha2\left\Vert\AA\mathbf{v}^M((L-1)h)\right\Vert^2\leq \frac\alpha2\Vert\AA\mathbf{v}_0^\alpha\Vert^2+C\int_0^{Lh}\left\Vert\AA\mathbf{v}^M\right\Vert^2dt\\&+C\int_0^{Lh}\left(\Vert \nabla\mathbf{v}^M\Vert^6+\Vert \nabla\ww^M\Vert^2\Vert\ww^M\Vert^2_{\mathbf{H}^2(\Omega)}\Vert \nabla \mathbf{v}^M\Vert^2+\Vert \nabla\pphi^M_h\Vert_{\mathbf{L}^6(\Omega)}^4\Vert \nabla\mathbf{v}^M\Vert^2\right.\nonumber\\&\left.+\Vert \mathbf{w}^M\Vert_{\mathbf{H}^3(\Omega)}^2\Vert \nabla\pphi^M_h\Vert^2+\Vert D\partial_t\widetilde{\mathbf{v}}^M\Vert^2\right)dt,
 \label{c2b}
 \end{align}
 with $C$ independent of $M,L$.
 Observing that 
 $$
 C\int_{(L-1)h}^{Lh}\left\Vert\AA\mathbf{v}^M\right\Vert^2dt=Ch\Vert \AA\mathbf{v}^M((L-1)h)\Vert^2,
 $$
 we can assume $h$ sufficently small so that $\frac\alpha2-Ch\geq \frac \alpha 4>0$, entailing by \eqref{c2b} that, for any $T>0$ and any $t\in[0,T]$,
 	\begin{align*}
 &\nonumber\frac\alpha4\left\Vert\AA\mathbf{v}^M(t)\right\Vert^2\leq \frac\alpha2\Vert\AA\mathbf{v}_0^\alpha\Vert^2+C(\alpha,T),
 \end{align*}
 for any $h$ sufficiently small, where we exploited \eqref{unif1}, \eqref{ext}, \eqref{regmu} and \eqref{la} to bound the last term in the right-hand side of \eqref{c2b}. This entails 
 \begin{align}
 \Vert \mathbf{v}^M\Vert_{L^\infty(0,T;\mathbf{W}_\sigma)}\leq C(T,\alpha),\quad \forall T>0.
 \label{1b}
 \end{align}
 In conclusion, by similar arguments applied to \eqref{c3}, we end up with 
 	\begin{align*}
 &\alpha\int_0^{Lh}\left\Vert{\AA\partial_t\widetilde{\mathbf{v}}^M}\right\Vert^2dt\\&\leq \frac{C}{\alpha}\int_0^{Lh}\left(\Vert \AA\mathbf{v}^M\Vert^2+\left\Vert D\partial_t\widetilde{\mathbf{v}}^M\right\Vert^2+\Vert \nabla\ww^M\Vert^2\Vert\ww^M\Vert^2_{\mathbf{H}^3(\Omega)}\Vert \nabla \mathbf{v}^M\Vert^2+\Vert \nabla\pphi^M_h\Vert_{\mathbf{L}^6(\Omega)}^4\Vert \nabla\mathbf{v}^M\Vert^2\right.\\&\left.+\Vert\nabla\mathbf{v}^M\Vert^6+\Vert \mathbf{w}^M\Vert_{\mathbf{H}^2(\Omega)}^2\Vert \nabla\pphi^M_h\Vert^2\right)ds,
 \end{align*}
 for any $L>0$. This implies, exploiting \eqref{unif1}, \eqref{ext}, \eqref{regmu}, \eqref{la} and \eqref{1b}, that 
 \begin{align}
 \Vert \partial_t\widetilde{\mathbf{v}}^M\Vert_{L^2(0,T;\mathbf{W}_\sigma)}\leq C(T,\alpha),\quad \forall T>0.
 \label{1c}
 \end{align}
 Following now the same arguments as in the proof of Theorem \ref{weakk}, we can pass to the limit as $M\to \infty$, obtaining the existence of a global weak solution $(\mathbf{v}_\alpha,\pphi_\alpha, \ww_\alpha)$ such that there exist $C>0$, independent of $\alpha$, and $C_\alpha>0$ depending on $\alpha$ such that
 \begin{align}
& 
 \sqrt{\alpha}\Vert\mathbf{v}_\alpha\Vert_{{L}^{\infty}(0,\infty;\mathbf{V}_\sigma)}+\Vert \mathbf{v}_\alpha\Vert_{L^\infty(0,\infty;\mathbf{H}_\sigma)}+ \label{reg1}\Vert \mathbf{v}_\alpha\Vert_{L^2(0,\infty;\mathbf{V}_\sigma)}\leq C,\\&\label{reg0}\Vert \pphi_\alpha\Vert_{L^\infty(0,\infty;\mathbf{H}^2(\Omega))}+\Vert \ww_\alpha\Vert_{L^\infty(0,\infty;\mathbf{H}^1(\Omega))}+\Vert \PsipA(\pphi_\alpha)\Vert_{L^\infty(0,\infty;\mathbf{L}^2(\Omega))}\leq C,\\&
\Vert \ww_\alpha\Vert_{L^2_{uloc}([0,\infty);\mathbf{H}^3(\Omega))}\leq C,
 \label{reg2}\\&
 \left\Vert {\mathbf{v}}_\alpha\right\Vert_{L^\infty(0,T;\mathbf{W}_\sigma)}+\left\Vert \partial_t{\mathbf{v}}_\alpha\right\Vert_{L^2(0,T;\mathbf{W}_\sigma)}\leq C_\alpha(T),\quad \forall T>0,
 \label{reg3} \\&
 \sum_{i=1}^N\varphi_{\alpha,i}\equiv 1,\quad \pphi_\alpha\in(0,1)^N\quad \text{a.e. in }\Omega\times(0,\infty).
 \end{align}
 The triple satisfies 
 \begin{align}
 &\mathbb{P}\left(\partial_t(\rho_\alpha\mathbf{v}_\alpha)-\alpha\Delta \partial_t\mathbf{v}_\alpha+\operatorname{div}(\rho_\alpha\mathbf{v}_\alpha\otimes\mathbf{v}_\alpha+\mathbf{v}_\alpha\otimes \mathbf{J}_{\rho,\alpha})-\operatorname{div}(2\nu(\pphi_\alpha)D\mathbf{v}_\alpha)\label{a2}-(\nabla\pphi_\alpha)^T\ww_\alpha\right)=\mathbf{0},\\&
 \partial_t\pphi_\alpha+(\nabla\pphi_\alpha)\mathbf{v}_\alpha=\operatorname{div}(\mathbf{M}\nabla\ww_\alpha),\\&
 \ww_\alpha= -\Delta\pphi_\alpha+\mathbf{P}(\Psi_{,\pphi}^1(\pphi_\alpha)-\mathbf{A}\pphi_\alpha),
 \end{align}
 almost everywhere on $\Omega\times(0,\infty)$, together with the initial conditions $(\mathbf{v}_\alpha,\pphi_\alpha)_{\vert t=0}=(\mathbf{v}_0^\alpha,\pphi_0)$ and the boundary conditions $\mathbf{v}=\mathbf{0}$, $\partial_\mathbf{n}\pphi=\partial_\mathbf{n}\ww=\mathbf{0}$ almost everywhere on $\partial\Omega\times(0,\infty)$.
  To pass to the limit in $\alpha$ we need to obtain a last estimate. In particular, by \eqref{reg3}, we can test \eqref{a2} with $\partial_t\mathbf{v}_\alpha$. This is possible since $\partial_t\mathbf{v}_\alpha\in L^2([0,T];\mathbf{W}_\sigma)$ for any $T>0$. Therefore we get
  \begin{align}
  &\nonumber\left(\partial_t(\rho_\alpha\mathbf{v}_\alpha),\partial_t\mathbf{v}_\alpha\right)_\Omega+\alpha\Vert D\partial_t\mathbf{v}_\alpha\Vert^2+\left(\operatorname{div}(\rho_\alpha\mathbf{v}_\alpha\otimes\mathbf{v}_\alpha),\partial_t\mathbf{v}_\alpha \right)_\Omega+\left(\operatorname{div}(\mathbf{v}_\alpha\otimes\RevA{\mathbf J_{\rho,\alpha}}),\partial_t\mathbf{v}_\alpha\right)_\Omega\\&+\left(2\nu(\pphi_\alpha)D\mathbf{v}_\alpha,D\partial_t\mathbf{v}_\alpha\right)_\Omega=\left((\nabla\pphi_\alpha)^T\ww_\alpha,\partial_t\mathbf{v}_\alpha\right)_\Omega.
  \label{o1}
  \end{align}    
  Now we have 
  \begin{align}
  \left(\partial_t(\rho_\alpha\mathbf{v}_\alpha),\partial_t\mathbf{v}_\alpha\right)_\Omega=\left(\rho_\alpha\partial_t\mathbf{v}_\alpha,\partial_t\mathbf{v}_\alpha\right)_\Omega+\left(\partial_t\rho_\alpha\mathbf{v}_\alpha,\partial_t\mathbf{v}_\alpha\right)_\Omega.
  \label{dec1}
  \end{align}
  Recalling Remark \ref{minrho}, we have
  $$
  \left(\rho_\alpha\partial_t\mathbf{v}_\alpha,\partial_t\mathbf{v}_\alpha\right)_\Omega\geq C_\rho\left\Vert\partial_t\mathbf{v}_\alpha\right\Vert^2,
  $$
  where $C_\rho=\min_{i=1,\ldots,N}\widetilde{\rho}_i>0$. 
  Observe that it holds, by Remark \ref{Jrho}, which is valid also at the continuous level,
  \begin{align*}
  &\left(\operatorname{div}(\rho_\alpha\mathbf{v}_\alpha\otimes\mathbf{v}_\alpha),\partial_t\mathbf{v}_\alpha \right)_\Omega+\left(\operatorname{div}(\mathbf{v}_\alpha\otimes\RevA{\mathbf J_{\rho,\alpha}}),\partial_t\mathbf{v}_\alpha\right)_\Omega+\left(\partial_t\rho_\alpha\mathbf{v}_\alpha,\partial_t\mathbf{v}_\alpha\right)_\Omega\\&=\left(\rho_\alpha(\mathbf{v}_\alpha\cdot\nabla)\mathbf{v}_\alpha,\partial_t\mathbf{v}_\alpha\right)_\Omega+\left((\RevA{\mathbf J_{\rho,\alpha}}\cdot\nabla)\mathbf{v}_\alpha,\partial_t\mathbf{v}_\alpha \right)_\Omega. 
  \end{align*}
  Now, by H\"{o}lder's and Young's inequalities together with 3D Sobolev-Gagliardo-Nirenberg's inequalities, we get, recalling that $\rho_\alpha$ is bounded uniformly in $\alpha$ due to $\pphi_\alpha\in[0,1]^N$,
  \begin{align*}
  &\left\vert\left(\rho_\alpha(\mathbf{v}_\alpha\cdot\nabla)\mathbf{v}_\alpha,\partial_t\mathbf{v}_\alpha\right)_\Omega\right\vert\leq C\Vert \mathbf{v}_\alpha\Vert_{\mathbf{L}^6(\Omega)}\Vert \nabla\mathbf{v}_\alpha\Vert_{\mathbf{L}^3(\Omega)}\left\Vert\partial_t\mathbf{v}_\alpha\right\Vert\\&\leq C\Vert \mathbf{v}_\alpha\Vert_{\mathbf{H}^2(\Omega)}^\frac 1 2\Vert\nabla \mathbf{v}_\alpha\Vert^{\frac 3 2}\left\Vert\partial_t\mathbf{v}_\alpha\right\Vert\leq \frac{C_\rho}8\left\Vert\partial_t\mathbf{v}_\alpha\right\Vert^2+\frac{\omega}{8}\Vert \mathbf{v}_\alpha\Vert_{\mathbf{H}^2(\Omega)}^2+C(\omega)\Vert\nabla\mathbf{v}_\alpha\Vert^6,
  \end{align*} 
  where $\omega>0$ is a positive constant to be specified later on. In the two-dimensional case we end up with 
  \begin{align*}
  &\left\vert\left(\rho_\alpha(\mathbf{v}_\alpha\cdot\nabla)\mathbf{v}_\alpha,\partial_t\mathbf{v}_\alpha\right)_\Omega\right\vert\leq C\Vert \mathbf{v}_\alpha\Vert_{\mathbf{L}^4(\Omega)}\Vert \nabla\mathbf{v}_\alpha\Vert_{\mathbf{L}^4(\Omega)}\left\Vert\partial_t\mathbf{v}_\alpha\right\Vert\\&\leq C\Vert \mathbf{v}_\alpha\Vert_{\mathbf{H}^2(\Omega)}^\frac 1 2\Vert \mathbf{v}_\alpha\Vert^{\frac 1 2}\Vert\nabla \mathbf{v}_\alpha\Vert\left\Vert\partial_t\mathbf{v}_\alpha\right\Vert\leq \frac{C_\rho}8\left\Vert\partial_t\mathbf{v}_\alpha\right\Vert^2+\frac{\omega}{8}\Vert \mathbf{v}_\alpha\Vert_{\mathbf{H}^2(\Omega)}^2+C(\omega)\Vert\nabla\mathbf{v}_\alpha\Vert^4\Vert \mathbf{v}_\alpha\Vert^2.
  \end{align*} 
  Then, concerning the other terms, recalling the definition of $\RevA{\mathbf J_{\rho,\alpha}}$, we get, in 3D,
  \begin{align*}
  &\left\vert\left((\RevA{\mathbf J_{\rho,\alpha}}\cdot\nabla)\mathbf{v}_\alpha,\partial_t\mathbf{v}_\alpha \right)_\Omega\right\vert\leq C\Vert \nabla\ww_\alpha\Vert_{\mathbf{L}^6(\Omega)}\Vert \nabla\mathbf{v}_\alpha\Vert_{\mathbf{L}^3(\Omega)}\left\Vert\partial_t\mathbf{v}_\alpha\right\Vert\\&\leq C\Vert \nabla\ww_\alpha\Vert^{\frac 1 2}\Vert \nabla\ww_\alpha\Vert_{\mathbf{H}^2(\Omega)}^\frac 1 2\Vert \nabla\mathbf{v}_\alpha\Vert^\frac 1 2 \Vert \mathbf{v}_\alpha\Vert^\frac 1 2_{\mathbf{H}^2(\Omega)}\left\Vert\partial_t\mathbf{v}_\alpha\right\Vert\\&\leq\frac {C_\rho} 8\left\Vert\partial_t\mathbf{v}_\alpha\right\Vert^2+\frac{\omega}{8}\Vert \mathbf{v}_\alpha\Vert^2_{\mathbf{H}^2(\Omega)}+C(\omega)\Vert \nabla\ww_\alpha\Vert^2\Vert \nabla\ww_\alpha\Vert_{\mathbf{H}^2(\Omega)}^2\Vert \nabla \mathbf{v}_\alpha\Vert^2.
  \end{align*}
  Here we exploited the 3D inequality \eqref{pp}.
  In the 2D case we have, by Agmon's inequality,
  \begin{align*}
  &\left\vert\left((\RevA{\mathbf J_{\rho,\alpha}}\cdot\nabla)\mathbf{v}_\alpha,\partial_t\mathbf{v}_\alpha \right)_\Omega\right\vert\leq C\Vert \nabla\ww_\alpha\Vert_{\mathbf{L}^\infty(\Omega)}\Vert \nabla\mathbf{v}_\alpha\Vert\left\Vert\partial_t\mathbf{v}_\alpha\right\Vert\\&\leq C\Vert \nabla\ww_\alpha\Vert^{\frac 1 2}\Vert \nabla\ww_\alpha\Vert_{\mathbf{H}^2(\Omega)}^\frac 1 2\Vert \nabla\mathbf{v}_\alpha\Vert\left\Vert\partial_t\mathbf{v}_\alpha\right\Vert\\&\leq \frac {C_\rho} 8\left\Vert\partial_t\mathbf{v}_\alpha\right\Vert^2+C(\omega)\Vert \nabla\ww_\alpha\Vert\Vert \nabla\ww_\alpha\Vert_{\mathbf{H}^2(\Omega)}\Vert \nabla \mathbf{v}_\alpha\Vert^2.
  \end{align*}
We then have
  \begin{align*}
  &\left(2\nu(\pphi_\alpha)D\mathbf{v}_\alpha,D\partial_t\mathbf{v}_\alpha\right)_\Omega=\ddt(\nu(\pphi_\alpha)D\mathbf{v}_\alpha,D\mathbf{v}_\alpha)-2((\nu_{,\pphi}(\pphi_\alpha)\cdot \partial_t\pphi_\alpha)D\mathbf{v}_\alpha,D\mathbf{v}_\alpha).
  \end{align*}
  This means that we only need the estimate, in 3D,
  \begin{align*}
  &\left\vert2((\nu_{,\pphi}(\pphi_\alpha)\cdot \partial_t\pphi_\alpha)D\mathbf{v}_\alpha,D\mathbf{v}_\alpha)\right\vert \leq C\left\Vert\partial_t\pphi_\alpha\right\Vert_{\mathbf{L}^6(\Omega)}\Vert D\mathbf{v}_\alpha\Vert\Vert D\mathbf{v}_\alpha\Vert_{\mathbf{L}^3(\Omega)}\\&\leq C\left\Vert\nabla\partial_t\pphi_\alpha\right\Vert\Vert \nabla \mathbf{v}_\alpha\Vert^\frac 3 2\Vert \mathbf{v}_\alpha\Vert_{\mathbf{H}^2(\Omega)}^\frac 1 2\leq {C}\left\Vert\nabla\partial_t\pphi_\alpha\right\Vert^2+\frac{\omega}{8}\Vert \mathbf{v}_\alpha\Vert_{\mathbf{H}^2(\Omega)}^2+C(\omega)\Vert\nabla\mathbf{v}_\alpha\Vert^6,
  \end{align*}
  where we exploited the assumption that $\nu_{,\pphi}$ is bounded uniformly. In 2D this becomes, 
  \begin{align*}
  &\left\vert((\nu_{,\pphi}(\pphi_\alpha)\cdot \partial_t\pphi_\alpha)D\mathbf{v}_\alpha,D\mathbf{v}_\alpha)\right\vert \leq C\left\Vert\partial_t\pphi_\alpha\right\Vert\Vert D\mathbf{v}_\alpha\Vert_{\mathbf{L}^4(\Omega)}^2\\&\leq C\left\Vert\nabla\partial_t\pphi_\alpha\right\Vert\Vert \nabla \mathbf{v}_\alpha\Vert\Vert \mathbf{v}_\alpha\Vert_{\mathbf{H}^2(\Omega)}\leq \frac{\omega}{8}\Vert \mathbf{v}_\alpha\Vert_{\mathbf{H}^2(\Omega)}^2+C(\omega)\left\Vert\nabla\partial_t\pphi_\alpha\right\Vert^2\Vert\nabla\mathbf{v}_\alpha\Vert^2.
  \end{align*}
  We are only left with one term: we have in 3D, by Agmon's inequality,
  \begin{align*}
  &\left\vert\left((\nabla\pphi_\alpha)^T\ww_\alpha,\partial_t\mathbf{v}_\alpha\right)_\Omega\right\vert\leq \Vert \nabla\pphi_\alpha\Vert\Vert\ww_\alpha\Vert_{\mathbf{L}^\infty(\Omega)}\left\Vert\partial_t\mathbf{v}_\alpha\right\Vert\leq C\Vert \nabla\pphi_\alpha\Vert\Vert\ww_\alpha\Vert_{\mathbf{H}^1(\Omega)}^\frac 1 2\Vert \ww_\alpha\Vert_{\mathbf{H}^2(\Omega)}^\frac 1 2\left\Vert\partial_t\mathbf{v}_\alpha\right\Vert\\&\leq \frac {C_\rho}{8}\left\Vert\partial_t\mathbf{v}_\alpha\right\Vert^2+C\Vert \ww_{\RevA{\alpha}}\Vert_{\mathbf{H}^2(\Omega)}^2+C \Vert \nabla\pphi_\alpha\Vert^4\Vert\ww_\alpha\Vert_{\mathbf{H}^1(\Omega)}^2,
  \end{align*}
  whereas analogously in 2D, by the embedding $\mathbf{H}^1(\Omega)\hookrightarrow \mathbf{L}^2(\Omega)$,
  \begin{align*}
  &\left\vert\left((\nabla\pphi_\alpha)^T\ww_\alpha,\partial_t\mathbf{v}_\alpha\right)_\Omega\right\vert\leq \Vert \nabla\pphi_\alpha\Vert\Vert\ww_\alpha\Vert_{\mathbf{L}^\infty(\Omega)}\left\Vert\partial_t\mathbf{v}_\alpha\right\Vert\leq C\Vert \nabla\pphi_\alpha\Vert\Vert\ww_\alpha\Vert^\frac 1 2\Vert \ww_\alpha\Vert_{\mathbf{H}^2(\Omega)}^\frac 1 2\left\Vert\partial_t\mathbf{v}_\alpha\right\Vert\\&\leq \frac {C_\rho}{8}\left\Vert\partial_t\mathbf{v}_\alpha\right\Vert^2+C\Vert \ww_\alpha\Vert_{\mathbf{H}^2(\Omega)}^2+C \Vert \nabla\pphi_\alpha\Vert^4\Vert\ww_\alpha\Vert_{\mathbf{H}^1(\Omega)}^2.
  \end{align*}
  To sum up, coming back to \eqref{o1}, in the end we obtain in 3D, recalling that $\nu\geq \nu_*>0$, together with the uniform estimates above,  
  \begin{align}
  \nonumber&\ddt\left(\nu({\pphi_\alpha})D\mathbf{v}_\alpha,D\mathbf{v}_\alpha\right)_\Omega+\alpha\Vert D\partial_t\mathbf{v}_\alpha\Vert^2+\frac{C_\rho}{2}\left\Vert\partial_t\mathbf{v}_\alpha\right\Vert^2\\&\leq C(\omega)\left(\left(\nu({\pphi_\alpha})D\mathbf{v}_\alpha,D\mathbf{v}_\alpha\right)_\Omega\right)_\Omega^3+\Vert \nabla\ww_\alpha\Vert^2\Vert \nabla\ww_\alpha\Vert_{\mathbf{H}^2(\Omega)}^2\left(\nu({\pphi_\alpha})D\mathbf{v}_\alpha,D\mathbf{v}_\alpha\right)_\Omega\nonumber\\&+ C\left\Vert\nabla\partial_t\pphi_\alpha\right\Vert^2+C\Vert \ww_\alpha\Vert_{\mathbf{H}^2(\Omega)}^2+C \Vert \nabla\pphi_\alpha\Vert^4\Vert\ww_\alpha\Vert_{\mathbf{H}^1(\Omega)}^2+\frac{3\omega}{8}\Vert \mathbf{v}_\alpha\Vert_{\mathbf{H}^2(\Omega)}^2,
  \label{dt}
  \end{align}
  whereas in 2D we have
  \begin{align}
  \nonumber&\ddt\left(\nu({\pphi_\alpha})D\mathbf{v}_\alpha,D\mathbf{v}_\alpha\right)_\Omega+\alpha\Vert D\partial_t\mathbf{v}_\alpha\Vert^2+\frac{C_\rho}{2}\left\Vert\partial_t\mathbf{v}_\alpha\right\Vert^2\\&\leq \nonumber C(\omega)\Vert\mathbf{v}_\alpha\Vert^2\left(\nu({\pphi_\alpha})D\mathbf{v}_\alpha,D\mathbf{v}_\alpha\right)_\Omega^2+\Vert \nabla\ww_\alpha\Vert\Vert \nabla\ww_\alpha\Vert_{\mathbf{H}^2(\Omega)}\left(\nu({\pphi_\alpha})D\mathbf{v}_\alpha,D\mathbf{v}_\alpha\right)_\Omega\\&+ C\left\Vert\nabla\partial_t\pphi_\alpha\right\Vert^2\left(\nu({\pphi_\alpha})D\mathbf{v}_\alpha,D\mathbf{v}_\alpha\right)_\Omega+C\Vert \ww_\alpha\Vert_{\mathbf{H}^2(\Omega)}^2+C \Vert \nabla\pphi_\alpha\Vert^4\Vert\ww_\alpha\Vert_{\mathbf{H}^1(\Omega)}^2+\frac{\omega}{4}\Vert \mathbf{v}_\alpha\Vert_{\mathbf{H}^2(\Omega)}^2.
  \label{dt2}
  \end{align}
  We only need to find a uniform $\mathbf{H}^2$-control on $\mathbf{v}_\alpha$. Let us test \eqref{a2} with $\AA\mathbf{v}_\alpha\in \mathbf{H}_\sigma$ and integrate over $\Omega$. We get 
  \begin{align}
  &\nonumber(\partial_t(\rho_\alpha\mathbf{v}_\alpha),\AA\mathbf{v}_\alpha)_\Omega+\frac{\alpha}{2}\ddt\Vert \AA\mathbf{v}_\alpha\Vert^2+(\operatorname{div}(\rho_\alpha\mathbf{v}_\alpha\otimes\mathbf{v}_\alpha),\AA\mathbf{v}_\alpha)_\Omega\\&+(\operatorname{div}(\mathbf{v}_\alpha\otimes \RevA{\mathbf J_{\rho,\alpha}}),\AA\mathbf{v}_\alpha)_\Omega-(2\operatorname{div}(\nu(\pphi_\alpha)D\mathbf{v}_\alpha),\AA\mathbf{v}_\alpha)_\Omega=((\nabla\pphi_\alpha)^T\ww_\alpha,\AA\mathbf{v}_\alpha)_\Omega.
  \label{h2}
  \end{align}
  Notice that this has been possible since $t\mapsto \Vert \AA\mathbf{v}_\alpha(t)\Vert^2$ is absolutely continuous on $[0,T]$, for any $T>0$.
  First, observe that 
  \begin{align*}
  &-(2\operatorname{div}(\nu(\pphi_\alpha)D\mathbf{v}_\alpha),\AA\mathbf{v}_\alpha)_\Omega\\&= -(2\nu(\pphi_\alpha)\Delta\mathbf{v}_\alpha,\AA\mathbf{v}_\alpha)_\Omega -(2D\mathbf{v}_\alpha\nu_{,\pphi}(\pphi_\alpha)\nabla\pphi_\alpha,\AA\mathbf{v}_\alpha)_\Omega\\&=(2\nu(\pphi_\alpha)\AA\mathbf{v}_\alpha,\AA\mathbf{v}_\alpha)_\Omega-(2\nu(\pphi_\alpha)\nabla p_\alpha,\AA\mathbf{v}_\alpha)_\Omega-(2D\mathbf{v}_\alpha\nu_{,\pphi}(\pphi_\alpha)\nabla\pphi_\alpha,\AA\mathbf{v}_\alpha)_\Omega, 
  \end{align*}
  where $p_\alpha\in H^1(\Omega)$ ($\overline{p}_\alpha=0$) satisfies $-\Delta\mathbf{v}_\alpha+\nabla p_\alpha=\AA\mathbf{v}_\alpha$ in $\Omega$. Observe that, since $\nu\geq \nu_*>0$,
  $$
  (2\nu(\pphi_\alpha)\AA\mathbf{v}_\alpha,\AA\mathbf{v}_\alpha)_\Omega\geq 2\nu_*\Vert \AA\mathbf{v}_\alpha\Vert^2.
  $$
  Moreover, by standard inequalities, being $\nu_{,\pphi}$ bounded, in 3D
  \begin{align*}
 & \vert(2D\mathbf{v}_\alpha\nu_{,\pphi}(\pphi_\alpha)\nabla\pphi_\alpha,\AA\mathbf{v}_\alpha)_\Omega\vert\leq C\Vert D\mathbf{v}_\alpha\Vert_{\mathbf{L}^3(\Omega)}\Vert \nabla\pphi_\alpha\Vert_{\mathbf{L}^6(\Omega)}\Vert \AA\mathbf{v}_\alpha\Vert\\&\leq C\Vert D\mathbf{v}_\alpha\Vert^\frac 1 2 \Vert \pphi_\alpha\Vert_{\mathbf{H}^2(\Omega)}\Vert \AA\mathbf{v}_\alpha\Vert^\frac 3 2\leq \frac{\nu_*}{4}\Vert \AA\mathbf{v}_\alpha\Vert^2+C \Vert \pphi_\alpha\Vert_{\mathbf{H}^2(\Omega)}^4\Vert D\mathbf{v}_\alpha\Vert^2,
  \end{align*}
 whereas in 2D
 \begin{align*}
  & \vert(2D\mathbf{v}_\alpha\nu_{,\pphi}(\pphi_\alpha)\nabla\pphi_\alpha,\AA\mathbf{v}_\alpha)_\Omega\vert\leq C\Vert D\mathbf{v}_\alpha\Vert_{\mathbf{L}^4(\Omega)}\Vert \nabla\pphi_\alpha\Vert_{\mathbf{L}^4(\Omega)}\Vert \AA\mathbf{v}_\alpha\Vert\\&\leq C\Vert D\mathbf{v}_\alpha\Vert^\frac 1 2 \Vert \pphi_\alpha\Vert_{\mathbf{H}^2(\Omega)}^\frac 1 2\Vert \nabla\pphi_\alpha\Vert^\frac 1 2\Vert \AA\mathbf{v}_\alpha\Vert^\frac 3 2\leq \frac{\nu_*}{4}\Vert \AA\mathbf{v}_\alpha\Vert^2+C \Vert \pphi_\alpha\Vert_{\mathbf{H}^2(\Omega)}^2\Vert \nabla\pphi_\alpha\Vert^2\Vert D\mathbf{v}_\alpha\Vert^2.
 \end{align*}
  Then we can write
 \RevA{ \begin{align*}
&(2\nu(\pphi_\alpha)\nabla p_\alpha,\AA\mathbf{v}_\alpha)_\Omega\\&=(\nabla(2\nu(\pphi_\alpha)p_\alpha), \AA\mathbf{v}_\alpha)_\Omega-(2(\nabla\pphi_\alpha)^T\nu_{,\pphi}(\pphi_\alpha)p_\alpha, \AA\mathbf{v}_\alpha)_\Omega\\&=-(2(\nabla\pphi_\alpha)^T\nu_{,\pphi}(\pphi_\alpha)p_\alpha, \AA\mathbf{v}_\alpha)_\Omega,
  \end{align*}}
  exploiting the fact that $\AA\mathbf{v}_\alpha\in \mathbf{H}_\sigma$ for any $t\geq0$, so that $\AA\mathbf{v}_\alpha\perp_{\mathbf{L}^2(\Omega)}\nabla(2\nu(\pphi_\alpha)p_\alpha)$, \RevA{since we have that } $2\nu(\pphi_\alpha)p_\alpha\in H^1(\Omega)$.
  Therefore, in 3D,
  \begin{align*}
  &\vert(2\nu(\pphi_\alpha)\nabla p_\alpha,\AA\mathbf{v}_\alpha)_\Omega\vert=\vert(2(\nabla\pphi_\alpha)^T\nu_{,\pphi}(\pphi_\alpha)p_\alpha, \AA\mathbf{v}_\alpha)_\Omega\vert\\&\leq C\Vert \nabla\pphi_\alpha\Vert_{\mathbf{L}^6(\Omega)}\Vert p_\alpha\Vert_{L^3(\Omega)}\Vert \AA\mathbf{v}_\alpha\Vert\leq C\Vert \pphi_\alpha\Vert_{\mathbf{H}^2(\Omega)}\Vert D\mathbf{v}_\alpha\Vert^\frac 1 2\Vert \AA\mathbf{v}_\alpha\Vert^\frac 3 2 \\&\leq 
  \frac{\nu_*}{4}\Vert \AA\mathbf{v}_\alpha\Vert^2+C\Vert \pphi_\alpha\Vert_{\mathbf{H}^2(\Omega)}^4\Vert D\mathbf{v}_\alpha\Vert^2,
  \end{align*}
  where we exploited \eqref{pr} with $\mathbf f=\AA\mathbf{v}_\alpha$. Analogously, in 2D we infer, again by \eqref{pr},
   \begin{align*}
  &\vert(2\nu(\pphi_\alpha)\nabla p_\alpha,\AA\mathbf{v}_\alpha)_\Omega\vert=\vert(2(\nabla\pphi_\alpha)^T\nu_{,\pphi}(\pphi_\alpha)p_\alpha, \AA\mathbf{v}_\alpha)_\Omega\vert\\&\leq C\Vert \nabla\pphi_\alpha\Vert_{\mathbf{L}^4(\Omega)}\Vert p_\alpha\Vert_{L^4(\Omega)}\Vert \AA\mathbf{v}_\alpha\Vert\leq C\Vert \pphi_\alpha\Vert_{\mathbf{H}^2(\Omega)}^\frac 1 2\Vert \nabla\pphi_\alpha\Vert^\frac 1 2\Vert D\mathbf{v}_\alpha\Vert^\frac 1 2\Vert \AA\mathbf{v}_\alpha\Vert^\frac 3 2 \\&\leq 
  \frac{\nu_*}{4}\Vert \AA\mathbf{v}_\alpha\Vert^2+C\Vert \pphi_\alpha\Vert_{\mathbf{H}^2(\Omega)}^2\Vert\nabla\pphi_\alpha\Vert^2\Vert D\mathbf{v}_\alpha\Vert^2.
  \end{align*}
  Proceeding in the estimates, notice that
  $$
  (\partial_t(\rho_\alpha\mathbf{v}_\alpha),\AA\mathbf{v}_\alpha)_\Omega=(\partial_t\rho_\alpha\mathbf{v}_\alpha,\AA\mathbf{v}_\alpha)_\Omega+(\rho_\alpha\partial_t\mathbf{v}_\alpha,\AA\mathbf{v}_\alpha)_\Omega,
  $$
  so that, recalling that $\rho_\alpha=\widetilde{\boldsymbol\rho}\cdot \pphi_\alpha$,
  \begin{align*}
  &\left\vert (\rho_\alpha\partial_t\mathbf{v}_\alpha,\AA\mathbf{v}_\alpha)\right\vert\leq C\Vert \partial_t\mathbf{v}_\alpha\Vert\Vert\AA\mathbf{v}_\alpha\Vert\leq \frac{\nu_*}{4}\Vert\AA\mathbf{v}_\alpha\Vert^2+C_0\Vert \partial_t\mathbf{v}_\alpha\Vert^2,
  \end{align*}
  where $C_0>0$ does not depend on $\alpha$. Furthermore, we have, as before, 
   \begin{align*}
  &\left(\operatorname{div}(\rho_\alpha\mathbf{v}_\alpha\otimes\mathbf{v}_\alpha),\AA\mathbf{v}_\alpha \right)_\Omega+\left(\operatorname{div}(\mathbf{v}_\alpha\otimes\RevA{\mathbf J_{\rho,\alpha}}),\AA\mathbf{v}_\alpha\right)_\Omega+\left(\partial_t\rho_\alpha\mathbf{v}_\alpha,\AA\mathbf{v}_\alpha\right)_\Omega\\&=\left(\rho_\alpha(\mathbf{v}_\alpha\cdot\nabla)\mathbf{v}_\alpha,\AA\mathbf{v}_\alpha\right)_\Omega+\left((\RevA{\mathbf J_{\rho,\alpha}}\cdot\nabla)\mathbf{v}_\alpha,\AA\mathbf{v}_\alpha \right)_\Omega.
  \end{align*}
    Now, by H\"{o}lder's and Young's inequalities together with 3D Sobolev-Gagliardo-Nirenberg's inequalities, we get, recalling that $\rho_\alpha$ is bounded uniformly in $\alpha$, since $\pphi_\alpha\in[0,1]^N$,
  \begin{align*}
  &\left\vert\left(\rho_\alpha(\mathbf{v}_\alpha\cdot\nabla)\mathbf{v}_\alpha,\AA\mathbf{v}_\alpha\right)_\Omega\right\vert\leq C\Vert \mathbf{v}_\alpha\Vert_{\mathbf{L}^6(\Omega)}\Vert \nabla\mathbf{v}_\alpha\Vert_{\mathbf{L}^3(\Omega)}\left\Vert\AA\mathbf{v}_\alpha\right\Vert\\&\leq C\Vert\nabla \mathbf{v}_\alpha\Vert^{\frac 3 2}\left\Vert\AA\mathbf{v}_\alpha\right\Vert^\frac 3 2\leq \frac{\nu_*}4\left\Vert\AA\mathbf{v}_\alpha\right\Vert^2+C\Vert\nabla\mathbf{v}_\alpha\Vert^6.
  \end{align*} 
  In the two-dimensional case we end up with  
  \begin{align*}
  &\left\vert\left(\rho_\alpha(\mathbf{v}_\alpha\cdot\nabla)\mathbf{v}_\alpha,\AA\mathbf{v}_\alpha\right)_\Omega\right\vert\leq C\Vert \mathbf{v}_\alpha\Vert_{\mathbf{L}^4(\Omega)}\Vert \nabla\mathbf{v}_\alpha\Vert_{\mathbf{L}^4(\Omega)}\left\Vert\AA\mathbf{v}_\alpha\right\Vert\\&\leq C\Vert \mathbf{v}_\alpha\Vert^{\frac 1 2}\Vert\nabla \mathbf{v}_\alpha\Vert\left\Vert\AA\mathbf{v}_\alpha\right\Vert^{\frac 3 2 }\leq \frac{\nu_*}4\left\Vert\AA\mathbf{v}_\alpha\right\Vert^2+C\Vert\nabla\mathbf{v}_\alpha\Vert^4\Vert \mathbf{v}_\alpha\Vert^2.
  \end{align*} 
  Then, concerning the other terms, recalling the definition of $\RevA{\mathbf J_{\rho,\alpha}}$, we get, in 3D,
  \begin{align*}
  &\left\vert\left((\RevA{\mathbf J_{\rho,\alpha}}\cdot\nabla)\mathbf{v}_\alpha,\AA\mathbf{v}_\alpha \right)_\Omega\right\vert\leq C\Vert \nabla\ww_\alpha\Vert_{\mathbf{L}^6(\Omega)}\Vert \nabla\mathbf{v}_\alpha\Vert_{\mathbf{L}^3(\Omega)}\left\Vert\AA\mathbf{v}_\alpha\right\Vert\\&\leq C\Vert \nabla\ww_\alpha\Vert^{\frac 1 2}\Vert \nabla\ww_\alpha\Vert_{\mathbf{H}^2(\Omega)}^\frac 1 2\Vert \nabla\mathbf{v}_\alpha\Vert^\frac 1 2\left\Vert\AA\RevA{\mathbf{v}_\alpha}\right\Vert^\frac 3 2\\&\leq\frac {\nu_*} 4\left\Vert\AA\RevA{\mathbf{v}_\alpha}\right\Vert^2+C\Vert \nabla\ww_\alpha\Vert^2\Vert \nabla\ww_\alpha\Vert_{\mathbf{H}^2(\Omega)}^2\Vert \nabla \mathbf{v}_\alpha\Vert^2.
  \end{align*}
  Here we exploited again the 3D inequality \eqref{pp}.
  In the 2D case we have, by Agmon's inequality,
  \begin{align*}
  &\left\vert\left((\RevA{\mathbf J_{\rho,\alpha}}\cdot\nabla)\mathbf{v}_\alpha,\AA\RevA{\mathbf{v}_\alpha} \right)_\Omega\right\vert\leq C\Vert \nabla\ww_\alpha\Vert_{\mathbf{L}^\infty(\Omega)}\Vert \nabla\mathbf{v}_\alpha\Vert\left\Vert\AA\RevA{\mathbf{v}_\alpha}\right\Vert\\&\leq C\Vert \nabla\ww_\alpha\Vert^{\frac 1 2}\Vert \nabla\ww_\alpha\Vert_{\mathbf{H}^2(\Omega)}^\frac 1 2\Vert \nabla\mathbf{v}_\alpha\Vert\left\Vert\AA\RevA{\mathbf{v}_\alpha}\right\Vert\\&\leq \frac {\nu_*} 4\left\Vert\AA\RevA{\mathbf{v}_\alpha}\right\Vert^2+C\Vert \nabla\ww_\alpha\Vert\Vert \nabla\ww_\alpha\Vert_{\mathbf{H}^2(\Omega)}\Vert \nabla \mathbf{v}_\alpha\Vert^2.
  \end{align*}
  In conclusion, in 3D, by Agmon's inequality,
  \begin{align*}
  &\left\vert\left((\nabla\pphi_\alpha)^T\ww_\alpha,\AA\RevA{\mathbf{v}_\alpha}\right)_\Omega\right\vert\leq \Vert \nabla\pphi_\alpha\Vert\Vert\ww_\alpha\Vert_{\mathbf{L}^\infty(\Omega)}\left\Vert\AA\RevA{\mathbf{v}_\alpha}\right\Vert\leq C\Vert \nabla\pphi_\alpha\Vert\Vert\ww_\alpha\Vert_{\mathbf{H}^1(\Omega)}^\frac 1 2\Vert \ww_\alpha\Vert_{\mathbf{H}^2(\Omega)}^\frac 1 2\left\Vert\AA\RevA{\mathbf{v}_\alpha}\right\Vert\\&\leq \frac {\nu_*}{4}\left\Vert\AA\RevA{\mathbf{v}_\alpha}\right\Vert^2+C \Vert \nabla\pphi_\alpha\Vert^4\Vert\ww_\alpha\Vert_{\mathbf{H}^1(\Omega)}^2+C\Vert \ww_\alpha\Vert_{\mathbf{H}^2(\Omega)}^2,
  \end{align*}
  whereas, analogously, in 2D,
  \begin{align*}
  &\left\vert\left((\nabla\pphi_\alpha)^T\ww_\alpha,\AA\RevA{\mathbf{v}_\alpha}\right)_\Omega\right\vert\leq \Vert \nabla\pphi_\alpha\Vert\Vert\ww_\alpha\Vert_{\mathbf{L}^\infty(\Omega)}\left\Vert\AA\RevA{\mathbf{v}_\alpha}\right\Vert\leq C\Vert \nabla\pphi_\alpha\Vert\Vert\ww_\alpha\Vert^\frac 1 2\Vert \ww_\alpha\Vert_{\mathbf{H}^2(\Omega)}^\frac 1 2\left\Vert\AA\RevA{\mathbf{v}_\alpha}\right\Vert\\&\leq \frac {\nu_*}{4}\left\Vert\AA\RevA{\mathbf{v}_\alpha}\right\Vert^2+C\Vert \ww_\alpha\Vert_{\mathbf{H}^2(\Omega)}^2+C \Vert \nabla\pphi_\alpha\Vert^4\Vert\ww_\alpha\Vert_{\mathbf{H}^1(\Omega)}^2.
  \end{align*}
  To sum up, we have obtained, in 3D,
  \begin{align}
  &\nonumber\frac \alpha 2\frac{d}{dt}\Vert \AA\mathbf{v}_\alpha\Vert^2+\frac{\nu_*}2\Vert \AA\mathbf{v}_\alpha\Vert^2\leq C\left(\left(\nu({\pphi_\alpha})D\mathbf{v}_\alpha,D\mathbf{v}_\alpha\right)_\Omega\right)^3\\&+C(\Vert \nabla\ww_\alpha\Vert^2\Vert \nabla\ww_\alpha\Vert_{\mathbf{H}^2(\Omega)}^2+\Vert \pphi_\alpha\Vert_{\mathbf{H}^2(\Omega)}^4)\left(\nu({\pphi_\alpha})D\mathbf{v}_\alpha,D\mathbf{v}_\alpha\right)_\Omega+ C\left\Vert\nabla\partial_t\pphi_\alpha\right\Vert^2\nonumber\\&+C\Vert \ww_\alpha\Vert_{\mathbf{H}^2(\Omega)}^2+C \Vert \nabla\pphi_\alpha\Vert^4\Vert\ww_\alpha\Vert_{\mathbf{H}^1(\Omega)}^2+C_0\Vert \partial_t\mathbf{v}_\alpha\Vert^2,
  \label{3d}
  \end{align}
  whereas, in 2D,
  \begin{align}
  \nonumber&\frac{\alpha}{2}\ddt\left\Vert\AA\RevA{\mathbf{v}_\alpha}\right\Vert^2+\frac{\nu_*}2\Vert \AA\mathbf{v}_\alpha\Vert^2\\&\leq \nonumber C\Vert\mathbf{v}_\alpha\Vert^2\left(\nu({\pphi_\alpha})D\mathbf{v}_\alpha,D\mathbf{v}_\alpha\right)_\Omega^2+C(\Vert \nabla\ww_\alpha\Vert\Vert \nabla\ww_\alpha\Vert_{\mathbf{H}^2(\Omega)}+\Vert \pphi_\alpha\Vert_{\mathbf{H}^2(\Omega)}^2\Vert\nabla\pphi_\alpha\Vert^2)\left(\nu({\pphi_\alpha})D\mathbf{v}_\alpha,D\mathbf{v}_\alpha\right)_\Omega\\&+ C\left\Vert\nabla\partial_t\pphi_\alpha\right\Vert^2\left(\nu({\pphi_\alpha})D\mathbf{v}_\alpha,D\mathbf{v}_\alpha\right)_\Omega+C\Vert \ww_\alpha\Vert_{\mathbf{H}^2(\Omega)}^2+C \Vert \nabla\pphi_\alpha\Vert^4\Vert\ww_\alpha\Vert_{\mathbf{H}^1(\Omega)}^2+C_0\Vert \partial_t\mathbf{v}_\alpha\Vert^2.
  \label{2d}
  \end{align}
  Now, in the 3D case we can sum \eqref{dt} with \eqref{3d} multiplied by $\frac{C_\rho}{4C_0}$ and fix $\omega=\frac{C_\rho\nu_*}{6C_0}$, to get, exploiting \eqref{reg1}-\eqref{reg2},
    \begin{align}
  &\nonumber\frac{d}{dt}\left(\frac\alpha 2\Vert \AA\mathbf{v}_\alpha\Vert^2+\left(\nu({\pphi_\alpha})D\mathbf{v}_\alpha,D\mathbf{v}_\alpha\right)_\Omega\right)+\alpha\Vert D\partial_t\mathbf{v}_\alpha\Vert^2+\frac{\nu_*C_\rho}{16C_0}\Vert \AA\mathbf{v}_\alpha\Vert^2+\frac{C_\rho}{4}\left\Vert\partial_t\mathbf{v}_\alpha\right\Vert^2\\&\leq C\left(\left(\nu({\pphi_\alpha})D\mathbf{v}_\alpha,D\mathbf{v}_\alpha\right)_\Omega\right)^3+C(\Vert \nabla\ww_\alpha\Vert_{\mathbf{H}^2(\Omega)}^2+1)\left(\nu({\pphi_\alpha})D\mathbf{v}_\alpha,D\mathbf{v}_\alpha\right)_\Omega\nonumber\\&+ C\left\Vert\nabla\partial_t\pphi_\alpha\right\Vert^2+C\Vert \ww_\alpha\Vert_{\mathbf{H}^2(\Omega)}^2+C,
  \label{3db}
  \end{align}
  where all the constants $C$ appearing on the left-hand side do not depend on $\alpha$.
  
 In the 2D case we can sum \eqref{dt2} with \eqref{2d} multiplied by $\frac{C_\rho}{4C_0}$ and fix $\omega= \frac{C_\rho\nu_*}{4C_0}$ to get, in the end, by means of \eqref{reg1}-\eqref{reg2},
   \begin{align}
 \nonumber&\frac{d}{dt}\left(\frac\alpha 2\Vert \AA\mathbf{v}_\alpha\Vert^2+\left(\nu({\pphi_\alpha})D\mathbf{v}_\alpha,D\mathbf{v}_\alpha\right)_\Omega\right)+\alpha\Vert D\partial_t\mathbf{v}_\alpha\Vert^2+\frac{\nu_*C_\rho}{16C_0}\Vert \AA\mathbf{v}_\alpha\Vert^2+\frac{C_\rho}{4}\left\Vert\partial_t\mathbf{v}_\alpha\right\Vert^2\\&\leq \nonumber C\left(\nu({\pphi_\alpha})D\mathbf{v}_\alpha,D\mathbf{v}_\alpha\right)_\Omega^2+C(1+\Vert \nabla\ww_\alpha\Vert_{\mathbf{H}^2(\Omega)})\left(\nu({\pphi_\alpha})D\mathbf{v}_\alpha,D\mathbf{v}_\alpha\right)_\Omega\\&+ C\left\Vert\nabla\partial_t\pphi_\alpha\right\Vert^2\left(\nu({\pphi_\alpha})D\mathbf{v}_\alpha,D\mathbf{v}_\alpha\right)_\Omega+C\Vert \ww_\alpha\Vert_{\mathbf{H}^2(\Omega)}^2+C.
 \label{2db}
 \end{align}
 Therefore, in the 3D case, we have also:
 \begin{align}
  u(t)\leq u(0)+\int_0^t f(s)w(u(s))ds,
 \label{En}
 \end{align}
 where 
 $$
 u:=\frac\alpha 2\Vert \AA\mathbf{v}_\alpha\Vert^2+\left(\nu({\pphi_\alpha})D\mathbf{v}_\alpha,D\mathbf{v}_\alpha\right)_\Omega
 $$
 $$
 f:=C(1+\Vert \nabla\ww_\alpha\Vert_{\mathbf{H}^2(\Omega)}^2+\left\Vert\nabla\partial_t\pphi_\alpha\right\Vert^2+\Vert \ww_{\RevA{\alpha}}\Vert_{\mathbf{H}^2(\Omega)}^2),
 $$
 $$
 w(s):=(1+s)^3.
 $$
  Now notice that, by \eqref{reg1}-\eqref{reg2} and \eqref{appr},
 \begin{align}
 u(0)=\frac\alpha 2\Vert \AA\mathbf{v}_\alpha(0)\Vert^2+\left(\nu(\pphi_\alpha(0))D\mathbf{v}_\alpha(0),D\mathbf{v}_\alpha(0)\right)_\Omega\leq C_A,
 \label{C_A}
 \end{align}
 uniformly in $\alpha$, as well as 
 $$
 \int_0^T f(s)ds\leq C_B(T),
 $$
 independently of $\alpha$. Notice that $C_B$ is such that $C_B(T)\to0$ as $T\to0$, but it also increases with $T$.
 
In the notation of Lemma \ref{bihari}, we also set
$$
G(x):=\int_{1}^x\frac{dy}{w(y)}dy= \frac 1 8 -\frac{1}{2(1+x)^2},\quad x\geq0,
$$
together with 
$$
G^{-1}(y)=\dfrac{1}{\sqrt{\frac 1 4 -2y}}-1,\quad -\frac 3 8<y<\frac1 8.
$$
Since $G^{-1}$ is increasing on its domain of definition, we deduce by Lemma \ref{bihari} that 
$$
\max_{t\in[0,T]}u(t)\leq G^{-1}\left(G(C_A)+C_B(T)\right)=C(T), 
$$
uniformly in $\alpha$, where $T<T_{M}$, with $T_{M}>0$ such that
$$
1+C_A\leq \dfrac{1}{\sqrt{2C_B(T_{M})}}. 
$$
Notice that $T_{M}>0$ does not depend on $\alpha$ thanks to \eqref{C_A}, and, since $C_B(T)\to 0$ as $T\to0$, there is some $T_M>0$ satisfying the inequality. Coming back to \eqref{3db}, have thus proved that
\begin{align}
\nonumber&\sqrt{\alpha}\Vert \partial_t\mathbf{v}_\alpha\Vert_{L^2(0,T;\mathbf{V}_\sigma)}+\sqrt{\alpha}\Vert \mathbf{v}_\alpha\Vert_{L^\infty(0,T;\mathbf{W}_\sigma)}+\Vert \mathbf{v}_\alpha\Vert_{L^\infty(0,T;\mathbf{V}_\sigma)}\\&+\Vert \mathbf{v}_\alpha\Vert_{L^2(0,T;\mathbf{W}_\sigma)}+\Vert \partial_t\mathbf{v}_\alpha\Vert_{L^2(0,T;\mathbf{H}_\sigma)}\leq C(T),\quad \forall T<T_{M}.
\label{w}
\end{align}

In the 2D case we can instead apply the classical Gronwall Lemma, to deduce, thanks again to \eqref{reg1}-\eqref{reg2} and \eqref{C_A}, that
\begin{align}
&\sqrt{\alpha}\Vert \partial_t\mathbf{v}_\alpha\Vert_{L^2(0,T;\mathbf{V}_\sigma)}\leq C(T),\\&
\sqrt{\alpha}\Vert \mathbf{v}_\alpha\Vert_{L^\infty_{uloc}([0,\infty);\mathbf{W}_\sigma)}+\Vert \mathbf{v}_\alpha\Vert_{L^\infty_{uloc}([0,\infty);\mathbf{V}_\sigma)}+\Vert \mathbf{v}_\alpha\Vert_{L^2_{loc}([0,\infty);\mathbf{W}_\sigma)}+\Vert \partial_t\mathbf{v}_\alpha\Vert_{L^2_{loc}([0,\infty);\mathbf{H}_\sigma)}\leq C,
\label{w2}
\end{align}
uniformly in $\alpha$.

We can now easily pass to the limit as $\alpha\to 0$, as in the proof of Theorem \ref{weakk}, recalling that, additionally, we have
$$
\alpha\mathbf{v}_\alpha\to 0 \quad \text{ in }L^\infty(I;\mathbf{W}_\sigma), \quad\text{and }\quad  \alpha\partial_t\mathbf{v}_\alpha\to 0\quad \text{ in }L^2(I;\mathbf{V}_\sigma),
$$
where $I=[0,T]$, for any $T<T_{M}$ in 3D, for any $T>0$ in 2D. Notice that in this context the strong convergence in a suitable topology for $\mathbf{v}_\alpha$ is much easier, since the sequence $\{\mathbf{v}_\alpha\}_{\alpha}$ is precompact in $BC(I;\mathbf{H}_\sigma)$ by Aubin-Lions Lemma. We have thus shown the existence of a triple $(\mathbf{v},\pphi,\ww)$ which is a local strong solution in 3D and a global strong solution in 2D according to Definition \ref{strong}. Note that, in the 3D case, the time $T_M$, depending only on suitable initial data norms, is a lower bound for the maximal time of existence of the obtained strong solution, $T_{max}$, so that $T_M\leq T_{max}$.  This concludes the proof of existence of strong solutions.
 \section{Proof of Theorem \ref{convective}}\label{sec:CH}
 Concerning uniqueness, one can easily adapt the proof of \cite[Sec.3]{EL}. In particular, we introduce the operator $\mathbb{L}: \mathbf{V}_0\to \mathbf{V}_0'$ as
 $$
 \langle\mathbb{L}\mathbf{u},\mathbf{v}\rangle=(\mathbf{M}\nabla\mathbf{u},\nabla\mathbf{v})_\Omega,
 $$ 
 which is invertible by the Lax-Milgram Lemma. Then considering two initial data $\pphi_{0,1}$ and $\pphi_{0,2}$, with the regularity stated in Theorem \ref{convective} and $\overline{\pphi}_{0,1}=\overline{\pphi}_{0,2}$, we define two corresponding strong solutions $\pphi_1$ and $\pphi_2$. By testing the equation satisfied by $\pphi:=\pphi_1-\pphi_2\in \mathbf{V}_0\hookrightarrow \mathbf{V}_0'$ with $\mathbb{L}^{-1}(\pphi_1-\pphi_2)$ and following \cite[Sec.3]{EL}, we get
 $$
 \frac 1 2\frac{d}{dt}\Vert \pphi\Vert^2_{\mathbf{V}_0'}+\Vert \nabla\pphi\Vert^2+\int_\Omega(\nabla\pphi)\textit{\textbf{v}}\cdot\mathbb{L}^{-1}\pphi dx+(\PsipA(\pphi_1)-\PsipA(\pphi_2),\pphi)_\Omega=(\mathbf{A}\pphi,\pphi)_\Omega
 $$
 where we considered an equivalent norm on $\mathbf{V}_0'$. Now, it is easy to see by monotonicity of $\PsipA$ that 
 $$
 (\PsipA(\pphi_1)-\PsipA(\pphi_2),\pphi)_\Omega\geq 0,
 $$
 and
 $$
 (\mathbf{A}\pphi,\pphi)_\Omega=(\mathbf{M}\nabla\mathbb{L}^{-1}\pphi, \nabla(\mathbf{A}\pphi))_\Omega\leq C\Vert \pphi\Vert_{\mathbf{V}_0'}\Vert \nabla \pphi\Vert\leq \frac{1}{2}\Vert{\nabla\pphi}\Vert^2+C\Vert \pphi\Vert_{\mathbf{V}_0'}^2.
 $$
 In conclusion, by standard inequalities and embeddings,
 $$
\left\vert \int_\Omega(\nabla\pphi)\textit{\textbf{v}}\cdot\mathbb{L}^{-1}\pphi dx\right\vert = \left\vert \int_\Omega(\nabla\mathbb{L}^{-1}\pphi)\textit{\textbf{v}}\cdot\pphi dx\right\vert\leq \Vert \pphi\Vert_{\mathbf{V}_0'}\Vert \textit{\textbf{v}}\Vert_{\mathbf{L}^3(\Omega)}\Vert \pphi\Vert_{\mathbf{L}^6(\Omega)}\leq C\Vert \pphi\Vert_{\mathbf{V}_0'}^2\Vert \textit{\textbf{v}}\Vert_{\mathbf{L}^3(\Omega)}^2+\frac{1}{4}\Vert \nabla\pphi\Vert^2,
 $$
 so that we end up with 
  $$
 \frac 1 2\frac{d}{dt}\Vert \pphi\Vert^2_{\mathbf{V}_0'}+\frac 1 4\Vert \nabla\pphi\Vert^2\leq C\Vert \pphi\Vert_{\mathbf{V}_0'}^2(1+\Vert \textit{\textbf{v}}\Vert_{\mathbf{L}^3(\Omega)}^2).
 $$
 Since $\textit{\textbf{v}}\in L^\infty(0,\infty;\mathbf{H}_\sigma)\cap L^2(0,\infty;\mathbf{V}_\sigma)\hookrightarrow L^2(0,\infty;\mathbf{L}^3(\Omega))$, we deduce the uniqueness of solutions by Gronwall's Lemma.
 
 The proof of the existence part of the theorem can be performed with the same time discretization scheme as in Section \ref{timediscr}, by considering only \eqref{phi0}$_{2}$-\eqref{phi0}$_{3}$ and \eqref{phi0}$_5$-\eqref{phi0}$_6$. First we need to discretize the prescribed velocity $\textit{\textbf{v}}$: in order to do so, we simply set, for $h:=\frac 1 M$, $M\in\N$ and $t_k:=kh$, $M\in\N_0$,
 $$
 \textit{\textbf{v}}^{k+1}:=\frac 1 h\int_{t_k}^{t_{k+1}}\textit{\textbf{v}}(s)ds\in \mathbf{V}_\sigma.
 $$
 Setting for simplicity $(\textit{\textbf{v}},\pphi,\ww):=(\textit{\textbf{v}}^{k+1},\pphi^{k+1},\ww^{k+1})$, we use the scheme: given the couple $(\textit{\textbf{v}},\pphi^{k})$, $k\in \N_0$, with $\pphi^k\in \mathcal{Z}_1$, find $(\pphi,\ww)$ such that $\pphi\in\mathcal{Z}_2$ and $\ww\in \mathbf{H}_\mathbf{n}^2$, satisfying
  \begin{align}
 \nonumber\dfrac{\pphi-\pphi^k}{h}+(\nabla\pphi^k)\textit{\textbf{v}}&=\operatorname{div}(\mathbf{M}\nabla\ww)&&\quad\text{ in }\Omega,\\
 \ww&= -\Delta\pphi+\mathbf{P}\left(\Psi_{,\pphi}^1(\pphi)-\frac  {\mathbf{A}\pphi}{2}-\frac{\mathbf{A}\pphi^k}{2}\right)&&\quad\text{ in }\Omega, \label{phi01}\\
 \nonumber\partial_\mathbf{n}\pphi&=\mathbf{0}&&\quad \text{ on }\partial\Omega,\\
 \partial_\mathbf{n}\ww &=\mathbf{0}&&\quad \text{ on }\partial\Omega.\nonumber
 \end{align} 
 In particular, we need to slightly adapt the proof of Theorem \ref{disc}, removing the equation for the velocity. We only highlight the main differences. First, it is immediate to see that we have the energy balance
 $$
 \mathcal{E}_{CH}(\pphi)+h\int_\Omega \mathbf{M}\nabla\ww :\nabla \ww dx+\frac 1 2\int_\Omega\vert \nabla\pphi-\nabla \pphi^k\vert^2dx+h\int_\Omega (\nabla\pphi^k)\textit{\textbf{v}}\cdot \ww dx\leq  \mathcal{E}_{CH}(\pphi^k),
 $$
 where we set
 $$
 \mathcal{E}_{CH}(\pphi):=\int_\Omega \Psi(\pphi)dx+\frac 1 2\int_\Omega\vert \nabla\pphi\vert^2dx.
 $$
Since then by assumption $\pphi_k\in[0,1]^N$ and $\textit{\textbf{v}}\in \mathbf{V}_\sigma$, we have
$$
\left\vert  \int_\Omega (\nabla\pphi^k)\textit{\textbf{v}}\cdot \ww dx\right\vert\leq \Vert \pphi^k\Vert_{\mathbf{L}^\infty(\Omega)}\Vert \textit{\textbf{v}}\Vert \Vert \nabla\ww\Vert\leq \frac 1 2 \int_\Omega \mathbf{M}\nabla\ww :\nabla \ww dx+C\Vert \textit{\textbf{v}}\Vert^2,
$$ 
 so that the energy balance implies 
  \begin{align}
 \mathcal{E}_{CH}(\pphi)+\frac h 2\int_\Omega \mathbf{M}\nabla\ww :\nabla \ww dx+\frac 1 2\int_\Omega\vert \nabla\pphi-\nabla \pphi^k\vert^2dx\leq  \mathcal{E}_{CH}(\pphi^k)+Ch\Vert \textit{\textbf{v}}\Vert^2.
 \label{enbal}
 \end{align} 
 We can now repeat almost word by word the fixed point argument adopted in the proof of Theorem \ref{disc} (removing the first row in the operators \eqref{L} and \eqref{F}, since now $\textit{\textbf{v}}$ is given). In particular, the only noticeable difference is in the study of \eqref{phia}, which is now given by
 \begin{align}
\nonumber -\operatorname{div}(\mathbf{M}(\pphi^k)\nabla\ww)+\int_\Omega\ww dx&=\lambda\left[-\frac{\pphi-\pphi^k}h-(\nabla\pphi^k)\textit{\textbf{v}}+\int_\Omega\ww dx\right],&&\quad\text{ in }\Omega,\\ \label{phiab}
 -\Delta\pphi+\mathbf{P}\PsipA(\pphi)&=\lambda\left[\ww+\mathbf{P}\mathbf{A}\frac{\pphi+\pphi^k}{2}\right],&&\quad\text{ in }\Omega,\\
 \partial_\mathbf{n}\pphi&=\mathbf{0},&&\quad\text{ on }\partial\Omega,\nonumber\\
 \partial_\mathbf{n}\ww&=\mathbf{0},&&\quad\text{ on }\partial\Omega,\nonumber
 \end{align}
 where now $\textit{\textbf{v}}\in \mathbf{V}_\sigma$ is given and fixed.
 Therefore, by a similar method as for obtaining the energy estimate \eqref{enbal}, we infer 
 \begin{align}
 \nonumber&\frac{1}{2}\int_\Omega\vert \nabla\pphi\vert^2-\lambda\int_\Omega\frac{\pphi^T\mathbf{A}\pphi}{2}+\int_\Omega\Psi^1(\pphi)dx +h\int_\Omega \mathbf{M}(\pphi^k)\nabla \ww:\nabla\ww dx\\&+(1-\lambda)h\left(\int_\Omega \ww dx\right)^2
 +\frac 1 2\int_\Omega\vert \nabla\pphi-\nabla \pphi^k\vert^2dx\nonumber\\&\leq \frac{1}{2}\int_\Omega\vert \nabla\pphi^k\vert^2-\lambda\int_\Omega\frac{{(\pphi^k)}^T\mathbf{A}\pphi^k}{2}+\int_\Omega\Psi^1(\pphi^k)dx-h\lambda\int_\Omega(\nabla\pphi^k)\textit{\textbf{v}}\cdot\ww dx.
 \label{ener22}
 \end{align}
As in the proof Theorem \ref{disc}, it holds $\pphi\in[0,1]^N$, but also $\left\vert \int_\Omega \Psi^1(\pphi)dx\right\vert\leq C$. Moreover, $\pphi^k\in[0,1]^N$ and thus $\left\vert \int_\Omega \Psi^1(\pphi^k)dx\right\vert\leq C$. Recall also that $\textit{\textbf{v}}\in \mathbf{V}_\sigma$. Therefore we have 
$$
\left\vert \int_\Omega(\nabla\pphi^k)\textit{\textbf{v}}\cdot\ww dx\right \vert\leq \Vert \pphi^k\Vert_{\mathbf{L}^\infty(\Omega)}\Vert \textit{\textbf{v}}\Vert \Vert \nabla\ww\Vert\leq \frac 1 2 \int_\Omega \mathbf{M}\nabla\ww :\nabla \ww dx+C_k, 
$$
so that, recalling $\lambda\in[0,1]$, we have the following
 \begin{align}
 & \frac{1}{2h}\int_\Omega\vert \nabla\pphi\vert^2dx+\frac 1 2\int_\Omega \mathbf{M}(\pphi^k)\nabla \ww:\nabla\ww dx+(1-\lambda)\left(\int_\Omega \ww dx\right)^2\leq \frac{C_k}{h},
 \label{ener2bis2}
 \end{align}
 where $h$ is fixed at this level. From this step on we can complete the proof exactly as in Theorem \ref{disc}, thus showing the existence of a solution to the discretized problem \eqref{phi01}. Concerning the uniform estimates, we have the energy balance \eqref{enbal} and we can repeat \textit{verbatim} the proof of Theorem \ref{disc} point (2), to obtain exactly the same result, i.e., that \eqref{h} holds also in this case.
 
 We now study how to pass to the limit as $M\to\infty$. Thanks to the validity of \eqref{enbal} and \eqref{h} we can follow the same arguments as in the proof of Theorem \ref{strong1}, omitting the equation for the velocity. In particular, after introducing the piecewise constant functions $\pphi^M,\ww^M,\textit{\textbf{v}}^M$ as in the aforementioned proof, we immediately obtain, as for \eqref{enerbis}, that 
  \begin{align}
 &\Vert \pphi^M\Vert_{L^\infty(0,\infty;\mathbf{H}^1(\Omega))}+ \frac{1}{\sqrt{h}}\Vert\pphi^M-\pphi^M_h\Vert_{L^2(0,\infty;\mathbf{H}^1(\Omega))}+\Vert \nabla\ww^M\Vert_{L^2(0,\infty;\mathbf{L}^2(\Omega))}\leq C.
 \label{unif2}
 \end{align}
 Moreover, by Theorem \ref{steaddy} point (3) applied to \eqref{ww}, with $\mathbf f=\ww^M+\mathbf{P}\left(\mathbf{A}\frac{\pphi^M+\pphi^M_h}{2}\right)$ and $\mathbf{m}=\overline{\pphi}_0$, we get, for any $T>0$,
 \begin{align}
 \Vert \pphi^M\Vert_{L^2(0,T;\mathbf{W}^{2,p}(\Omega))}+\Vert \PsipA(\pphi^M)\Vert_{L^2(0,T;\mathbf{L}^p(\Omega))}+\left\Vert\overline{\ww}^M\right\Vert_{L^2(0,T)}\leq C(T),
 \label{wwm3}
 \end{align}
 for $p\in[2,\infty)$ when $n=2$ and $p\in[2,6]$ if $n=3$.
 Defining $\widetilde{\pphi}^M(t):=\frac 1 h\int_{t-h}^t\pphi^M(s)ds$, which is the piecewise linear interpolant of $\pphi^M(kh)$, $k\in \N_0$, we have $\partial_{t}\widetilde{\pphi}^M=\partial_{t,h}^-\pphi^M$.
 Thus, it is easy to see from \eqref{dtp} and \eqref{unif} that 
 \begin{align}
 \Vert \partial_t\widetilde{\pphi}^M\Vert_{L^2(0,\infty;(\mathbf{H}^1(\Omega))')}\leq C,
 \label{dt11b}
 \end{align}
 uniformly in $M$. More precisely, we have, for any $k\in\N_0$, $t\in[kh,(k+1)h)$,
\begin{align}
\Vert \partial_t\widetilde{\pphi}^M(t)\Vert_{\mathbf{V}_0'}^2\leq C(\Vert \textit{\textbf{v}}^M(t)\Vert^2+ \left(\mathbf{M}\nabla \ww^M(t),\nabla\ww^M(t)\right)_\Omega).
\label{phit}
\end{align} 
  Note now that we have, by construction of $\textit{\textbf{v}}^M$, 
 \begin{align}
 \int_{0}^t\Vert \nabla\textit{\textbf{v}}^M(s)\Vert^2ds=\sum_{m=0}^{k-1 }h\Vert\nabla\textit{\textbf{v}}^{m+1} \Vert^2\leq \sum_{m=0}^{k-1}\int_{mh}^{(m+1)h}\Vert\nabla \textit{\textbf{v}}(s) \Vert^2ds\leq \int_0^\infty \Vert\nabla \textit{\textbf{v}}(s) \Vert^2ds\leq C,
 \label{ot}
 \end{align}
 for any $t=kh$, $k\in\N$, by assumption on $\textit{\textbf{v}}$. Furthermore, we have 
 \begin{align}
 \sup_{t\in[0,Lh)}\Vert \textit{\textbf{v}}^M(t)\Vert=\sup_{k=0,\ldots,L-1}\Vert \textit{\textbf{v}}^{k+1}\Vert=\sup_{k=0,\ldots,L-1}\frac 1 h\left\Vert \int_{kh}^{(k+1)h}\textit{\textbf{v}}(s)ds\right\Vert
\leq \Vert \textit{\textbf{v}}\Vert_{L^\infty(0,\infty;\mathbf{L}^2(\Omega))}\leq C, \label{ca}
 \end{align}
 for any $L\in\N$, by assumption. By summing \eqref{h} over $k=0,\ldots,L$, for any fixed $L>0$ we get (see also \eqref{d1})
 \begin{align}
 &\nonumber\frac 1 {2} \left(\mathbf{M}\nabla \ww^M((L-1)h),\nabla\ww^M((L-1)h)\right)_\Omega+\frac 1 2\int_{0}^{Lh} \left\Vert\nabla\partial_t\widetilde{\pphi}^M\right\Vert^2ds\\&\leq \frac 1 {2} \left(\mathbf{M}\nabla \ww^M_0,\nabla\ww^M_0\right)_\Omega+C(\min_{i=1,\ldots,N}\overline{\varphi}^0_i)\int_{0}^{Lh} \Vert \textit{\textbf{v}}^M\Vert\Vert \nabla\textit{\textbf{v}}^M\Vert\Vert \nabla\ww^M\Vert^2ds\nonumber\\&+C(\min_{i=1,\ldots,N}\overline{\varphi}^0_i)\int_{-h}^{(L-1)h}\Vert \nabla\textit{\textbf{v}}^M(h+s)\Vert^2\Vert \nabla\ww^M(s)\Vert^2ds+C\int_{0}^{Lh}\left\Vert \partial_t\widetilde{\pphi}^M\right\Vert^{2}_{\mathbf{V}_0'}ds,
 \label{d1b}
 \end{align}
 where in the last step we performed a change of variables. We can then exploit \eqref{phit} to get  
  \begin{align}
 &\nonumber\frac 1 {2} \left(\mathbf{M}\nabla \ww^M((L-1)h),\nabla\ww^M((L-1)h)\right)_\Omega+\frac 1 2\int_{0}^{Lh} \left\Vert\nabla\partial_t\widetilde{\pphi}^M\right\Vert^2ds\\&\leq \frac 1 {2} \left(\mathbf{M}\nabla \ww^M_0,\nabla\ww^M_0\right)_\Omega+C(\min_{i=1,\ldots,N}\overline{\varphi}^0_i)\int_{0}^{Lh} \Vert \textit{\textbf{v}}^M\Vert\Vert \nabla\textit{\textbf{v}}^M\Vert\Vert \nabla\ww^M\Vert^2ds\nonumber\\&+C(\min_{i=1,\ldots,N}\overline{\varphi}^0_i)\int_{-h}^{(L-1)h}\Vert \nabla\textit{\textbf{v}}^M(h+s)\Vert^2\Vert \nabla\ww^M(s)\Vert^2ds\nonumber\\&+C\int_{0}^{Lh}\left(\Vert \textit{\textbf{v}}^M\Vert^2+ \left(\mathbf{M}\nabla \ww^M,\nabla\ww^M\right)_\Omega\right)ds\nonumber\\&= \frac 1 {2} \left(\mathbf{M}\nabla \ww^M_0,\nabla\ww^M_0\right)_\Omega+C_0(\min_{i=1,\ldots,N}\overline{\varphi}^0_i)\int_{0}^{Lh} \Vert \textit{\textbf{v}}^M\Vert\Vert \nabla\textit{\textbf{v}}^M\Vert\Vert \nabla\ww^M\Vert^2ds\nonumber\\&+C_1(\min_{i=1,\ldots,N}\overline{\varphi}^0_i)\left(\int_{0}^{(L-1)h}\Vert \nabla\textit{\textbf{v}}^M(h+s)\Vert^2\Vert \nabla\ww^M(s)\Vert^2ds+\Vert \nabla\ww_0^M\Vert^2\int_0^h\Vert \nabla\textit{\textbf{v}}^M(s)\Vert^2ds\right)\nonumber\\&+C_2\int_{0}^{Lh}\left(\Vert \textit{\textbf{v}}^M\Vert^2+ \left(\mathbf{M}\nabla \ww^M,\nabla\ww^M\right)_\Omega\right)ds,
 \label{d1b1}
 \end{align}
 for some $C_i>0$, $i=0,1,2$.
 By the definition of $\ww^M$, Cauchy-Schwartz inequality and \eqref{ot}-\eqref{ca}, 
 \begin{align*}
 &\int_{(M-1)h}^{Mh}\Vert \textit{\textbf{v}}^M\Vert\Vert \nabla\textit{\textbf{v}}^M\Vert\Vert \nabla\ww^M\Vert^2dsds\\&
 =\Vert \nabla\ww^M((M-1)h)\Vert^2\Vert \textit{\textbf{v}}^M((M-1)h)\Vert\int_{(M-1)h}^{Mh}\Vert \nabla\textit{\textbf{v}}^M\Vert ds
 \\&\leq C\sqrt{h}\Vert \nabla\ww^M((M-1)h)\Vert^2\left(\int_{(M-1)h}^{Mh}\Vert \nabla\textit{\textbf{v}}^M\Vert^2 ds\right)^\frac 1 2\\&\leq \sqrt{h}C_3\left(\mathbf{M}\nabla \ww^M((M-1)h),\nabla\ww^M(
 (M-1)h)\right)_\Omega,
 \end{align*}
 for $C_3>0$ uniform in $M,L$, where we exploited property \eqref{nondeg} of $\mathbf{M}$.

  Recalling now that $\ww^M(t)=\ww^M((L-1)h)$ for any $t\in[(L-1)h,Lh)$, we deduce, by Poincar\'{e}'s inequality and \eqref{ot},
 \begin{align}
 &\nonumber \left(\frac 1 {2}-\sqrt{h}C_3\right)\left(\mathbf{M}\nabla \ww^M(t),\nabla\ww^M(t)\right)_\Omega+\frac 1 2\int_{0}^{t} \left\Vert\nabla\partial_t\widetilde{\pphi}^M\right\Vert^2ds\\&\leq \frac 1 {2} \left(\mathbf{M}\nabla \ww^M_0,\nabla\ww^M_0\right)_\Omega+C_0(\min_{i=1,\ldots,N}\overline{\varphi}^0_i)\int_{-h}^{t} \Vert \textit{\textbf{v}}^M\Vert\Vert \nabla\textit{\textbf{v}}^M\Vert\left(\mathbf{M}\nabla \ww^M,\nabla\ww^M\right)_\Omega ds\nonumber\\&+C(\min_{i=1,\ldots,N}\overline{\varphi}^0_i)\left(\int_{0}^{t}\Vert \nabla\textit{\textbf{v}}^M(h+s)\Vert^2\left(\mathbf{M}\nabla \ww^M(s),\nabla\ww^M(s)\right)_\Omega ds+\Vert \nabla\ww_0^M\Vert^2\int_0^h\Vert \nabla\textit{\textbf{v}}(s)\Vert^2ds\right)\nonumber\\&+C\int_{0}^{\infty}\Vert \nabla\textit{\textbf{v}}\Vert^2 ds+ C\int_{0}^{T}\left(\mathbf{M}\nabla \ww^M,\nabla\ww^M\right)_\Omega ds,
 \label{uloc1}
 \end{align}
 for any $t\in[(L-1)h,Lh)$ and for any $L>0$, $T>t$, implying that this inequality holds for any $0\leq t\leq T$ and any $T>0$. We can now choose, e.g., $h<\overline{h}$, with $\overline{h}:=	\frac 1 {16C_3^2}$, to apply Gronwall's Lemma to get 
 \begin{align}
 \nonumber&\frac{1}{4}\left(\mathbf{M}\nabla \ww^M(t),\nabla\ww^M(t)\right)_\Omega+\int_{0}^{t} \left\Vert\nabla\partial_t\widetilde{\pphi}^M\right\Vert^2ds\\&\nonumber\leq C\left(\left(\mathbf{M}\nabla \ww^M_0,\nabla\ww^M_0\right)_\Omega+\Vert \nabla\ww_0^M\Vert^2\int_0^h\Vert \nabla\textit{\textbf{v}}(s)\Vert^2ds+\int_{0}^{T}\left(\mathbf{M}\nabla \ww^M,\nabla\ww^M\right)_\Omega ds+\int_{0}^{\infty}\Vert \nabla\textit{\textbf{v}}\Vert^2 ds\right)\\&\nonumber\times e^{C\int_{0}^t\left(\Vert \textit{\textbf{v}}^M\Vert\Vert \nabla\textit{\textbf{v}}^M\Vert+\Vert \nabla\textit{\textbf{v}}^M\Vert^2\right)ds}\\&\nonumber\leq
 C\left(\left(\mathbf{M}\nabla \ww^M_0,\nabla\ww^M_0\right)_\Omega+\Vert \nabla\ww_0^M\Vert^2\int_0^h\Vert \nabla\textit{\textbf{v}}(s)\Vert^2ds+\int_{0}^{T}\left(\mathbf{M}\nabla \ww^M,\nabla\ww^M\right)_\Omega ds+\int_{0}^{\infty}\Vert \nabla\textit{\textbf{v}}\Vert^2 ds\right)\\&\times e^{C\int_{0}^\infty\Vert \nabla\textit{\textbf{v}}\Vert^2ds}, 
 \label{n1}
  \end{align}
 where in the last inequality we applied \eqref{ot} together with Poincar\'{e}'s inequality.
 %Clearly, from \eqref{d1b1} we end up with
 % \begin{align}
 % &\nonumber\left(\mathbf{M}\nabla \ww^M(t),\nabla\ww^M(t)\right)+\int_{0}^{t} \left\Vert\nabla\partial_t\widetilde{\pphi}^M\right\Vert^2ds\\&\leq
 % C\left(\left(\mathbf{M}\nabla \ww^M_0,\nabla\ww^M_0\right)+\Vert \nabla\ww_0^M\Vert^2\int_0^h\Vert \nabla\textit{\textbf{v}}(s)\Vert^2ds+\int_{0}^{T}\left(\mathbf{M}\nabla \ww^M,\nabla\ww^M\right)ds+\int_{0}^{\infty}\Vert \nabla\textit{\textbf{v}}\Vert^2 ds\right)e^{C\int_{0}^\infty\Vert \nabla\textit{\textbf{v}}\Vert^2ds}, 
 % \label{n1}
 % \end{align}
 % for any $0\leq t\leq T$ and any $T>0$. 
 Therefore, thanks to the assumptions on $\textit{\textbf{v}}$ and \eqref{unif2}, we get 
 \begin{align*}
 \Vert \nabla\ww^M\Vert_{L^\infty(0,\infty;\mathbf{L}^2(\Omega))}+\Vert \nabla\partial_t\widetilde{\pphi}^M\Vert_{L^2_{uloc}([0,\infty);\mathbf{L}^2(\Omega)}\leq C,
 \end{align*} 
 independently of $M$. A further application of Theorem \ref{steaddy} point (3), with the choice $\mathbf f=\ww^M+\mathbf{P}\left(\mathbf{A}\frac{\pphi^M+\pphi^M_h}{2}\right)$ and $\mathbf{m}=\overline{\pphi}_0$, then gives
 \begin{align}
 \Vert \pphi^M\Vert_{L^\infty(0,\infty;\mathbf{H}^2(\Omega))}+\Vert \PsipA(\pphi^M)\Vert_{L^\infty(0,\infty;\mathbf{L}^2(\Omega))}+\left\Vert\overline{\ww}^M\right\Vert_{L^\infty(0,\infty)}\leq C,
 \label{wwm2b}
 \end{align}	
 uniformly in $M$, so that we can obtain
 \begin{align}
 \Vert\ww^M\Vert_{L^\infty(0,\infty;\mathbf{H}^1(\Omega))}\leq C.
 \label{exta}
 \end{align}
  By comparison in \eqref{phiab}$_1$, it now easy to see that 
 \begin{align}
 \Vert \ww^M\Vert_{L^2_{uloc}([0,\infty);\mathbf{H}^3(\Omega))}\leq C,
 \label{regmu1}
 \end{align}
 uniformly in $M$. 
 By compactness arguments and arguing in a similar way as in the proof of Theorems \ref{weakk} and \ref{strong1}, we can pass to the limit in $M\to \infty$ and show that there exists a couple $(\pphi,\ww)$ which is a strong solution to \eqref{Ch}. In particular, this is possible since we have (see e.g.\ \cite[Lemma 4.11]{LT})
 $$
 \textit{\textbf{v}}^M\to \textit{\textbf{v}}\quad \text{ in }L^2(0,\infty;\mathbf{L}^2(\Omega)),
 $$
 for any $T>0$.
The proof of the existence is thus concluded. In order to get the estimate \eqref{w1}, we need to pass to the limit in \eqref{n1} as $M\to \infty$. To do so, we need that, up to subsequences,
\begin{align}
\int_{0}^{T}\left(\mathbf{M}\nabla \ww^M,\nabla\ww^M\right)_\Omega ds\to \int_{0}^{T}\left(\mathbf{M}\nabla \ww,\nabla\ww\right)_\Omega ds,
\label{co}
\end{align}
for almost any $T>0$. We thus exploit the energy identity. First notice that, by the above convergences, for almost any $t\in(0,T)$ and any $T>0$,
$$
\int_0^t\int_\Omega (\nabla\ww^M)\textit{\textbf{v}}^M\cdot \pphi^M dx ds\to \int_0^t\int_\Omega (\nabla\ww)\textit{\textbf{v}}\cdot \pphi dx ds.
$$
Moreover, it is easy do deduce, by weak lower semicontinuity of the norms involved and by Fatou's Lemma, that 
$$
\mathcal{E}_{CH}(\pphi(t)) \leq \liminf_{M\to\infty}\mathcal{E}_{CH}(\pphi^M(t)),
$$
for almost any $t\in(0,T)$ and any $T>0$. Therefore we have the following
\begin{align*}
&\mathcal{E}_{CH}(\pphi(t))+\liminf_{M\to\infty}\int_0^t\left(\mathbf{M}\nabla \ww^M,\nabla\ww^M\right)_\Omega ds\\&\leq\liminf_{M\to\infty}\left(\mathcal{E}_{CH}(\pphi^M(t))+\int_0^t\left(\mathbf{M}\nabla \ww^M,\nabla\ww^M\right)_\Omega ds\right)\\&\leq \liminf_{M\to\infty}\left(\mathcal{E}_{CH}(\pphi^M_0)+\int_0^t\int_\Omega (\nabla\ww^M)\textit{\textbf{v}}^M\cdot \pphi^M dx ds\right)\\&= \mathcal{E}_{CH}(\pphi_0)+\int_0^t\int_\Omega (\nabla\ww)\textit{\textbf{v}}\cdot \pphi dx ds,
\end{align*} 
for almost any $t\in(0,T)$ and any $T>0$. But we also know that the following energy identity holds
$$
\mathcal{E}_{CH}(\pphi(t))+\int_0^t\left(\mathbf{M}\nabla \ww,\nabla\ww\right)_\Omega ds= \mathcal{E}_{CH}(\pphi_0)+\int_0^t\int_\Omega (\nabla\ww)\textit{\textbf{v}}\cdot \pphi dx ds,
$$
so that by the chain of inequalities above we infer
$$ \int_0^t\left(\mathbf{M}\nabla \ww,\nabla\ww\right)_\Omega ds\geq \liminf_{M\to\infty}\int_0^t\left(\mathbf{M}\nabla \ww^M,\nabla\ww^M\right)_\Omega ds,
$$
which implies, together with the fact that $\nabla\ww^M\to \nabla\ww$ weakly in $L^2(0,\infty;\mathbf{L}^2(\Omega))$, that
$$
\liminf_{M\to\infty}\int_0^t\left(\mathbf{M}\nabla \ww^M,\nabla\ww^M\right)_\Omega ds= \int_0^t\left(\mathbf{M}\nabla \ww,\nabla\ww\right)_\Omega ds,
$$
for almost any $t\in(0,T)$ and any $T>0$. From this we clearly deduce \eqref{co} up to a subsequence. Thanks to this result, together with the above convergences, by choosing $T>0$ such that \eqref{co} holds, we can pass to the limit as $M\to\infty$ (and thus $h\to 0$), entailing \eqref{w1}. Indeed, we immediately obtain that, for almost any $t\in(0,T)$ and almost any $T>0$,
 \begin{align*}
&\nonumber\left(\mathbf{M}\nabla \ww(t),\nabla\ww(t)\right)_\Omega +\int_{0}^{t} \left\Vert\nabla\partial_t\widetilde{\pphi}\right\Vert^2ds\\&
\leq
C\left(\left(\mathbf{M}\nabla \ww_0,\nabla\ww_0\right)_\Omega +\int_{0}^{\infty}\left(\mathbf{M}\nabla \ww,\nabla\ww\right)_\Omega ds+\int_{0}^{\infty}\Vert \nabla\textit{\textbf{v}}\Vert^2 ds\right)e^{C\int_{0}^\infty\Vert \nabla\textit{\textbf{v}}\Vert^2ds},
\end{align*}
% so that 
%  \begin{align*}
% &\nonumber\left(\mathbf{M}\nabla \ww(t),\nabla\ww(t)\right)+\int_{0}^{t} \left\Vert\nabla\partial_t\widetilde{\pphi}\right\Vert^2ds
% \leq
% C\left(\left(\mathbf{M}\nabla \ww_0,\nabla\ww_0\right)+\int_{0}^{\infty}\left(\mathbf{M}\nabla \ww,\nabla\ww\right)ds+\int_{0}^{\infty}\Vert \nabla\textit{\textbf{v}}\Vert^2 ds\right)e^{C\int_{0}^\infty\Vert \nabla\textit{\textbf{v}}\Vert^2ds},
% \end{align*}
for almost any $t\geq 0$. Therefore, we conclude
 \begin{align*}
&\nonumber\Vert\left(\mathbf{M}\nabla \ww,\nabla\ww\right)_\Omega\Vert_{L^\infty(0,\infty)}+\int_{0}^{\infty} \left\Vert\nabla\partial_t\widetilde{\pphi}\right\Vert^2ds\\&
\leq
C\left(\left(\mathbf{M}\nabla \ww_0,\nabla\ww_0\right)_\Omega+\int_{0}^{\infty}\left(\mathbf{M}\nabla \ww,\nabla\ww\right)_\Omega ds+\int_{0}^{\infty}\Vert \nabla\textit{\textbf{v}}\Vert^2 ds\right)e^{C\int_{0}^\infty\Vert \nabla\textit{\textbf{v}}\Vert^2ds},
\end{align*}
entailing \eqref{w1}.
Concerning estimates \eqref{w3}-\eqref{w4}, they can be retrieved by the same arguments as in \cite[Corollary 2.6]{AGGio}. Moreover, in case $\pphi$ is strictly separated, the additional regularity on $\pphi$, i.e., $\pphi\in L^\infty(0,\infty;\mathbf{H}^3(\Omega))\cap L^2_{uloc}((0,\infty];\mathbf{H}^4(\Omega))$ can be retrieved by elliptic regularity on \eqref{Ch}$_2$, exploiting the separation property and the regularity on $\ww$.  The proof is then concluded.
\section{Proof of Theorem \ref{long}}\label{sec:long}
\subsection{Longtime behavior of solutions and strict separation property}
\label{ll1}
To prove Theorem \ref{strong} part (1), it is enough to notice that, by the regularity of a weak solution according to Definition \ref{weaks}, it follows that, for any $\tau>0$, there exists $0<\widetilde{\tau}\leq \tau$ such that $\pphi(\widetilde{\tau})$ satisfies all the assumptions on the initial datum required in Theorem \ref{convective}. This means that there exists a strong solution $(\widetilde{\pphi},\widetilde{\ww})$, defined in the statement of the aforementioned theorem, defined on $[\widetilde{\tau},+\infty)$. Since the strong solutions to the advective Cahn-Hilliard equations are unique, this implies that, for $t\geq \widetilde{\tau}$, and thus for any $t\geq \tau$, the weak solution $(\pphi,\mathbf{u})$ instantaneously regularizes in the $\pphi$ variable, since $\pphi\equiv \widetilde{\pphi}$ almost everywhere.
For the other parts of Theorem \ref{long} we need to work a little more. We first have the following
\begin{lemma}
	Let $(\uu,\pphi)$ be a weak solution according to Definition \ref{weaks}. Then $\mathbf{v}(t)\to \mathbf{0}$ in $\mathbf{L}^2(\Omega)$ as $t\to \infty$.
	\label{convergence}
\end{lemma}
\begin{proof}
	The proof of this Lemma is exactly the same as the one of \cite[Lemma 3.2]{AGGio}, as long as we substitute $\phi$ with $\pphi$ and $\mu$ with $\ww$. Indeed, the definition of the energy $\mathcal{E}$ is formally the same. 
\end{proof}
For a given $\mathbf{m}\in \Sigma$,
such that ${m}_{i}\in (0,1),$ for any $i=1,\ldots ,N$, we set
\begin{equation*}
\mathcal{V}_{\mathbf{m}}:=\Big\{\mathbf{u}\in \mathbf{H}^{1}(\Omega ):\ 0\leq
\mathbf{u}(x)\leq 1, \quad \text{ for a.a. }x\in \Omega ,\quad \overline{%
	\mathbf{u}}=\mathbf{m},\quad \sum_{i=1}^{N}u_i = 1\Big\},
\end{equation*}%
endowed with the $\mathbf{H}^{1}$-topology. We now introduce the $\omega-$limit set of a weak solution $(\mathbf{u},\pphi)$ with initial data in $\mathbf{H}_\sigma\times \mathcal{V}_{\mathbf{m}}$:
\begin{equation*}
\omega (\mathbf{v},\pphi)=\{(\mathbf{v}',\pphi')\in \mathbf{H}_\sigma\times\mathbf{H}^{2r}(\Omega )\cap
\mathcal{V}_{\mathbf{m}}:\exists\,t_{n} \nearrow \infty \text{ s.t. }(\mathbf{%
	v}(t_{n}),\pphi(t_n))\rightarrow (\mathbf{v}',\pphi')\text{ in }\mathbf{H}_\sigma\times \mathbf{H}^{2r}(\Omega )\},
\end{equation*}%
where $r\in \lbrack \tfrac{1}{2},1)$. In particular, for later purposes, we
fix $r\in (\tfrac{n}{4},1)$.
Thanks to Theorem \ref{long}, part (1), we have that, for any $\tau>0$, $\pphi\in L^{\infty }(\tau ,\infty ;%
\mathbf{H}^{2}(\Omega ))$ for any $\tau >0$. Moreover, we know from Lemma \ref{convergence} that $\mathbf{v}\to 0$ in $\mathbf{L}^2(\Omega)$ as $t\to\infty$. This means that $\omega (\mathbf{v},\pphi)$ is  non-empty and compact in $\mathbf{H}_\sigma\times\mathbf{H}%
^{2r}(\Omega )$. Moreover, it is easy to show that
\begin{equation}
\lim_{t\rightarrow \infty }\operatorname{dist}_{\mathbf{H}_\sigma\times\mathbf{H}^{2r}(\Omega )}((\mathbf{v}(t),\pphi(t)),\omega (\mathbf{v},\pphi))=0.
\label{conv3}
\end{equation}%
Let us set

\begin{equation*}
Z=\{\mathbf{u}\in \mathbf{H}^{1}(\Omega ):\mathcal{E}(\mathbf u)<+\infty \},
\end{equation*}%
where $\mathcal{E}$ is defined as 
 \begin{align}
 \mathcal{E}(\pphi):=\frac 1 2 \int_\Omega \vert \nabla \pphi\vert^2 dx +\int_\Omega \Psi(\pphi)dx,
 \label{enerCH}
 \end{align}
and introduce the
notion of stationary point. Given
\begin{equation*}
\mathbf{f}_{1}\in \mathcal{G}:=\{\mathbf{v}\in \mathbf{L}^{\infty }(\Omega ):\
\mathbf{v}(x)\in T\Sigma \text{ for almost any }x\in \Omega \},
\end{equation*}%
we say that $\mathbf{z}\in \mathbf{H}^{2}(\Omega )\cap Z$ is a stationary point if it solves the boundary value problem
\begin{equation}
\begin{cases}
-\Delta \mathbf{z}+\mathbf{P}\Psi _{,\mathbf{z}}^{1}(\mathbf{z})=\mathbf f_{1}+%
\mathbf{P}\mathbf{A}\mathbf{z},\quad &\text{ a.e. in }\Omega , \\
\partial _{\mathbf{n}}\mathbf{z}=0,\quad &\text{ a.e. on }\partial \Omega , \\
\sum_{i=1}^{N}\mathbf{z}_{i}= 1, \quad &\text{ in } \Omega.
\end{cases}
\label{steady2b}
\end{equation}%
Let then $\mathcal{W}$ be the set of all the stationary points. Theorem \ref{steaddy}, applied with $\mathbf f:=\mathbf f_{1}+%
\mathbf{P}\mathbf{A}\mathbf{z}\in \mathbf{L}^\infty(\Omega)$, guarantees that $\mathbf{z}$ is \textit{strictly separated} from
the pure phases, i.e., there exists $0<{\delta }=\delta (\mathbf f)<\frac{1}{N}$
such that
\begin{equation}
\delta <\mathbf{z}(x)  \label{prop2b}
\end{equation}%
for any $x\in \overline{\Omega }$. Thus all the stationary points in $%
\mathcal{W}$ are strictly separated, but possibly \textit{not} uniformly. About uniform strict separation, we have the following
\begin{lemma}
	\label{stationary}
	Let $(\mathbf{v},\pphi)$ be a weak solution according to Definition \ref{weaks}. Then 
	$$
	\omega(\mathbf{v},\pphi)\subset \left\{(\mathbf{0},\pphi'): \pphi'\in \mathbf{H}^2(\Omega)\cap \mathcal{V}_{\mathbf{m}} \text{ is a solution to the steady-state equation %
	\eqref{steady2b}, with } \mathbf{f}_1=\overline{\mathbf{P}\Psi_{,\mathbf{u}}({%
			\mathbf{u}_\infty})}\right\},
	$$
where $\mathbf{m}=\overline{\pphi}_0$. Moreover, there exists $\delta>0$ so that
\begin{equation*}
\delta<\pphi'(x), \quad\forall x\in \overline{\Omega },
\end{equation*}
for any $(\mathbf{0},\pphi')\in 	\omega(\mathbf{v},\pphi) $. In conclusion, there exists $\mathcal{E}_\infty$ such that 
$$
\mathcal{E}(\mathbf{0},\pphi')=\mathcal{E}_{CH}(\pphi')=\mathcal{E}_\infty,\quad\forall(\mathbf{0},\pphi')\in 	\omega(\mathbf{v},\pphi),\text{ and }  \lim_{t\to\infty}\mathcal{E}(\mathbf{v}(t),\pphi(t))=\mathcal{E}_\infty.
$$
\end{lemma}
\begin{proof}
	The proof can be carried out as in \cite[Lemma 3.3]{AGGio} (see also \cite[Lemma 3.11]{GGPS}). 
	\end{proof}
\begin{remark}
	Thanks to the characterization of the $\omega$-limit set, it is now immediate to deduce, by the regularity on $\pphi$, that $\omega(\mathbf{v},\pphi)$ is connected in the topology of $\mathbf{H}_\sigma\times \mathbf{H}^{2r}(\Omega)$.
\end{remark}
Thanks to the choice of $r\in (\tfrac{n}{4},1)$, $\omega (\mathbf{v},\pphi)$ is
compact in $\mathbf{H}_\sigma\times \mathbf{L}^{\infty }(\Omega )$. Arguing as in \cite[Sec.3.3]{GGPS}, there exists an open set $U_\varepsilon$ in $\mathbf{L}^\infty(\Omega)$ such that, for almost any $x\in \Omega $, we have 
\begin{equation}
0<\delta -\varepsilon \leq \mathbf{z}(x)\leq 1-((N-1)\delta -\varepsilon
)<1,\quad \forall \mathbf{z}\in U_{\varepsilon },
\label{sea}
\end{equation}%
where $\delta>0$ is the one given in Lemma \ref{stationary}.
Furthermore, by \eqref{conv1} and the embedding $\mathbf{H}^{2r}(\Omega
)\hookrightarrow \mathbf{L}^{\infty }(\Omega )$ we deduce that there exists $%
t^{\ast }\geq 0$, depending on $(\mathbf{v},\pphi)$, such that $\pphi(t)\in
U_{\varepsilon }$ for any $t\geq t^{\ast }$, i.e., the strict separation property holds for any $t\geq t^*$, concluding the proof of Theorem \ref{long} point (2).
\subsection{A weak-strong uniqueness result}
We now show a weak-strong uniqueness result, which is the analog of \cite[Thm.\ 4.1]{AGGio}:
\begin{theorem}
	Assume that, given $T>0$ and the mobility matrix $\mathbf{M}$ constant, and  $(\mathbf{u}_0, \pphi_0)$ complying with the regularity assumptions of Definition \ref{strong}, there exists a weak solution $(\mathbf{v},\pphi,\ww)$ on $[0,T]$ according to Definition \ref{weaks}, departing from those data,  such that, additionally, the regularity of Theorem \ref{convective} holds for $(\pphi,\ww)$. Assume also that there exists a strong solution  $(\mathbf{V},\Phi,\mathbf{W},\Pi)$ on $[0,T]$ according to Definition \ref{strong}. Moreover, assume that there exists $\delta\in(0,\tfrac 1 N)$ such that 
	$$
	\pphi(x,t)\geq \delta,\quad \Phi(x,t)\geq \delta\quad\forall x\in\overline{\Omega}\times[0,T],
	$$
	and also that 
	$$
	\pphi,\Phi\in L^2(0,T;\mathbf{H}^4(\Omega)).
	$$
	Then, $\mathbf{v}=\mathbf{U}$ and $\pphi=\Phi$ on $[0,T].$\label{ws}
\end{theorem}
\begin{proof}
	The proof is actually very similar to the one of \cite[Thm.4.1]{AGGio}. In particular, the treatment of the equation for the velocity is exactly the same, as long as we substitute the scalar phase field with the vectorial one. We can then deduce 
	\begin{align}
	&\int_\Omega \frac 1 2 \rho(\pphi(t))\vert \mathbf{v}(t)-\mathbf{V}(t)\vert^2dx+\int_0^t \int_\Omega \nu(\pphi)\vert D(\mathbf{v}-\mathbf{V})\vert^2dxd\tau\nonumber\\&\leq
	\frac{\nu_*}{2}\int_0^t\Vert D(\mathbf{v}-\mathbf{V})\Vert^2d\tau+\frac 1 4 \int_0^t\Vert \Delta^2(\pphi-\Phi)\Vert^2d\tau\nonumber\\&
	+C\int_0^t\left(1+\Vert \mathbf{V}\Vert^2_{\mathbf{H}^2(\Omega)}+\Vert \partial_t\mathbf{V}\Vert^2\right)\left(\Vert \mathbf{v}-\mathbf{V}\Vert^2+\Vert \Delta(\pphi-\Phi)\Vert^2\right).
	\label{ene}
	\end{align}
	Concerning the Cahn-Hilliard part, we need to test the equation for $\pphi-\Phi$ by $\Delta^2(\pphi-\Phi)$. Since then, clearly, $\mathbf{P}\Delta^2(\pphi-\Phi)=\Delta^2(\pphi-\Phi)$, we obtain, exactly as in \cite[(4.27)]{AGGio}, the differential equation
	\begin{align*}
	&\frac 1 2 \ddt \Vert \Delta(\pphi-\Phi)\Vert^2+(\mathbf{M} \Delta^2(\pphi-\Phi),\Delta^2(\pphi-\Phi))\\&=\int_\Omega (\nabla(\pphi-\Phi))\mathbf{v}\cdot \Delta^2(\pphi-\Phi)dx+\int_\Omega(\nabla\Phi)(\mathbf{v}-\mathbf{V})\cdot \Delta^2(\pphi-\Phi)dx\\&
	+\int_\Omega\mathbf{M}\Delta(\Psi_{,\pphi}(\pphi)-\Psi_{,\pphi}(\Phi))\cdot \Delta^2(\pphi-\Phi)dx.
	\end{align*}
	We can thus follow word by word the proof of \cite[Thm.4.1]{AGGio} and conclude by Gronwall's Lemma that $\mathbf{v}=\mathbf{V}$ and $\pphi=\Phi$ on $[0,T]$.
\end{proof}
 \subsection{Large time regularity of the velocity}
 Again this result is identical to the one in \cite[Sec.5]{AGGio}. In particular, we need to exploit the asymptotic strict separation property given in point (1) of Theorem \ref{long} and the weak-strong uniqueness of Theorem \ref{ws}. Clearly the only differences are related to the fact that we need to substitute the scalar quantities $\phi,\mu$ with the vector quantities $\pphi,\ww$. The proof can thus be carried out identically, also exploiting estimates \eqref{w1}-\eqref{w4}. Therefore, we can consider point (3) of Theorem \ref{long} to be proven.
\subsection{Convergence to equilibrium}  
We are now left to prove point (4) of Theorem \ref{long}, i.e., the convergence of each weak solution to a single equilibrium. 
Due to Theorem \ref{long} point (2), the singularities of $\psi $ and its derivatives no
longer play any role in our analysis as we are only interested in the
behavior of the solution $\pphi(t)$, as $t\rightarrow \infty $. Thus we
can alter the function $\psi $ outside the interval $I_{\varepsilon
}=[\delta -\varepsilon ,1-((N-1)\delta -\varepsilon )]$ ($\varepsilon>0$ is given in Section \ref{ll1}) in such a way that
the extension $\widetilde{\psi }$ is of class $C^{3}(\mathbb{R}^{N})$ and
additionally $|\widetilde{\psi }^{(j)}(s)|$, $j=1,2,3$, are uniformly
bounded on $\mathbb{R}$. Correspondingly we define $\widetilde{\Psi }(%
\mathbf{s}):=\sum_{i=1}^{N}\widetilde{\psi }(s_{i})-\frac{1}{2}%
\mathbf{As}\cdot \mathbf{s}$. Observe that $\widetilde{\psi }%
_{I_{\varepsilon }}=\psi $ and $\psi $ is analytic in $I_{\varepsilon }$ by
assumption (\textbf{E2}). We then introduce the "reduced" energy
$\widetilde{\mathcal{E}}_{CH} :\mathbf{V}_{0}\rightarrow \mathbb{R}$
by setting
\begin{equation*}
\widetilde{\mathcal{E}}_{CH}(\mathbf{z}):=\frac{1}{2}\int_{\Omega }|\nabla
\mathbf{z}|^{2}dx+\int_{\Omega }\widetilde{\Psi }(\mathbf{z}+\mathbf{m})dx.
\end{equation*}%
Observe that $\widetilde{\mathcal{E}}_{CH}(\pphi%
_{0}-\mathbf{m})=\mathcal{E}(\pphi_{0})$ for all $%
\pphi_{0}\in \mathcal{V}_{\mathbf{m}}\cap U_{\varepsilon }$, thanks to
\eqref{sea} and to the definition of the extension $\widetilde{\psi }$.
The proof is then based on the validity of the  {\L}ojasievicz-Simon inequality, which was proven in \cite[Lemma 3.15]{GGPS}.
\begin{lemma}[{\L}ojasiewicz-Simon Inequality]
	Assume $(\mathbf{0},\pphi')\in \omega (\mathbf{v},\pphi)$.
	Then there exist $\theta \in (0,\tfrac{1}{2}]$, $C,\sigma >0$ such that
	\begin{equation*}
	\left\vert \widetilde{\mathcal{E}}_{CH}(\mathbf{u})-\widetilde{\mathcal{E}}_{CH}(%
	\pphi'-\mathbf{m})\right\vert ^{1-\theta
	}\leq C\Vert \widetilde{\mathcal{E}}_{CH}^{\prime }(\mathbf{u})\Vert _{\mathbf{V}%
		_{0}^{\prime }},
	\end{equation*}%
	whenever $\mathbf{u}\in \mathbf{V}_0$ is such that $\Vert \mathbf{u}-\pphi'+\mathbf{m}%
	\Vert _{\mathbf{V}_{0}}\leq \sigma $. \label{Loja}
\end{lemma}
 Now, if we consider the kinetic energy $$\mathcal{E}_{kin}(\mathbf{v},\pphi):=\frac1 2 \int_\Omega  \rho(\pphi)\vert \mathbf{v}\vert^2dx,$$ 
 we immediately see that 
 \begin{align}
 \label{en}
 \left\vert \mathcal{E}_{kin}(\mathbf{V},\Phi) -\mathcal{E}_{kin}(\mathbf{0},\pphi')\right\vert^{\frac 1 2}\leq \sqrt{\frac{\max_{i=1,\ldots,N}\widetilde{\rho_i}}{2}}\Vert \mathbf{V}\Vert:=C_\rho\Vert \mathbf{V}\Vert,
 \end{align}
 for any $(\mathbf{V},\Phi)\in \mathbf{H}_\sigma\times \mathbf{L}^\infty(\Omega)$, with $\Phi\in[0,1]^N$. 
 As already noticed, thanks to Lemma \ref{stationary}, $\widetilde{\mathcal{E}}_{CH}(\mathbf{u})={\mathcal{E}}_{CH}(\mathbf{u}+\mathbf{m})$ for any $\mathbf{u}\in \omega(\mathbf{v},\pphi)-\mathbf{m}$. 
  Since $\mathbf{H}^{2r}(\Omega )\hookrightarrow \hookrightarrow \mathbf{H}%
 ^{1}(\Omega )$, $\omega (\mathbf{v},\pphi)$ is compact in $\mathbf{H}_\sigma\times \mathbf{H}%
 ^{1}(\Omega )$, thus we can find a finite number $M_{1}$ of $\mathbf{H}%
 ^{1}(\Omega )$-open balls $B_{m}$, $m=1,\ldots ,M_{1}$ of radius $\sigma $, centered at $\pphi_m\in \mathcal{P}_2\omega (\mathbf{v},\pphi)$ ($\mathcal{P}_2$ is the projector on the second component of a two-components vector),
 such that
 \begin{equation*}
 \mathcal{P}_2\omega (\mathbf{v},\pphi)\subset \widetilde{U}:= \bigcup_{m=1}^{M_{1}}B_{m}.
 \end{equation*}%
Therefore, also exploiting \eqref{en} and Lemma \ref{stationary} (so that $\mathcal{E}_{CH}(\pphi')=\mathcal{E}_\infty$ for any $\pphi'\in \mathcal{P}_2\omega(\mathbf{v},\pphi)$), we immediately obtain that there exist uniform $C>0,\theta\in(0,\tfrac 1 2]$ such that, for any $0<R\leq C_\rho^{-1}$,
 	\begin{equation}
 \left\vert \mathcal{E}_{kin}(\mathbf{V},\Phi)+\widetilde{\mathcal{E}}_{CH}(\mathbf{U}-\mathbf{m})-\mathcal{E}_\infty\right\vert ^{1-\theta
 }\leq C(\Vert \widetilde{\mathcal{E}}_{CH}^{\prime }(\mathbf{U}-\mathbf{m})\Vert _{\mathbf{V}%
 	_{0}^{\prime }}+\Vert \mathbf{V}\Vert),
 \label{lla}
 \end{equation}%
 for any $\mathbf{U}\in \mathbf{V}_0+\mathbf{m}\cap \widetilde{U}$ and $\Vert \mathbf{V}\Vert\leq R$. From \eqref{conv1} and the embedding $\mathbf{H}%
 ^{2r}(\Omega )\hookrightarrow \hookrightarrow \mathbf{H}^{1}(\Omega )$, we
 deduce that there exists $\tilde{t}>0$ such that $(\mathbf{v}(t),\pphi(t))\in \mathbf{H}_\sigma\times \widetilde{U}$ for any $t\geq \widetilde{t}$. Recalling the
 definition of $U_{\varepsilon }$ given in Section \ref{ll1}, we have
 \begin{equation*}
 (\mathbf{v}(t),\pphi(t))\in \mathbf{H}_\sigma\times {U}_\varepsilon,\quad \forall\, t\geq
 t^{\ast }.
 \end{equation*}%
 Moreover, since $\mathbf{v}(t)\to0$ in $\mathbf{H}_\sigma$, there exists $\overline{\overline{t}}>0$ such that $\Vert\mathbf{v}(t)\Vert\leq R$ for any $t\geq \overline{\overline{t}}$.
 therefore we can choose $\overline{t}:=\max \{\widetilde{t},t^{\ast },\overline{\overline{t}},t_R\}$ ($t_R$ defined in Theorem \ref{long} point (3)) and
 $\mathbf{U}=\widetilde{U}\cap U_{\varepsilon }$ such that $(\mathbf{v}(t),\pphi(t))\in \mathbf{B}_R\times {\mathbf{U}}$ for any $t\geq \overline{t}$ ($\mathbf{B}_R$ being the $\mathbf{H}_\sigma$-ball of radius $R$). Since then $\pphi(t)\in \mathbf{V}_{0}+\mathbf{m}\cap \widetilde{U}%
 $, it holds, recalling Poincar\'{e}'s inequality,
 	\begin{equation}
 \left\vert \mathcal{E}(\mathbf{v}(t),\pphi(t))-\mathcal{E}_\infty\right\vert ^{1-\theta
 }\leq C\left(\Vert \widetilde{\mathcal{E}}_{CH}^{\prime }(\pphi(t)-\mathbf{m})\Vert _{\mathbf{V}%
 	_{0}^{\prime }}+\left(\int_\Omega 2\nu(\pphi(t))\vert D\mathbf{v}(t)\vert^2dx\right)^{\frac 1 2}\right),\quad\forall t\geq \overline{t},
 \label{lla1}
 \end{equation}%
 since $\pphi(t)\in U_{\varepsilon }$ for any $t\geq
 \overline{t}$ and $\widetilde{\mathcal{E}}_{CH}(\pphi(t)-\mathbf{m}%
 )={\mathcal{E}}_{CH}(\pphi(t))$. Observe now that, for any $t\geq \overline{t}$%
 ,
 \begin{align*}
 & \langle \widetilde{\mathcal{E}}_{CH}^{\prime }(\pphi(t)-\mathbf{m}),\mathbf{h}\rangle _{\mathbf{V}_{0}^{\prime },\mathbf{V}_{0}}
 \\
 & =\int_{\Omega }\nabla \pphi(t)\cdot \nabla \mathbf{h}\,dx+\int_{\Omega }%
 \widetilde{\Psi }_{,\pphi}(\pphi(t))\cdot \mathbf{h}\,dx \\
 & =\int_{\Omega }\nabla \pphi(t)\cdot \nabla \mathbf{h}\,dx+\int_{\Omega }%
 {\Psi _{,\pphi}}(\pphi(t))\cdot \mathbf{h}\,dx \\
 & =\int_{\Omega }\left( -\Delta \pphi(t)+\mathbf{P}{{\Psi _{,\pphi}%
 }}(\pphi(t))\right) \cdot \mathbf{h}\,dx \\
 & =\int_{\Omega }\left( -\Delta \pphi(t)+P_{0}(\mathbf{P}{{\Psi _{\pphi}}}(\pphi(t)))\right) \cdot \mathbf{h}\,dx
 =(\ww(t)-\overline{\ww(t)},\mathbf{h}) \\
 & \leq \Vert \nabla \ww(t)\Vert \Vert \mathbf{h}\Vert \leq C\sqrt{(\mathbf{M}\nabla \ww(t),\nabla \ww(t))_\Omega}\Vert \mathbf{h}\Vert _{%
 	\mathbf{V}_{0}},\quad \forall\, \mathbf{h}\in \mathbf{V}_{0},
 \end{align*}%
 where $P_0$ is the $\mathbf{L}^2$-orthogonal projector on the space of functions with zero integral mean. Note that we exploited
 Poincar\'{e}'s inequality. This means that
 \begin{equation}
 \Vert \widetilde{\mathcal{E}}_{CH}^{\prime }(\pphi(t)-\mathbf{m}%
 )\Vert _{\mathbf{V}_{0}^{\prime }}\leq C\sqrt{(\mathbf{M} \nabla
 	\ww(t),\nabla {\ww}(t))_\Omega},\quad \forall\, t\geq \overline{t}.
 \label{ep1}
 \end{equation}%
 Recalling now that $\overline{t}\geq t_R$, and due to Theorem \ref{long} point (3), the energy identity holds: setting then $H(t):=\left\vert {%
 	\mathcal{E}}(\mathbf{v}(t),\pphi(t))-\mathcal{E}_{\infty }\right\vert ^{\theta }$, by \eqref{lla1} and \eqref{ep1},
 we have that
 \begin{align*}
  -\dfrac{d}{dt}H(t)&=-\theta \dfrac{d{\mathcal{E}}(\mathbf{u}(t))}{dt}%
 \left\vert {\mathcal{E}}(\mathbf{u}(t))-E_{\infty }\right\vert ^{\theta -1}
 \\
 & \geq \theta \dfrac{(\mathbf{M} \nabla {\ww}(t),\nabla {\ww}(t))_\Omega+\int_{\Omega}2\nu(\pphi(t))\vert D\mathbf{v}(t)\vert^2dx}{%
 	C\left(\Vert \widetilde{\mathcal{E}}^{\prime }(\mathbf{u}(t)-\mathbf{m}%
 	)\Vert _{\mathbf{V}_{0}^{\prime }}+\left(\int_\Omega 2\nu(\pphi(t))\vert D\mathbf{v}(t)\vert^2dx\right)^{\frac 1 2}\right)} \\
 & \geq C\sqrt{(\mathbf{M} \nabla {\ww}(t),\nabla {\ww}(t))_\Omega+\int_{\Omega}2\nu(\pphi(t))\vert D\mathbf{v}(t)\vert^2dx}, \quad \forall\, t\geq \overline{t},
 \end{align*}%
 where we exploited the inequality $(a+b)^{\frac 1 2}\leq a^\frac 1 2+b^\frac 1 2$, for $a,b\geq0$.
Since $H$ is a non nonincreasing nonnegative function such that $%
 H(t)\rightarrow 0$ as $t\rightarrow \infty $, we can integrate from $%
 \overline{t}$ to $+\infty $ and deduce that 
 \begin{equation*}
C\int_{\overline{t}}^{+\infty}\sqrt{(\mathbf{M} \nabla {\ww}(t),\nabla {\ww}(t))_\Omega+\int_{\Omega}2\nu(\pphi(t))\vert D\mathbf{v}(t)\vert^2dx}dt\leq CH(\overline{t})<+\infty ,
 \end{equation*}%
 i.e., $\nabla\ww\in L^{1}(\overline{t},+\infty ;\mathbf{L}%
 ^{2}(\Omega ))$ and $\mathbf{v}\in L^1(\overline{t},+\infty;\mathbf{V}_\sigma)$, entailing by comparison $\partial _{t}\pphi\in L^{1}(%
 \overline{t},+\infty ;\mathbf{H}^{1}(\Omega )^{\prime })$. Hence, there exists $(\mathbf{0},\pphi') \in \omega(\mathbf{v},\pphi)$  such that
 \begin{equation*}
 \pphi(t)=\pphi(\overline{t})+\int_{\overline{t}}^{t}\partial _{t}%
 \pphi(\tau )d\tau {\longrightarrow }%
 \ {\pphi}'\quad \text{ in }\mathbf{H}^{1}(\Omega )^{\prime },\quad \text{ as } t \to \infty,
 \end{equation*}%
 and, by uniqueness of the limit, we conclude that  $\omega(\mathbf{v},\pphi)=\{(\mathbf{0},{\pphi}')\}$.
 We also have $\lim_{t\rightarrow \infty }\pphi(t)={\pphi}'$
 in $\mathbf{H}^{2r}(\Omega )$ for a fixed $r\in (\tfrac{n}{4},1)$ (the one
 used in the definition of the $\omega $-limit set). On the other hand,
 thanks to the embedding $\mathbf{H}^{2r}(\Omega ) \hookrightarrow \mathbf{H}^{1}(\Omega )^{\prime }$,
 which is valid for all $r\in (0,1)$, we deduce that the convergence to the equilibrium actually holds for any $r\in (0,1)$. Moreover, since $(\mathbf{v},\pphi)\in L^\infty(t_R,+\infty;\mathbf{V}_\sigma\times \mathbf{W}^{2,p}(\Omega))$, by uniqueness of the limit we also deduce that $(\mathbf{v}(t),\pphi(t))\to (\mathbf{0},\pphi')$ weakly in $\mathbf{V}_\sigma\times \mathbf{W}^{2,p}(\Omega)$, for $p\in[2,6]$ when $n=3$ and $p\in[2,\infty)$ when $n=2$. The proof is thus finished.
 
\section{appendix}
\subsection{Bihari's Lemma}
\begin{lemma}(\cite{Bihari})
	\label{bihari}
Let $u$ and $f$ be non-negative continuous functions defined on $[0,T']$, where $T'\geq T$ and $T>0$ is such that 
$$
G(\alpha)+\int_0^t\,f(s)\,ds\in \mathfrak{D}(G^{-1}),\qquad \forall \, t \in [0,T],
$$
where $\alpha$ is a non-negative constant, the function $G$ is defined by
$$
G(x)=\int_{x_0}^x \frac{dy}{w(y)},\qquad \text{for } x \geq 0, 
$$
for a fixed $x_0>0$, and $G^{-1}$ denotes its inverse. Moreover, $w$ is a continuous non-decreasing function defined on $[0,\infty)$, with $w(u)>0$ on $(0,\infty)$.
If $u$ satisfies the following integral inequality
$$
u(t)\leq \alpha+ \int_0^t f(s)\,w(u(s))\,ds,\qquad \forall t\in[0,T'],
$$
then
$$
u(t)\leq G^{-1}\left(G(\alpha)+\int_0^t\,f(s) \, ds\right),\qquad \forall t\in[0,T].
$$
\end{lemma}
\textbf{Acknowledgments.} 
\RevB{We thank the anonymous referees for their valuable comments and remarks, which significantly improved the clarity of our work.}
Part of this work was done while AP was visiting HA and HG at the Department of Mathematics of the University of Regensburg, whose hospitality is kindly acknowledged. AP is a member of Gruppo Nazionale per l’Analisi Matematica, la Probabilitá e le loro Applicazioni (GNAMPA) of Istituto Nazionale per l’Alta Matematica (INdAM) and is supported by the MUR grant Dipartimento di Eccellenza 2023-2027. AP has also been partially funded by MIUR-PRIN research grant n. 2020F3NCPX.\\

\textbf{Conflict of interests.} On behalf of all authors, the corresponding author states that there is no conflict of interest.\\

\textbf{Data availability.} Data sharing not applicable to this article as no datasets were generated or analysed during the current study.
\bibliography{Bibliography}

\def\ocirc#1{\ifmmode\setbox0=\hbox{$#1$}\dimen0=\ht0 \advance\dimen0
  by1pt\rlap{\hbox to\wd0{\hss\raise\dimen0
  \hbox{\hskip.2em$\scriptscriptstyle\circ$}\hss}}#1\else {\accent"17 #1}\fi}
\begin{thebibliography}{10}

\bibitem{ADG}
H.~Abels, D.~Depner, and H.~Garcke.
\newblock Existence of weak solutions for a diffuse interface model for
  two-phase flows of incompressible fluids with different densities.
\newblock {\em J. Math. Fluid Mech.}, 15(3):453--480, 2013.

\bibitem{AGGio}
H.~Abels, H.~Garcke, and A.~Giorgini.
\newblock Global regularity and asymptotic stabilization for the incompressible
  {N}avier–{S}tokes-{C}ahn–{H}illiard model with unmatched densities.
\newblock {\em Math. Ann., \RevA{55pp., doi:10.1007/s00208-023-02670-2}}, 2023.

\bibitem{AGG}
H.~Abels, H.~Garcke, and G.~Gr{\"{u}}n.
\newblock Thermodynamically consistent, frame indifferent diffuse interface
  models for incompressible two-phase flows with different densities.
\newblock {\em Math. Models Methods Appl. Sci.}, 22(3):1150013 (40 pages),
  2012.

\bibitem{Banasnurnberg}
L.~Ba\v{n}as and R.~N\"{u}rnberg.
\newblock Numerical approximation of a non-smooth phase-field model for
  multicomponent incompressible flow.
\newblock {\em ESAIM Math. Model. Numer. Anal.}, 51(3):1089--1117, 2017.

\bibitem{Bihari}
I.~Bihari.
\newblock A generalization of a lemma of {B}ellman and its application to
  uniqueness problems of differential equations.
\newblock {\em Acta Math. Acad. Sci. Hungar.}, 7:81--94, 1956.

\bibitem{DongJCP2014}
S.~Dong.
\newblock An efficient algorithm for incompressible {N}-phase flows.
\newblock {\em J. Comput. Phys.}, 276:691--728, 2014.

\bibitem{DongJCP2018}
S.~Dong.
\newblock Multiphase flows of {$N$} immiscible incompressible fluids: a
  reduction-consistent and thermodynamically-consistent formulation and
  associated algorithm.
\newblock {\em J. Comput. Phys.}, 361:1--49, 2018.

\bibitem{DunbarLamStinner}
O.~R.~A. Dunbar, K.~F. Lam, and B.~Stinner.
\newblock Phase field modelling of surfactants in multi-phase flow.
\newblock {\em Interfaces Free Bound.}, 21(4):495--547, 2019.

\bibitem{EL}
C.~M. Elliott and S.~Luckhaus.
\newblock A generalized equation for phase separation of a multi-component
  mixture with interfacial free energy.
\newblock {\em IMA Preprint Series \# 887}, 1991.

\bibitem{GGPS}
C.~Gal, M.~Grasselli, A.~Poiatti, and J.~Shomberg.
\newblock Multi-component {C}ahn-{H}illiard systems with singular potentials:
  Theoretical results.
\newblock {\em \RevA{Appl. Math. Optim.}}, 88(73), 2023.

\bibitem{GGGP}
C.~G. Gal, A.~Giorgini, M.~Grasselli, and A.~Poiatti.
\newblock Global well-posedness and convergence to equilibrium for the
  abels-garcke-gr\"un model with nonlocal free energy.
\newblock {\em J. Math. Pures Appl. (9)}, \RevA{178}:\RevA{46--109}, 2023.

\bibitem{Galdi}
G.~P. Galdi.
\newblock {\em An Introduction to the Mathematical Theory of the
  {N}avier-{S}tokes Equations, Volume 1}.
\newblock Springer, Berlin - Heidelberg - New York, 1994.

\bibitem{GarckeElasticMisfit}
H.~Garcke.
\newblock On a {C}ahn-{H}illiard model for phase separation with elastic
  misfit.
\newblock {\em Ann. Inst. H. Poincar\'e Anal. Non Lin\'eaire}, 22(2):165--185,
  2005.

\bibitem{GGW}
A.~Giorgini, M.~Grasselli, and H.~Wu.
\newblock On the mass-conserving {A}llen-{C}ahn approximation for
  incompressible binary fluids.
\newblock {\em J. Funct. Anal.}, 283(9):Paper No. 109631, 2022.

\bibitem{GP}
M.~Grasselli and A.~Poiatti.
\newblock Multi-component conserved {A}llen-{C}ahn equations.
\newblock {\em Interfaces Free Bound., in press; preprint arXiv:2304.03363},
  2023.

\bibitem{GPV96}
M.~E. Gurtin, D.~Polignone, and J.~Vi{\~n}als.
\newblock Two-phase binary fluids and immiscible fluids described by an order
  parameter.
\newblock {\em Math. Models Methods Appl. Sci.}, 6(6):815--831, 1996.

\bibitem{HH77}
P.~Hohenberg and B.~Halperin.
\newblock Theory of dynamic critical phenomena.
\newblock {\em Rev. Mod. Phys.}, 49:435--479, 1977.

\bibitem{KimKangLowengrub}
J.~Kim, K.~Kang, and J.~Lowengrub.
\newblock Conservative multigrid methods for {C}ahn-{H}illiard fluids.
\newblock {\em J. Comput. Phys.}, 193(2):511--543, 2004.

\bibitem{KimLowengrub}
J.~Kim and J.~Lowengrub.
\newblock Phase field modeling and simulation of three-phase flows.
\newblock {\em Interfaces Free Bound.}, 7(4):435--466, 2005.

\bibitem{LT}
Y.~Liu and D.~Trautwein.
\newblock On a diffuse interface model for incompressible viscoelastic
  two-phase flows.
\newblock {\em preprint arXiv:2212.13507}, 2022.

\bibitem{LowengrubQuasiIncompressible}
J.~Lowengrub and L.~Truskinovsky.
\newblock Quasi-incompressible {C}ahn-{H}illiard fluids and topological
  transitions.
\newblock {\em R. Soc. Lond. Proc. Ser. A Math. Phys. Eng. Sci.},
  454(1978):2617--2654, 1998.

\bibitem{UnifiedNSCH}
M.~Ten~Eikelder, K.~van~der Zee, I.~Akkerman, and D.~Schillinger.
\newblock A unified framework for {N}avier-{S}tokes {C}ahn-{H}illiard models
  with non-matching densities.
\newblock {\em Math.\ Models Meth.\ Appl.\ Sci.}, 33(1):175--221, 2023.

\bibitem{YangDongJCP2018}
Z.~Yang and S.~Dong.
\newblock Multiphase flows of {$N$} immiscible incompressible fluids: an
  outflow/open boundary condition and algorithm.
\newblock {\em J. Comput. Phys.}, 366:33--70, 2018.

\end{thebibliography}
\bibliographystyle{abbrv}

\end{document}